\documentstyle{article}
\input amssymb.sty
\input epsf

\newtheorem{theorem}{Theorem}[section]
\newtheorem{lemma}[theorem]{Lemma}
\newtheorem{proposition}[theorem]{Proposition}
\newtheorem{corollary}[theorem]{Corollary}
\newtheorem{conjecture}[theorem]{Conjecture}
\newtheorem{definition}[theorem]{Definition}

\newtheorem{theorem-construction}[theorem]{Theorem--Construction}
\newtheorem{lemma-construction}[theorem]{Lemma--Construction}

\begin{document}
 
\newcommand{\Z}{{\Bbb Z}}
\newcommand{\R}{{\Bbb R}}
\newcommand{\Q}{{\Bbb Q}}
\newcommand{\C}{{\Bbb C}}
\newcommand{\lra}{\longrightarrow}
\newcommand{\lms}{\longmapsto}
\newcommand{\AAA}{{\Bbb A}}
\newcommand{\Alt}{{\rm Alt}}
\newcommand{\wg}{\wedge}
\newcommand{\ol}{\overline}
\newcommand{\CP}{{\Bbb C}P}
\newcommand{\bwg}{\bigwedge}
\newcommand{\caL}{{\cal L}}
\newcommand{\PP}{{\Bbb P}}
\newcommand{\HH}{{\Bbb H}}
\newcommand{\LL}{{\Bbb L}}

\begin{titlepage}
\title{Polylogarithms, regulators, and Arakelov motivic complexes}
\author{A. B.  Goncharov}
\date{}
\end{titlepage}
\stepcounter{page}
\maketitle
\tableofcontents

\section  {Introduction}      

{\bf Abstract}. We construct an explicit regulator map 
from the weight $n$ Bloch Higher Chow group  complex to the weight $n$ 
Deligne complex 
of a regular projective complex algebraic variety $X$.  
We define the weight $n$ 
Arakelov motivic complex
  as the cone of this map 
shifted by one.
 Its last cohomology group is (a version of) the 
Arakelov Chow group defined by H. Gillet and C. Soul\'e ([GS]). 
 
We 
relate the 
Grassmannian $n$--logarithms (defined as in [G5])  
to the geometry of the symmetric space  
$SL_n(\C)/SU(n)$. For  $n=2$  we recover  
Lobachevsky's formula expressing the volume of an ideal geodesic 
simplex in the hyperbolic space via the dilogarithm. 
Using the relationship with symmetric spaces 
we construct the Borel regulator 
on $K_{2n-1}(\C)$ via the 
Grassmannian $n$--logarithms. 

We study the Chow dilogarithm and prove a 
reciprocity law 
which strengthens Suslin's reciprocity law for Milnor's group $K^M_3$ on curves. 

Our note [G5] can serve as an introduction to this paper.  

{\bf 1. Beilinson's conjectures on special values of L-functions}. 
Let $X$ be a regular scheme. 
A.A. Beilinson [B1] defined the rational motivic cohomology 
of $X$ via Quillen's 
algebraic $K$-theory of $X$ by the following formula:
\begin{equation} \label{4.29.02.1}
H^i_{\cal M}(X, \Q(n)) = K^{(n)}_{2n-i}(X)_{\Q}
\end{equation}
 where on the right stays the weight $n$ eigenspace for the Adams operations 
acting on $K_{2n-i}(X)\otimes \Q$. 
For a regular complex algebraic variety $X$  Beilinson defined the regulator map 
to the weight $n$ Deligne cohomology of 
$X(\C)$: 
$$
r_B: H_{{\cal M}}^{i}(X, \Q(n)) \lra H_{{\cal D}}^{i}(X(\C), \R(n))
$$

Let $X$ be a regular projective scheme over 
$\Q$. Let  $L(h^i(X),s)$  be the  $L$-function related to its i-dimensional
cohomology. 
Beilinson  conjectured that for any integer $n > 1+i/2$ 
its  special value at $s=n$  is described, up to a nonzero rational factor,
  by  the regulator map to the weight $n$ {\it real} Deligne cohomology of 
$X$ 
$$
r_B: H_{{\cal M}}^{i+1}(X, \Q(n))_{\Z} \lra H_{{\cal D}}^{i+1}(X\otimes_{\Q}\R_{/\R}, \R(n))
$$
Here the subscript $\Z$ on the left indicates the subspace of $H_{{\cal M}}^{i+1}(X, \Q(n))$ 
coming from a regular model of the scheme $X$ over $\Z$, see [B1] and [RSS] for details.

This conjecture is fully established only when $X= {\rm Spec}(F)$ where
$F$ is a number field. In this case the regulator map $r_B$ coincides, up to a 
non-zero rational factor, 
 with 
 the Borel regulator ([B1]), 
and the relation with special values of the Dedekind zeta-function of $F$ was given
by the Borel theorem [Bo]. 

Although Beilinson's conjectures are far from being
proved, it is interesting to see what kind of information about the
special values of $L$-functions they 
suggest. So we come to the problem of {\it explicit} calculation of
Beilinson's regulator. This problem is already very interesting for the
Borel regulator.

{\bf 2. Regulator maps on motivic complexes and Arakelov motivic cohomology}.
Beilinson [B2] and S. Lichtenbaum [L1] conjectured that the weight $n$ 
integral 
motivic cohomology 
of $X$ should appear as  the cohomology of certain complexes of abelian groups 
$\Gamma(X;n)$, called the weight $n$ motivic complexes:
$$
H^i_{\cal M}(X, \Z(n)):= H^i(\Gamma(X;n))
$$
These  complexes are well defined as objects of the 
derived category. They 
should appear as the hypercohomology of certain complexes of Zarisky 
sheaves. 

The first  motivic complexes satisfying  
Beilinson's formula (\ref{4.29.02.1}) were Bloch's  Higher Chow group complexes 
${\cal Z}^{\bullet}(X; n)$ [Bl1]. Later on A.A. Suslin and V.A. Voevodsky 
defined several important versions of these complexes. 
For another candidates for motivic complexes, called the polylogarithmic motivic 
complexes,  see [G1-2]. They are very explicit and the smallest among 
possible candidates, however Beilinson's formula (\ref{4.29.02.1}) 
is far from being established for them.

The real Deligne cohomology arises also as the cohomology  
of certain complexes. 
It was suggested in [G5] and [G7] that the regulator map should be {\it explicitly} 
defined on the level of complexes. 

Let $X$ be a regular projective variety over $\C$. 
In Chapter 2 we construct 
a  homomorphism of complexes
\begin{equation} \label{4.29.02.2}
\mbox{Bloch's weight $n$ Higher Chow group complex ${\cal Z}^ {\bullet}(X; n)$ of $X$} 
\quad \lra 
\end{equation}
$$
\mbox{the weight $n$ real Deligne complex 
${\cal C}^{\bullet}_{\cal D}(X(\C); n)$ of $X$}
$$
This construction is a version of the one given in [G5]. 
The  complex ${\cal C}^{\bullet}_{\cal D}(X(\C); n)$ is the  
truncation $\tau_{\leq 2n}$ of the complex 
proposed by 
Deligne [Del]. The 
weight $n$ Arakelov motivic  complex $\Gamma_{\cal A}^{\bullet}(X; n)$ is the cone of 
the map (\ref{4.29.02.2}), shifted by $-1$:
\begin{equation} \label{6.11.02.14}
\Gamma_{\cal A}^{\bullet}(X; n):= {\rm Cone}\Bigl({\cal Z}^ {\bullet}(X; n)
\stackrel{(\ref{4.29.02.2})}{\lra} 
 {\cal C}^{\bullet}_{\cal D}(X(\C); n)\Bigr)[-1]
\end{equation}
For a regular projective variety   $X$ over $\R$ the image of  map 
(\ref{4.29.02.2}) lies in the subcomplex
$$
{\cal C}^{\bullet}_{\cal D}(X_{/\R}; n) := 
{\cal C}^{\bullet}_{\cal D}(X(\C); n)^{\overline F_{\infty}}
$$
where $\overline F_{\infty}$ is the De Rham involution provided by the action of 
complex conjugation. 
The weight $n$ 
{\it real} Arakelov motivic complex 
is defined as 
\begin{equation} \label{6.11.02.4}
\Gamma_{\cal A}^{\bullet}(X_{/\R}; n):= {\rm Cone}\Bigl({\cal Z}^ {\bullet}(X; n)
\stackrel{(\ref{4.29.02.2})}{\lra} 
 {\cal C}^{\bullet}_{\cal D}(X_{/\R}; n)\Bigr)[-1]
\end{equation}
Let $X$ be a regular projective variety  $X$ over a number field $F$. We view $X$ over $\Q$: 
$X \lra {\rm Spec}(F) \lra {\rm Spec}(\Q)$, and set
\begin{equation} \label{6.11.02.4q}
\Gamma_{\cal A}^{\bullet}(X_{/F}; n):= {\rm Cone}\Bigl({\cal Z}^ {\bullet}(X; n)
\stackrel{}{\lra} 
 {\cal C}^{\bullet}_{\cal D}(X\otimes_{\Q} \R_{/\R}; n)\Bigr)[-1]
\end{equation}
The weight $n$ Arakelov motivic cohomology is the cohomology of this complex. 
Our construction 
works equally well for the Suslin-Voevodsky versions of the 
motivic complexes. 

 Taking the cohomology we get a construction of the 
regulator map on motivic cohomology. 
For a different construction see [Bl3].

The regulator map on the polylogarithmic motivic complexes 
was defined in [G7] explicitly via the classical polylogarithms. 
The Arakelov motivic complexes constructed using  regulator maps 
on different motivic complexes are supposed to lead to 
the same object of the derived category.  
However a precise relationship between the construction 
given in [G7] and the one 
in Chapter 2 is not clear.

{\it Higher Arakelov Chow groups}. The last group of the 
complex ${\cal C}^{\bullet}_{\cal D}(X(\C); n)$ consists of 
closed distributions of a certain type on $X(\C)$. Replacing it 
by the quotient modulo smooth closed 
forms of the same type we get the quotient complex 
$\widetilde {\cal C}^{\bullet}_{\cal D}(X(\C); n)$. 
Changing  ${\cal C}$ to  
$\widetilde {\cal C}$ in (\ref{6.11.02.14}) 
we define the weight $n$ 
Higher Arakelov Chow group complex. Its last cohomology group is isomorphic to 
the Arakelov Chow group $\widehat {CH}^n(X(\C))$ as defined by Gillet and Soul\'e 
[GS], [S]. 

{\bf Problems}. a) Show that taking
 cohomology of
 the map (\ref{4.29.02.2}) and using the isomorphism between 
the rational Bloch's Higher Chow groups of $X$ and the corresponding part of the rational 
$K$-theory of $X$ ([Bl2], [Lev]) 
we get a non-zero rational multiple of the Beilinson's regulator map. 

b) To generalize the arithmetic Riemann--Roch theorem 
proved by Gillet and Soul\'e to the case of Higher Arakelov Chow groups. 
 
  {\bf Remark}. The weight $n$ Arakelov motivic complex should be considered as an ingrediant 
of a definition of the weight $n$ {\it arithmetic motivic complex}. The latter 
is related to the regulator maps on 
$H_{{\cal M}}^{\bullet}(X, \Q(n))_{\Z}$, while the former is related  
to the ones on $H_{{\cal M}}^{\bullet}(X, \Q(n))$. 
Ideally one should 
have for every place $p$ of $\Q$ a map from the left hand side of (\ref{4.29.02.2}) 
to a certain complex, which for the Archimedian  place should be given by our map. 
Then one should take the shifted by $-1$ cone of the sum of these maps. 
  
{\bf 3. The Chow $n$--logarithm function}. 
Let us  describe the  regulator map (\ref{4.29.02.2})
in the simplest case when $X = {\rm Spec}(\C)$ is a point. 

Let us choose in $\PP^m$  homogeneous coordinates $(z_0: ... : z_m)$. 
The union of the coordinate hyperplanes is a simplex $L$. 
Let $\AAA^m$ be the complement to
 the hyperplane $z_1 + ... + z_m = z_0$ in $\PP^m$. The abelian group 
${\cal Z}_m({\rm Spec}(\C); n)$ is freely generated by 
the codimension $n$ irreducible algebraic cycles in $\AAA^{m}$ intersecting  properly 
the faces of the simplex $L$. The intersection with codimension one faces $L_j$ of $L$ provide homomorphisms 
$$
\partial_j: {\cal Z}_m({\rm Spec}(\C); n) \lra {\cal Z}_{m-1}({\rm Spec}(\C); n); \quad 
\partial:= \sum_{j=0}^{m} (-1)^j\partial_j 
$$
The weight $n$ Higher Chow group complex over ${\rm Spec}(\C)$, where $n>0$, 
written as a homological complex,  looks as follows:
$$
... \stackrel{\partial}{\lra} {\cal Z}_2({\rm Spec}(\C); n) \stackrel{\partial}{\lra}
{\cal Z}_1({\rm Spec}(\C); n) \stackrel{\partial}{\lra} {\cal Z}_0({\rm Spec}(\C); n) 
$$

The Deligne complex of a point is the complex $(2\pi i)^{n}\R \lra \C$, with $\C$ is 
in the degree $+1$. So it is quasiisomorphic to the group $\R(n-1):= (2\pi i)^{n-1}\R$ placed in degree $+1$. 
The
 regulator map (\ref{4.29.02.2}) 
boils down to a construction  of a homomorphism 
$$
{\cal Z}_{2n-1}({\rm Spec}(\C); n)  \stackrel{{\cal P}_n}{\lra} \R(n-1), \quad 
\mbox{such that 
${\cal P}_n \circ \partial = 0$} 
$$
It is provided by  
a function 
${\cal P}_n$ on the space of 
codimension $n$ cycles in $\C\PP^{2n-1}$ intersecting properly 
 faces of a simplex $L$. This function, called 
the {\it Chow $n$--logarithm function}, was constructed in [G5]. 
To recall its construction, observe that 
a codimension $n$ cycle given by an irreducible 
 subvariety $X$ in $\PP^{2n-1} -L$ provides the 
$(n-1)$-dimensional variety $X$ with $2n-1$ rational functions 
$f_1, ..., f_{2n-1}$:  These functions are obtained by 
restriction of the coordinate functions 
$z_i/z_0$ to the cycle $X$. We define a natural $(2n-2)$-form 
$r_{2n-2}(f_1, ..., f_{2n-1})$ on $X(\C)$ and set
\begin{equation} \label{6.17.02.1}
{\cal P}_n(X; f_1, ..., f_{2n-1}):= (2\pi i)^{1-n}\int_{X(\C)}r_{2n-2}(f_1, ..., f_{2n-1})
\end{equation}

{\bf 4. An example: the Chow dilogarithm}. 
Let $f_1,f_2,f_3$ be three arbitrary 
rational functions on a complex curve $X$. 
Set 
$$
r_2(f_1,f_2,f_3):= 
$$
$$
{\rm Alt}_3 \Bigl(\frac{1}{6}\log|f_1| d \log|f_2| 
\wedge d
\log|f_3|
-\frac{1}{2}
\log|f_1| d\arg f_2 \wedge d\arg f_3 \Bigr) 
$$
 where 
${\rm Alt}_3$ is the alternation of $f_{1}, f_{2}, f_{3}$. 
 Consider the space of
quadruples $(X;f_1,f_2,f_3)$. It is a union of finite dimensional algebraic varieties. 
The Chow dilogarithm is a real function 
on its complex points  defined by the formula
$$
{\cal P}_2(X;f_1,f_2,f_3):= \frac{1}{2\pi i}\int_{X(\C)}r_2(f_1,f_2,f_3)
$$
The integral converges. 
The Chow dilogarithm provides a homomorphism 
\begin{equation} \label{chowh}
\Lambda^3 \C(X)^* \to \R,\quad  f_1\wedge f_2\wedge f_3 \lms {\cal P}_2(X;f_1,f_2,f_3)
\end{equation}

Why does the dilogarithm appear in the name of the function ${\cal P}_2$? 
Recall  the classical dilogarithm 
$$
Li_2(z):= - \int_0^z\log(1-z)d\log z
$$ 
It has a single-valued cousin, the Bloch-Wigner  function:
$$
{\cal L}_2(z) := {\rm Im} Li_2(z) + \arg (1-z) \log \vert z \vert 
$$

The Chow dilogarithm
is defined by a two-dimensional integral
 over $X(\C)$, while ${\cal L}_2(z)$ 
is given by an 
integral over a path in ${\C}{\Bbb P}^1$.  In Chapter 6 we
 show that nevertheless the Chow dilogarithm can be expressed by the function ${\cal L}_2(z)$. Here is how it 
works when $X = \C\PP^1$. 
For $f \in \C(X)$ let $v_x(f)$ be the order of zero of $f$ at $x \in X(\C)$. 
Choose a point $\infty$ on ${\Bbb P}^1$. Then 
\begin{equation} \label{6.14.02.2}
{\cal P}_2(\C\PP^1; f_1, f_2, f_3) =  \sum_{x_i \in {\Bbb P}^1(\C)} 
v_{x_1}(f_1)v_{x_2}(f_2)v_{x_3}(f_3){\cal L}_2(r(x_1, x_2, x_3, \infty))
\end{equation}where $r(...)$ denotes the cross-ratio of four points on $\PP^1$. 
A formula for the Chow dilogarithm on elliptic curves is given 
in Chapter 6. 

The function ${\cal L}_2$  satisfies  Abel's five term functional equation: 
\begin{equation} \label{6.14.02.1}
\sum_{i=1}^{5}(-1)^i {{\cal L}}_2(r(x_1,...,\widehat  x_i,...,x_{5})) =0 
\end{equation}
The Chow dilogarithm also satisfies functional equations. They appear as a reformulation 
of 
the fact that the composition 
$$
{\cal Z}_4(Spec(\C); 2) \stackrel{\partial}{\lra} {\cal Z}_3(Spec(\C); 2) 
\stackrel{{\cal P}_2}{\lra} \R(1)
$$ 
is zero. Namely,  
let $Y$ be an algebraic surface with four rational functions 
$g_1, ..., g_4$ on it corresponding to an element of ${\cal Z}_4(Spec(\C); 2)$. 
To evaluate the composition on this element 
we do the following. 
Take the divisor ${\rm div} (g_i)$ and 
restrict the other functions $g_i$ 
to it. Then applying the Chow dilogarithm to the obtained data and 
taking the alternating sum over $1 \leq i \leq 4$ we get zero. 
In the special case when $Y = \C\PP^2$  and ${\rm div}g_i = (l_i) - (l_{5})$, 
where $l_1, ...,  l_{5}$ are five  lines in the plane, this 
functional equation  plus (\ref{6.14.02.2}) is equivalent to   Abel's equation 
(\ref{6.14.02.1}).

{\bf 5. The Grassmannian $n$--logarithm and  symmetric space 
$SL_n(\C)/SU(n)$}. 
Restricting the Chow $n$-logarithm function to the subvariety of $(n-1)$--planes in 
$\C\PP^{2n-1}$ in general position with respect to the simplex $L$ 
we get the Grassmannian $n$--logarithm function ${\cal L}^G_n$. 

Let $G$ be a group and $X$ a $G$-set. {\it Configurations} of $n$ points in $X$ are 
by definition the points of the quotient $X^n/G$. There is a  natural bijection
$$
\mbox{ $\{ (n-1)$--planes in 
 $\PP^{2n-1}$ in generic position with respect to a simplex $L$\}}/({\Bbb G}_m^*)^{2n-1} \quad 
$$
$$  <--> \quad \left \{ \mbox{Configurations of} \quad 2n \quad\mbox {generic hyperplanes in } \PP^{n-1}\right \}
$$
given by intersecting of an $(n-1)$--plane $h$ with the codimension one faces of $L$. 

\begin{figure}[ht]
\centerline{\epsfbox{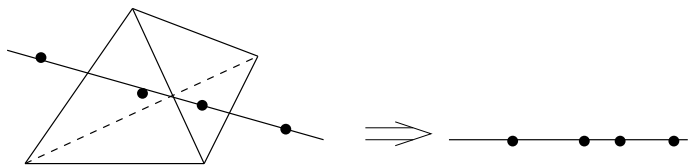}}
\caption{Toric quotients of Grassmannians and configurations of hyperplanes}
\label{gpol6}
\end{figure}

Using it 
we can view ${\cal L}^G_n$ as a function on the configurations of $2n$ hyperplanes in $\C \PP^{n-1}$. Applying the projective duality we can consider it as a function 
on configurations of $2n$ points in $\C \PP^{n-1}$.

In fact one can define the Grassmannian $n$-logarithm 
${{\cal L}}^G_n(x_1,...,x_{2n})$  as a 
function 
 on configurations of {\it arbitrary}  $2n$ points in $\C\PP^{n-1}$, see 
 Chapter 4. 
It is a measurable function which is 
real analytic on generic configurations. 
It satisfies the two functional equations
\begin{equation} \label{6.13.02.100}
\sum_{i=0}^{2n}(-1)^i {{\cal L}}^G_n(x_0,...,\widehat  x_i,...,x_{2n}) =0, \quad \sum_{j=0}^{2n}(-1)^j {{\cal L}}^G_n(y_j|y_0,...,\widehat  y_j,...,y_{2n}) = 0
\end{equation}
In the second formula $(y_0,...,y_{2n})$ is a configuration
of $2n+1$ points in $\C\PP^{n}$ and
$(y_j|y_0,...,\widehat  y_j,...,y_{2n})$ 
is a configuration of $2n$ points in $\C\PP^{n-1}$
obtained by projection from  $y_j$. 

\begin{figure}[ht]
\centerline{\epsfbox{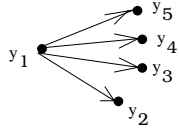}}
\caption{The configuration $(y_1|y_2, y_3, y_4, y_5)$ on $\PP^1$}
\label{gpol20}
\end{figure}

It follows from (\ref{6.14.02.2}) that 
the Grassmannian dilogarithm is given by the Bloch-Wigner function:
\begin{equation} \label{6.13.02.101}
{\cal L}^G_2(z_1, ..., z_4) = {\cal L}_2(r(z_1, ..., z_4))
\end{equation}
Abel's five term equation coincides with 
(\ref{6.13.02.100}). (The two functional equations (\ref{6.13.02.100}) 
are equivalent when $n=2$). 

Lobachevsky discovered that the dilogarithm 
appears in the computation of volumes of geodesic 
simplices in the three dimensional hyperbolic space ${\cal H}_3$.  
Let $I(z_1,...,z_4)$ be  the ideal geodesic simplex with vertices at the 
points $z_1,...,z_4$ on the absolute of ${\cal H}_3$. 
The absolute is naturally identified with 
$\C\PP^1$. 
\begin{figure}[ht]
\centerline{\epsfbox{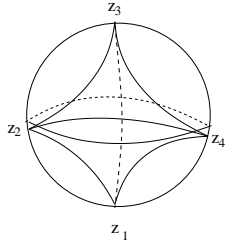}}
\caption{An ideal simplex in the hyperbolic $3$-space}
\label{gpol21}
\end{figure}
Lobachevsky's formula 
relates its volume to the Bloch-Wigner  function:
$$
{\rm vol}\Bigl(I(z_1,...,z_4)\Bigr) =   {\cal L}_2(r(z_1,...,z_4))
$$
The volume function ${\rm vol}I(z_1,...,z_4)$ is invariant under the group 
$SL_2(\C)$ of isometries of ${\cal H}_3$. So it  
depends only on the cross ratio 
of the points $z_1,...,z_4$. It   satisfies the five term equation (\ref{6.14.02.1}).
Indeed, $\sum(-1)^iI(z_1,..., 
\widehat z_i, ..., z_5) = \emptyset$. 
By Bloch's theorem  [Bl2] any measurable function $f(z)$ on $\C$
satisfying the five term equation  is proportional to ${\cal L}_2(z)$.
So we get the formula up to a constant. 

We generalize this picture as follows. 
 $\C\PP^{n-1}$ is realized as 
the smallest boundary 
stratum of the symmetric space ${\Bbb H}_n := SL_n(\C)/SU(n)$. We  define a function 
$\psi_n(x_1, ..., x_{2n})$  on configurations of $2n$ points 
of the symmetric space. The function $\psi_n$ 
is defined by an integral over $\C\PP^{n-1}$ similar to (\ref{6.17.02.1}).  
We show that it can be naturally extended to a function $\overline \psi_n$ 
on configurations of $2n$ points in a compactification $\overline {\Bbb H}_n$ of 
the symmetric space. The Grassmannian $n$--logarithm function 
turns out to be the value of the function $\overline \psi_n$ 
on configurations of $2n$ points at the smallest boundary strata, 
which is identified with $\C\PP^{n-1}$.

Now let  $n=2$. Then 
$SL_2(\C)/SU(2)$ is identified with the hyperbolic $3$-space. 
 We prove in Chapter 7 that  $\psi_2(x_1, x_2, x_3, x_{4})$ is 
 the volume of the geodesic simplex with vertices at the points 
$x_1, ..., x_4$. Restricting to the ideal geodesic simplices 
and using the relation to the Grassmannian dilogarithm 
plus  (\ref{6.13.02.101}) we 
get a 
new  proof of   Lobachevsky's formula.

{\bf 6. The Grassmannian $n$--logarithms and the Borel regulator}.
For any point $x \in \C\PP^{n-1}$ the function
\begin{equation} \label{6.17.02.2}
c_{2n-1}^n(g_1,...,g_{2n}):= { {\cal L}}^G_n(g_1x,...,g_{2n}x)
\end{equation} 
is a measurable $(2n-1)$-cocycle of the Lie group $GL_n(\C)$. 
Indeed, it is invariant
under the diagonal action of $GL_n(\C)$ and  the cocycle
condition is just the first functional equation for the function
$ {\cal L}^G_n$. Different points $x$ give canonically cohomologous cocycles. 
However {\it a priori} it is not clear that the corresponding   cohomology class 
  is non-zero.

Let $H^{2n-1}_{m}(GL_n(\C), \R)$ be the space of measurable 
cohomology of the Lie group $GL_n(\C)$. It is known that
$$
H^{\ast}_{m}(GL_n(\C), \R) = \Lambda_{\Bbb
R}^{\ast}(b_1,b_3,...,b_{2n-1}) 
$$
where $b_{2k-1} \in H^{2k-1}_{m}(GL_n(\C), \R)$ are certain
canonical generators  called the Borel classes ([Bo1]).
\begin {theorem}  
\label {0.4}
The cohomology class of the Grassmannian cocycle  (\ref{6.17.02.2}) is a non
zero rational multiple of the Borel class $b_{2n-1}$. 
\end{theorem}

For normalization of the Borel classes and precise relationship between the Grassmannian polylogarithms and the Borel regulator 
see Chapter 5, especially Sections 5.4 and 5.5. 

The essential role in  the proof is played by 
the fact that the Grassmannian $n$--logarithm function 
 ${\cal L}^G_n$ is a boundary value of the function  
$\overline \psi_n$. The function 
$\overline \psi_n(x_1,...,x_{2n})$ is not continuous at certain boundary points, but 
always satisfies the cocycle condition. 
So taking any point $x \in \overline {\Bbb H}_n$ we get a cocycle
$$
c_x(g_1, ..., g_{2n-1}):= \overline \psi_n(g_1x,..., g_{2n}x)
$$
of the group $GL_n(\C)$. Its cohomology class 
does not depend on  $x$. If $x \in {\Bbb H}_n$ 
the corresponding  cocycle  is smooth. We can differentiate it, getting   
a  cohomology class of the Lie algebra $gl_n$,  and relate it to the 
Borel class. On the other hand taking $x$ to be a point on the boundary stratum 
$\C\PP^{n-1}$ we recover the Grassmannian cocycle (\ref{6.17.02.2}). 
So we get the theorem. 

Combining it with the technique developed in [G1-2] we 
get  a simple explicit construction of the Borel regulator 
$$
K_{2n-1}(\C) \longrightarrow \R
$$
in terms of the Grassmannian $n$-logarithms.
 The second functional equation for ${\cal L}^G_n$ plays an important role in the proof. 
Therefore, thanks to the Borel theorem [Bo2], this allows 
to  express the  special values of Dedekind $\zeta$--functions at $s=n$ 
via the Grassmannian $n$--logarithms. 

The definition of the Higher Chow groups of a variety $X$ 
is much simpler than 
the definition of algebraic K-groups of $X$. The situation with 
the regulator maps is similar. However relating the 
special values of the Dedekind $\zeta$-functions to motivic 
cohomology of the corresponding  number fields
we need to work with the algebraic K-theory (or homology of $GL_n(F)$)  
of number fields.  

Chapter 2 is the main core of the paper. In Chapters 3, 4 and 6  
the main construction of Chapter 2 is investigated from different points of view. 
Chapters 4 and 5 are rather 
independent from the other Chapters. 

{\bf Acknowledgement}. 
I am grateful to Spencer Bloch for discussions of Arakelov 
motivic complexes. 
The results of this  paper were discussed in  my course at Brown (Spring 1998), 
and   at 
the Newton Institute (Cambridge, March 1998). 
I am grateful to participants for 
their interest and useful comments. I am very grateful to the referees 
for many comments which greatly improved the exposition.

During the work on this paper I enjoyed  the hospitality and support of the 
MPI(Bonn) and 
IHES (Bures-sur-Yvette). I am grateful to these 
institutions. 
This work was supported by the NSF grants DMS-9800998 
and  DMS-0099390.

 \section { Arakelov motivic complexes}

{\bf 1.  The  Higher Chow group complex}. A (non-degenerate) simplex in 
$\PP^m$ is an ordered
collection of hyperplanes $L_0,...,L_m$ in generic position, i.e. with empty
intersection. 
Let us choose in $\PP^m$ a simplex $L$ and a generic hyperplane $H$.
We might think about this data as of a simplex in the $m$-dimensional
affine space $\AAA^m:= \PP^m - H$.    For any two non-degenerate 
simplices in $\AAA^m$ there is a unique affine transformation sending one 
simplex to the other.

Let $I=(i_1,...,i_k)$ and $L_I:= L_{i_1} \cap ... \cap L_{i_k}$. 
Let $X$ be a regular projective variety over a field $F$. Let 
${\cal Z}_{m}(X; n)$ be
the free abelian group generated by irreducible codimension $n$
algebraic 
subvarieties in $X \times \AAA^m$ which  intersect properly (i.e.  with 
the  right codimension) all faces
$X \times L_I$.

{\bf Warning}. We use the notation  ${\cal Z}_{m}(X; n)$ for the group 
denoted  ${\cal Z}^{n}(X; m)$ by  Bloch. This allows us to use upper and lower 
indices to distinguish between the 
homological and cohomological notations, see below. 

For a given codimension 1 face $L_i$ of a simplex $L$ in $\AAA^m$ the  other 
faces $L_j$ cut a simplex $\widehat  L_i:= \{L_i \cap L_j\}$,  in $L_i$. 
So the  intersection with codimension 1 faces $X \times L_i$ provides 
 group homomorphisms 
$$
\partial_i: {\cal Z}_{m}(X; n) \longrightarrow 
{\cal Z}_{m-1}(X; n); \qquad \partial:= \sum_{i=0}^m (-1)^i \partial_i
$$
 Then   $\partial^2 =0$, so 
$({\cal Z}_{\bullet}(X; n); \partial)$    
is a homological complex. Its homology groups  
are Bloch's Higher
Chow groups. 
By the fundamental theorem of Bloch ([Bl1-2], [Lev])
$$
H_i({\cal Z}_{\bullet}(X; n)\otimes \Q)
 = K^{(n)}_{i}(X)\otimes \Q
$$
Let us cook up  a cohomological complex setting
$$
{\cal Z}^{\bullet}(X; n):= {\cal Z}_{2n-\bullet}(X; n)
$$ 
Its cohomology provides  a definition of the integral motivic cohomology of $X$:
$$
H_{{\cal M}}^i(X, \Z(n)) := H^i({\cal Z}^{\bullet}(X; n))
$$
Bloch's theorem guarantees Beilinson's formula (\ref{4.29.02.1}) 
for  the rational motivic cohomology.

{\bf 2.  The  Beilinson--Deligne complex}. 
Recall that an  $n$-distribution, sometimes also called an $n$-form with generalized function coefficients, or an $n$-current,  on a smooth  oriented manifold $X$ is 
 a continuous linear  functional on the space of 
$(\dim_{\R}X-n)$--forms with compact support.  Denote by 
${\cal D}^{n}_{X}$ the space of all {\it real} $n$--distributions on $X$.  
The  
De Rham complex of distributions 
$({\cal D}^{\bullet}_{X},d)$ is a resolution of the 
constant sheaf $\R$.  The space ${\cal A}^{n}_{X}$ of all smooth $n$--forms 
on $X$ is a subspace of 
${\cal D}^{n}_{X}$.

Let $X$ be a regular projective variety over $\C$.  The 
standard weight $n$  Beilinson-Deligne complex 
${\R}^{\bullet}(X; n)_{{\cal D}}$ 
is the total complex associated with the following
bicomplex: 
$$
\begin{array}{ccccccccccc} \label{del}
\Bigl({\cal D}_{X}^{0}&\stackrel{d}{\longrightarrow}&{\cal
D}_{X}^{1}&\stackrel{d}{\longrightarrow}&\ldots&\stackrel{d}{\longrightarrow}&{\cal
D}^{n}_{X}&\stackrel{d}{\longrightarrow}&{\cal
D}_{X}^{n+1}&\stackrel{d}{\longrightarrow}&\ldots\Bigr) \otimes \R(n-1)\\
&&&&&&&&&&\\
&&&&&&\uparrow\pi_{n} &&\uparrow\pi_{n}&&\\
&&&&&&&&&&\\
&&&&&&\Omega^{n}_{X}
&\stackrel{\partial}{\longrightarrow}&\Omega_{X}^{n+1}&\stackrel
{\partial}{\longrightarrow}&
\end{array}
$$
Here  $\R(n):= (2\pi i)^n\R$ and 
$$
\pi_n: {\cal D}_{X}^{p}\otimes \C \longrightarrow
{\cal D}_{X}^{p}\otimes \R(n-1)
$$ is the projection induced by the one  $\C = \R(n-1) \oplus \R(n)  \longrightarrow 
\R(n-1)$.  Further,  ${\cal D}^{0}_{X}$ placed in
degree 1 and  
$(\Omega^{\bullet}_{X}, \partial)$ is the De Rham complex of 
holomorphic forms. 

The Beilinson-Deligne complex ${\R}^{\bullet}(X; n)_{{\cal D}}$ 
is quasiisomorphic to the complex
$$
\R(n) \lra {\cal O}_X \lra \Omega^{1}_{X} 
\lra \Omega^{2}_{X}\lra ... \lra \Omega^{n-1}_{X}
$$


{\bf 3.  The  truncated Deligne complex}. 
Let ${\cal D}_X^{ p,q} = {\cal D}^{ p,q}$ be the abelian group of complex valued distributions of type $(p,q)$ on $X(\C)$. Consider the following cohomological 
``bicomplex'',  where ${\cal D}^{n,n}_{\rm cl} $ 
is the subspace of the space ${\cal D}^{n,n}$ of closed currents, and  
 ${\cal D}^{0,0}$ is in degree $1$:
$$
\begin{array}{ccccccccc}
&&&&&&&&{\cal D}_{cl}^{n,n}\\
&&&&&&&&\\
&&&&&&&2 \overline \partial \partial\nearrow&\\
&&&&&&&&\\
{\cal D}^{0,n-1}&\stackrel{\partial}{\longrightarrow}&{\cal D}^{1,n-1}&\stackrel{\partial}{\longrightarrow}&
...&\stackrel{\partial}{\longrightarrow}&{\cal D}^{n-1,n-1}&&\\
&&&&&&&&\\
\overline \partial \uparrow &&\overline \partial \uparrow &&&&\overline \partial \uparrow&&\\
 ...&...&...&...&...&...&...&&\\
 \overline \partial \uparrow&&\overline \partial \uparrow&&&&\overline \partial \uparrow&&\\
&&&&&&&&\\
{\cal D}^{0,1}&\stackrel{\partial}{\longrightarrow}& {\cal D}^{1,1}&\stackrel{\partial}{\longrightarrow}&... &\stackrel{\partial}{\longrightarrow}& {\cal D}^{n-1,1}&&\\
&&&&&&&&\\
\overline \partial \uparrow&&\overline \partial \uparrow&&&&\overline \partial \uparrow&&\\
&&&&&&&&\\
{\cal D}^{0,0}&\stackrel{\partial}{\longrightarrow} &{\cal D}^{1,0}& \stackrel{\partial}{\longrightarrow} &...&\stackrel{\partial}{\longrightarrow} &{\cal D}^{n-1,0}&&
\end{array}
$$
Properly speaking, it is not a bicomplex due to the presence of the operator 
$2\overline \partial \partial$, but we can handle it the same way we handle the bicomplexes. Namely, we define its total complex   
$Tot^{\bullet}$. It is concentrated 
in degrees $[1,2n]$. The   complex $C^{\bullet}_{{\cal D}}(X(\C); n) 
= C^{\bullet}_{{\cal D}}(n)$ is a  subcomplex 
of the complex 
$Tot^{\bullet}$ defined as follows. 
Take  the intersection of the part of the complex $Tot^{\bullet}$  
  coming from the $n \times n$ square in the diagram (and concentrated in degrees 
$[1,2n-1]$) with the complex of distributions with values in $\R(n-1)$. Consider 
the   subgroup    ${\cal D}_{\R, cl}^{n, n}(n) \subset {\cal D}_{cl}^{n, n}$  of
 the $\R(n)$-valued distributions of type $(n, n)$. They  form a subcomplex in 
$Tot^{\bullet}$ because
 $\overline \partial \partial$ sends $\R(n-1)$-valued distributions to 
$\R(n)$--valued distributions.
This is the complex $C^{\bullet}_{{\cal D}}(n)$. It is a truncation of the 
complex   considered by Deligne ([Del]). Its cohomology is the absolute Hodge cohomology defined by Beilinson [B3].

\begin{proposition} \label{dcom} Let $X$ be a regular complex projective 
variety. 
Then the    complex $C^{\bullet}_{{\cal D}}(X; n)$  is quasiisomorphic to the 
truncated Beilinson-Deligne complex $\tau_{\leq 2n}\R^{\bullet}(X; n)_{{\cal D}}$. 
\end{proposition}

{\bf Proof}. We need the following general construction. 
Let $f^{\bullet}: X^{\bullet} \lra Y^{\bullet}$ be a morphism of complexes such that
the map $f^i $ is injective for $i \leq p$ and surjective for $i \geq p$ (and hence is an isomorphism for $i=p$). Consider a complex
$$
Z^{\bullet}:= \quad {\rm Coker}f^{<p}[-1] \stackrel{D}{\lra} {\rm Ker}f^{>p}
$$
where the differential $D: {\rm Coker} f^{p-1} \lra {\rm Ker} f^{p+1}[1]$ is defined via the following diagram (the vertical sequences are exact):
$$
\begin{array}{ccccccccc}
&&&&&& 0 &&0\\
&&&&&&\downarrow &&\downarrow \\
0&&0 &&&&{\rm Ker}f^{p+1}&\lra &{\rm Ker}f^{p+2}\\
\downarrow &&\downarrow &&&&\downarrow &&\downarrow \\
X^{p-2} &\lra &X^{p-1} &\lra &X^{p}&\lra & X^{p+1}&\lra &X^{p+2}\\
\downarrow &&\downarrow &&f^p \downarrow = && \downarrow &&\downarrow \\
Y^{p-2}&\lra &Y^{p-1} &\lra &Y^{p} &\lra & Y^{p+1}&\lra &Y^{p+2}\\
\downarrow &&\downarrow &&&&\downarrow &&\downarrow \\
{\rm Coker} f^{p-2}&\lra &{\rm Coker} f^{p-1}&&&& 0&&0\\
\downarrow &&\downarrow &&&&&&\\
0&&0 &&&&&&
\end{array}
$$

\begin{lemma} \label{decom}
The complex $Z^{\bullet}$ is canonically quasiisomorphic to ${\rm Cone} (X^{\bullet} \stackrel{f^{\bullet}}{\lra} Y^{\bullet})$.
\end{lemma}

{\bf Proof}. Let
$$
\widetilde  \tau_{<p}X^{\bullet}:= \quad ... \stackrel{d_X}{\lra}  X^{p-2} \stackrel{d_X}{\lra}  X^{p-1} \stackrel{d_X}{\lra}  {\rm Im} d_X
$$
$$
\widetilde  \tau_{\geq p}Y^{\bullet}:= \quad Y^{p}/{\rm Im} d_Y  \stackrel{d_Y}{\lra}  Y^{p+1} 
\stackrel{d_Y}{\lra}  Y^{p+2} \stackrel{d_Y}{\lra} ...
$$
Then there is an exact sequence of complexes $0 \lra \widetilde  \tau_{<p}X^{\bullet} \lra X^{\bullet} \lra \widetilde  \tau_{\geq p}X^{\bullet} \lra 0$. The conditions on the maps $f^{\bullet}$ imply that 
$$
\widetilde  \tau_{<p}f^{\bullet}: \quad \widetilde  \tau_{<p}X^{\bullet} \lra \widetilde  \tau_{<p}Y^{\bullet} \quad \mbox{is injective}
$$
$$
\widetilde  \tau_{\geq p}f^{\bullet}: \quad \widetilde  \tau_{\geq p}Y^{\bullet} \lra \widetilde  \tau_{\geq p}Y^{\bullet} \quad \mbox{is surjective}
$$
We get maps of complexes
$$
{\rm Cone}(\widetilde  \tau_{<p}X^{\bullet} \lra f^{\bullet}(\widetilde  \tau_{<p}X^{\bullet})) \stackrel{\alpha}{\hookrightarrow}
{\rm Cone}(X^{\bullet} \lra Y^{\bullet}) \stackrel{\beta}{\longrightarrow} 
{\rm Cone}(\widetilde  \tau_{\geq p}Y^{\bullet} \lra \widetilde  \tau_{\geq p}Y^{\bullet}) 
$$
where $\alpha$ is injective and $\beta$ is surjective.
The complex ${\rm Ker} (\beta)/{\rm Im}(\alpha)$ looks as follows:
$$
\begin{array}{ccccccccc}
0&\lra &0 &\lra & {(f^p)}^{-1}{\rm Im}(d_Y)/{\rm Im}(d_X)&\lra &{\rm Ker} f^{p+1}&\lra &{\rm Ker} f^{p+2}\\
\downarrow &&\downarrow &&\downarrow &&\downarrow &&\downarrow \\
{\rm Coker} f^{p-2}&\lra &{\rm Coker} f^{p-1}&\lra & {\rm Im}(d_Y)/{f^p}{\rm Im}(d_X)  &\lra &0&\lra &0
\end{array}
$$
Since the map $f^p : {\rm Im}(d_Y)/{\rm Im}(d_X) \to {\rm Im}(d_Y)/{f^p}{\rm Im}(d_X)$ is an isomorphism it is quasiisomorphic to $Z^{\bullet}$. The lemma is proved.

Applying the lemma to the morphism of complexes 
$$
{\rm Tot}({\cal D}^{\geq n, \bullet}) \stackrel{\pi_n}{\longrightarrow} {\cal D}^{n+\bullet}\otimes_{\R}\R(n-1)
$$
we see that the complex $\R(n)_{{\cal D}}$ is canonically quasiisomorphic to the following complex:

\begin{figure}[ht]
\centerline{\epsfbox{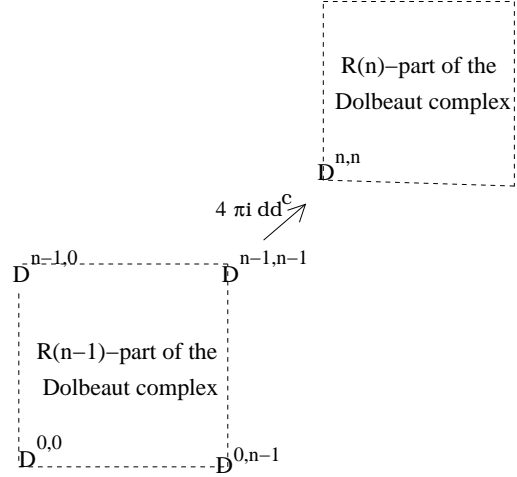}}
\caption{The weight $n$ real Deligne complex }
\label{gpol1}
\end{figure}
To compute the differential ${\cal D}_{\R}^{n-1, n-1}(n-1) \lra {\cal D}_{\R}^{n, n}(n)$ we proceed as follows. Take  $\alpha \in {\cal D}_{\R}^{n-1, n-1}(n-1)$, so $\alpha = (-1)^{n-1}\overline \alpha$. Then 
$$
d \alpha \quad = \quad \partial \alpha + \overline \partial \alpha \quad = \quad \partial \alpha + (-1)^{n-1}\overline {\partial \alpha } \quad = \quad 2 \pi_n(\partial \alpha )
$$
Applying $d = \partial + \overline \partial$ again and taking the $(n,n)$-component we get $2\overline \partial \partial (\alpha )$. 
Truncating this complex we obtain the proof of the proposition. 
(Note that $d^{C}:= 
(4\pi i)^{-1}(\partial - \overline \partial)$, 
so $dd^{C} = (2\pi i)^{-1}\overline \partial\partial$.)

Now if $X$ is a variety over $\R$, then we set
$$
C^{\bullet}_{{\cal D}}(X_{/{\R}}; n):= C^{\bullet}_{{\cal D}}(X; n)^{\overline F_{{\infty}}}; \qquad H_{{\cal D}}^i(X_{/{\R}}; \R(n)) := 
H^i\Bigl(  C^{\bullet}_{{\cal D}}(X_{/{\R}}; n)\Bigr)
$$
where $\overline F_{{\infty}}$ is the De Rham involution,  i.e.  the composition  of the involution $F_{{\infty}}$   on $X(\C)$ induced by the complex conjugation  with the complex conjugation of coefficients.

\begin{theorem-construction} \label{6.11.02.1} 
Let $X$ be a regular complex projective variety. Then 
there exists a {\rm canonical} homomorphism of complexes 
$$
{\cal P}^{\bullet}(n): {\cal Z}^{\bullet}(X; n) 
\longrightarrow C^{\bullet}_{{\cal D}}(X; n)
$$
If $X$ is defined over $\R$ then the image of the map  
${\cal P}^{\bullet}(n)$  lies in  
the subcomplex $C^{\bullet}_{{\cal D}}(X_{/\R}; n)$. 
\end{theorem-construction}



To construct this homomorphism we need to define certain 
homomorphism $r_{n-1}$ ([G5]). In the next section we 
recall its definition  and establish  its basic properties. 
Using it we define an $(m-1)$-form 
$r_{m-1}(L; H)$ canonically attached to the pair $(\AAA^m; L) = 
(\PP^m - H, L)$, and then define the homomorphism ${\cal P}^{\bullet}(n)$.

{\bf 4. The homomorphism $r_{m-1}$}.    
Let $X$ be a variety over $\C$. Let $f_1,...,f_m$ be $m$ rational functions on $X$.
 We attach to them the   $(m-1)$-form  
\begin {equation} \label{1wq}
r_{m-1}(f_1,..., f_m) :=
\end {equation}
$$
 {\rm Alt}_m \sum_{j\geq 0, 2j+1\leq 2m+1} c_{j,m}\log|f_1|d\log|f_2|
\wedge ... \wedge d\log|f_{2j+1}|\wedge di\arg f_{2j+2}\wedge ... \wedge
di\arg f_{m}
$$
Here $c_{j,m}:= \frac{1}{(2j+1)!(m-2j-1)!} $ and ${\rm Alt}_m$ is the operation of alternation: 
$$
{\rm
Alt}_m F(x_1,...,x_m):= \sum_{\sigma \in S_m}(-1)^{|\sigma|}F(x_{\sigma
(1)},...,x_{\sigma (m)})
$$

So $r_{m-1}(f_1,..., f_m)$ is an $\R(m-1)$-valued $(m-1)$-form and it
is easy to check that 
\begin{equation} \label{6.16.04.3}
d r_{m-1}(f_1,..., f_m) = \pi_m \Bigl(
d \log f_1 \wedge ... \wedge d \log f_m\Bigr) \quad  
\end{equation}

The form (\ref{1wq}) is a part of a cocycle representing the product in real Deligne cohomology of
1-cocycles $(\log|f_i|, d\log f_i)$. 

Here is a yet another, a bit more general 
 way to look at the homomorphism $r_{m-1}$.  Let 
${\cal A}^i(M)$ be the space of smooth $i$-forms on a real smooth manifold $M$. 
Consider the following map
\begin{equation} \label{6.16.04.1}
\omega_{m-1}: \Lambda^{m}{\cal A}^0(M) \lra {\cal A}^{m-1}(M)
\end{equation}
$$
\omega_{m-1}(\varphi_1 \wedge ... \wedge \varphi_{m}) := 
$$
$$
\frac{1}{m!}{\rm Alt}_{m}
\Bigl(\sum_{k=1}^{m}
(-1)^{k-1}\varphi_1\partial \varphi_2 \wedge ... \partial \varphi_k \wedge \overline 
\partial \varphi_{k+1} \wedge ... \wedge \partial \varphi_{m}\Bigr)
$$
For example 
$$
\omega_0(\varphi_1) = \varphi_1;\qquad 
\omega_1(\varphi_1 \wedge \varphi_{2}) = \frac{1}{2}\Bigl(\varphi_1\partial \varphi_2  - 
\varphi_2\partial \varphi_1 - \varphi_1\overline \partial \varphi_2  +
\varphi_2\overline \partial \varphi_1   \Bigr)
$$
Then one easily checks that 
\begin{equation} \label{6.16.04.2}
d\omega_{m-1}(\varphi_1 \wedge ... \wedge \varphi_{m}) = 
\partial \varphi_1 \wedge ... \wedge \partial \varphi_{m} 
+(-1)^m \overline \partial \varphi_1 \wedge ... \wedge \overline \partial \varphi_{m} +
\end{equation}
$$
\sum_{i=1}^{m}(-1)^{i}\overline \partial \partial \varphi_i \wedge 
\omega_{m-2}(\varphi_1 \wedge ... \wedge \widehat \varphi_i \wedge ... \wedge \varphi_{m})
$$

Now let $f_i$ be rational functions on 
a complex algebraic variety $X$. Set $M:= X^0(\C)$, where $X^0$ is the open part 
of $X$ where the functions $f_i$ are regular. Then $\varphi_i := \log|f_i|$ 
are smooth functions on $M$, and we have an identity 
$$
\omega_{m-1}(\log |f_1| \wedge ... \wedge \log |f_{m}|) = r_{m-1}(f_1 \wedge ... \wedge f_{m})
$$

Observe that $\overline \partial \partial \log|f| = 0$ 
on $X^0(\C)$. Therefore the second term in  the formula (\ref{6.16.04.2}) is zero, and so this formula 
is consistent with the one (\ref{6.16.04.3}). Notice however that 
if we understood
$\overline \partial \partial \log|f| $ as a distribution on $X(\C)$, then 
by the Poincar\'e-Lelong formula one has 
\begin{equation} \label{6.16.04.4}
2\overline \partial \partial \log|f| = 2\pi i \delta(f)
\end{equation}  

Our next goal is to interpret the form $r_{m-1}(f_1 \wedge ... \wedge f_{m})$ as a distribution 
on $X(\C)$ and calculate the differential of this distribution, taking into account 
formula  (\ref{6.16.04.4}).

{\bf 5. The distribution $r_{m-1}(f_1 \wedge ... \wedge f_{m})$}. 
Recall (see for instance [S]) 
that for a subvariety $Y$ of a smooth complex variety $X$ 
we define the $\delta$-distribution $\delta_{Y}$ by setting 
$$
<\delta_{Y}, \omega>:= \int_{Y^0(\C)} \omega
$$ 
where $Y^0$ is the nonsingular part of $Y$. 

\begin{theorem} \label{8.6.02.2}
 Let $Y$ be an arbitrary irreducible subvariety of a smooth complex variety 
$X$ and $f_1, ..., f_m \in {\C}^*(Y)$. Then for any smooth differential 
form $\omega$
with compact support on $X(\C)$ the following integral is convergent:
$$
\int_{Y^0(\C)}r_{m-1}(f_1, ..., f_m) \wedge i_Y^0\omega
$$ 
Here $Y^0$ is the nonsingular part of $Y$ and $i_Y^0\omega$ is the
restriction of the form $\omega$ to $Y^0(\C)$. 
Thus the form $r_{m-1}(f_1,..., f_m)$ defines a distribution 
$r_{m-1}(f_1, ...,f_m)\delta_Y$ on 
$X(\C)$ given by 
$$
<r_{m-1}(f_1, ...,f_m)\delta_Y, \omega>:= 
\int_{Y^0(\C)}r_{m-1}(f_1, ..., f_m) \wedge i_Y^0\omega
$$
 It provides  a group homomorphism 
\begin{equation} \label{2}
r_{m-1}: \Lambda^m\C(Y)^{\ast} \longrightarrow {\cal D}^{m-1}_{X(\C)}(m-1) 
\end{equation} 
\end{theorem}

{\bf Proof}. We need the following lemma

\begin{lemma} \label{8.6.02.112} 
Let $Y$ be a smooth complex projective variety. Then 
for any non zero rational functions $f_1,...,f_m$ on $X$ and 
for any smooth form $\omega$
with compact support on $Y(\C)$ the integral 
$$
\int_{Y(\C)}r_{m-1}(f_1, ..., f_m) \wedge \omega
$$
is convergent. So the form
$r_{m-1}(f_1, ..., f_m)$ defines a distribution 
on $Y(\C)$. 
\end{lemma}

{\bf Basic example}. The integral $\int_{\C}\log|z|d\log(z-a) \wedge
d\log \overline {(z-b)}$ is 
divergent at infinity, where all the functions $z, z-a, z-b$ have a
simple pole, since $\int_{\C}\log|z|\frac{dz\wedge d\overline z}{|z|^2}$ 
is divergent (both near zero and infinity). However 
$$
4\cdot \int_{\C}\log|z|d\log|z-a| \wedge
d\log |z-b| = 
$$
$$
\int_{\C}\log|z|\Bigl( d\log(z-a) \wedge d\log \overline {(z-b)} \quad + \quad d\log \overline {(z-a)} \wedge d\log  (z-b)\Bigr)
$$ 
is convergent: the divergent parts cancel each other. 
In $r_{m-1}(f_1,...,f_{m})$ such divergences
cancel because of multiplicativity and skew-symmetry of $r_{m-1}$.

{\bf Proof of Lemma \ref{8.6.02.112}}.    
Resolving singularities we reduce the statement of the lemma to the case when divisors 
${\rm div}f_i$ 
have normal crossing. Using the fact that  $r_{m-1}$ is a homomorphism 
to differential forms 
we may suppose that these divisors are 
different.  Our statement is local, so we can assume  that in
local coordinates $z_1,...,z_m$ one has $f_1 = z_1,...,f_k = z_k$ and
${\rm div}f_j$ for $j>k$ does not intersect  the origin.
After this the statement of lemma is obvious: each term in (\ref{1wq})
defines
a distribution near the origin. For instance the worst possible singularities have 
the term $\log|z_1|d\log|z_2|
\wedge ... \wedge d\log|z_{k}|\wedge \omega$ where $\omega$ is smooth
near the origin. It is clearly integrable with a smooth test form. The lemma is proved.

{\bf Remark}. In particular if ${\rm dim}_{\C}X=n$  the integral 
\begin{equation} \label{a5}
\int_{X(\C)}r_{2n}(f_1,...,f_{2n+1})
\end{equation}
is convergent.

Below we use the following form of the resolution of singularities
theorem. Recall that the proper preimage $\widetilde Y$ 
is the closure in $\widetilde X$ 
of the preimage of
the generic part of $Y$.

\begin{theorem} \label{8.9.02.4} 
Let $Y$ be an arbitrary 
 subvariety of a regular variety $X$ over a characteristic
zero field, and $Z$ 
a divisor of $Y$. Then there exists a sequence of blow ups 
$\widetilde X \lra X_1 \lra ... \lra  X$ providing a projection 
 $\pi: \widetilde X \to
X$ such that the proper preimage $\widetilde Y$ of $Y$ 
is non singular, 
$\widetilde Z:= p^*Z$, where $p= \pi|_{\widetilde Y}:\widetilde Y \to Y$, 
  is a normal crossing divisor in $\widetilde Y$, 
and restriction of  $\pi$ to the nonsingular part 
of $\widetilde Y- \widetilde Z$ 
is an isomorphism. 
\end{theorem}

{\bf Proof of Theorem \ref{8.6.02.2}}. By Theorem \ref{8.9.02.4} 
there exists a sequence of blow ups providing a projection 
$\pi: \widetilde X \lra X$ such that 
the proper preimage  $\widetilde Y$ of $Y \subset X$ is
smooth.  
By the above lemma the integral 
\begin{equation} \label{8.6.02.1}
\int_{\widetilde Y(\C)}r_{m-1}(\pi^*f_1, ..., \pi^*f_m) \wedge \pi^*\omega
\end{equation}
is convergent. Therefore the similar integral
over any Zariski dense subset of $\widetilde Y(\C)$ is also convergent and  
coincides with (\ref{8.6.02.1}). Since 
$\pi$ is an isomorphism on the  nonsingular part of 
 $\widetilde Y- \widetilde Z$, we are done.  
The Theorem \ref{8.6.02.2} is proved.

Below we employ notation $r_{m-1}(f_1 \wedge ... \wedge f_m)$ 
for the distribution given by (\ref{2}).

{\bf 6. Differential of the distribution 
$r_{n-1}(f_1 \wedge ... \wedge f_n)$}. 
Let $X$ be a normal variety. Then there is the residue homomorphism
$$
{\rm Res}: \Lambda^n\C(X)^* \lra \oplus_{Y \subset X^{(1)}}\Lambda^{n-1}\C(Y)^*, \qquad 
$$
where the sum is over all irreducible divisors of $X$. 

Here is its definition. Let $K$ be a field with 
 a discrete valuation $v$ and  the residue field $
k_v$.  The group of units $U$ has a natural 
homomorphism $U\longrightarrow k_v^{\ast}\; , \; u\mapsto 
\overline u$.  An element $\pi \in K^{\ast}$ is prime if 
${\rm ord}_{v}\pi = 1$. There is a homomorphism ${\rm res}_v:\Lambda^{n}K^{\ast} 
\longrightarrow\Lambda^{n-1} k_v^{\ast}$ uniquely defined 
by the  properties $(u_{i}\in U)$: 
$$
{\rm res}_v (\pi\wedge u_{1}\wedge \cdots\wedge u_{n-1}) 
= \overline u_{1}\wedge\cdots \wedge \overline u_{n-1}\quad  \mbox{and}  \quad 
{\rm res}_v
(u_{1}\wedge \cdots \wedge u_{n}) = 0
$$ 
It  does not depend on the choice of $\pi$.

Observe that if $X$ is normal then 
the local ring of any irreducible divisor of 
$X$ is a discrete valuation ring, so we can apply the above construction. 
We set ${\rm Res}:= \sum {\rm res}_{v}$ where the sum is over 
all valuations of the field $\C(X)$ 
corresponding to the codimension one points of $X$. 

 {\bf Remark}. If for any $i$ the restrictions of 
the functions $f_j$ for $j \not = i$ 
to the generic points of all irreducible components of 
the divisor  ${\rm div}f_i$ are non zero, then 
\begin{equation} \label{8.12.02.100}
{\rm Res}(f_1 \wedge ... \wedge f_n) =  \sum_{Y \in X^{(1)}}\sum_{i=1}^n(-1)^{i-1}v_Y(f_i) 
\cdot {f_1}_{|_Y} \wedge ... \wedge \widehat  {f_i}_{|_Y}\wedge ... 
\wedge {f_n}_{|_Y}
\end{equation}
where $v_Y(f)$ is the order of zero of $f$ at the generic point of $Y$.

Let $f_1 \wedge ...\wedge  f_n \in \Lambda^n \C(Y)^*$ where $Y$ is a
subvariety of a regular complex variety $X$. We define 
a distribution 
\begin{equation} \label{8.9.02.1}
(r_{n-2}\circ {\rm Res})(f_1 \wedge ...\wedge  f_n)
\end{equation}
on $X(\C)$ as follows.

i) If $Y$ is a normal then we have defined 
the residue map ${\rm Res}$ on $\Lambda^n \C(Y)^*$. 
So we define 
(\ref{8.9.02.1}) as 
$$
\sum_{Z \in Y^{(1)}}r_{n-2}{\rm res}_Z(f_1 \wedge ...\wedge  f_n)\delta_Y
$$
i.e. for any smooth form $\omega$ on $X(\C)$ 
$$
<(r_{n-2}\circ {\rm Res})(f_1 \wedge ...\wedge  f_n), \omega> := 
\sum_Z \int_{Z^0(\C)}r_{n-2}{\rm res}_Z(f_1\wedge  ... \wedge f_n)\wedge \omega
$$

ii) For an arbitrary $Y$ we take the normalization 
$\pi: Y^{\nu} \to Y$ and define (\ref{8.9.02.1}) as 
$\pi_*((r_{n-2}\circ {\rm Res}) (\pi^*f_1 \wedge ...\wedge \pi^* f_n))$

\begin{lemma} \label{Lemma 2.7} Let $Y$ be a subvariety of a regular
complex variety $X$,  $f_1 \wedge ...\wedge  f_n \in \Lambda^n
\C(Y)^*$, and  $Z:= \cup_i {\rm div}f_i$. 
 Let $\pi: \widetilde X \to X$ is a blow up of $X$ 
as in Theorem \ref{8.9.02.4}.  Then 
\begin{equation} \label{8.9.02.3}
\pi_*\Bigl((r_{n-2}\circ {\rm Res})(\pi^* f_1 \wedge
... \wedge \pi^*  f_n) \Bigr) = 
(r_{n-2}\circ {\rm Res})(f_1 \wedge
... \wedge f_n) 
\end{equation}
\end{lemma}

{\bf Proof}. Thanks to the definition  ii) 
of distribution  (\ref{8.9.02.1}) for non normal varieties 
we may assume without loss 
of generality 
that $Y$ is normal.

Let $Z = \cup_{i \in I}Z_i$ be the decomposition of the divisor $Z$ into
irreducible components parametrised by a set $I$. 
By the very definition the right hand side of (\ref{8.9.02.3}) 
is a sum over  $i \in
I$ of 
distributions $\psi_i= r_{n-2}(F_i)$ on
$Z_i$   corresponding to certain elements 
$F_i \in \Lambda^{n-2}\C(Z_i)^*$. 

Let $\widetilde Z = \cup_{j \in J}\widetilde Z_j$ 
be the decomposition of  $\widetilde Z$ into
irreducible components parametrised by a set $J$.  
The left hand side of (\ref{8.9.02.3}) 
is a sum over $j \in J$ of distributions $\widetilde \psi_j 
 = r_{n-2}(\widetilde F_j)$  corresponding to certain elements 
$\widetilde F_j \in \Lambda^{n-2}\C(\widetilde Z_j)^*$. 
One has $I \subset
J$ since the  proper preimage $\widetilde Z_i$ 
of $Z_i$ is an irreducible component
of 
$\widetilde Z$.

The lemma follows from the following two claims:
$$
\pi_*(\widetilde \psi_i) = \psi_i, \quad i\in I; \qquad 
\pi_*(\widetilde \psi_j) = 0, \quad j\in J - I
$$
The first one is obvious since both distributions 
$\widetilde \psi_i$ and $\psi_i$ are defined  
by their restriction
to any nonsingular Zariski dense open part of the corresponding
divisor, and $\pi$, being  restricted to such a sufficiently small part 
of $\widetilde Z_i$, is an isomorphism. Let us prove the second claim. 
Observe that for $j \in J-I$ the subvariety $\pi(\widetilde Z_j)$ is
of
 codimension 
$2$ in $Y$, and restriction of $\pi$ to a Zariski open part of 
$\widetilde Z_j$ is a
fibration with fibers of positive dimension. We need to show that for
any smooth form $\omega$ on $X(\C)$ 
$$
\int_{\widetilde Z_j(\C)} \widetilde \psi_j \wedge \pi^*\omega 
:= \int_{\widetilde Z'_j(\C)} \widetilde \psi_j\wedge \pi^*\omega =0
$$
where $Z_j'$ is a (sufficiently small) 
nonsingular  Zariski dense open part of $Z_j$. 
Observe that $\widetilde \psi_j \wedge \pi^*\omega$ 
is a smooth form on $\widetilde
Z'_j(\C)$. For any vector field
$v$ on $\widetilde Z_j'$ tangent to the fibers of $\pi$  we  have 
$i_v\pi^*\omega =0$. So it suffices to show that 
$i_v \widetilde \psi_j =0$. Since the statement is local, 
we can choose a local equation of an open part of $\widetilde Z'_j$
in the form $\pi^*z =0$. Using it as a local
parameter at the definition of ${\rm res}_{\pi^*z =0}$, we see that 
${\rm res}_{\pi^*z =0}(\pi^*f_1 \wedge ... \wedge \pi^* f_n)$ is
lifted from $Y$, and hence so is 
$$
\widetilde \psi_j = r_{n-2}
 {\rm
  res}_{\pi^*z =0}
(\pi^*f_1 \wedge ... \wedge \pi^* f_n)
$$ 
Thus  $i_v \widetilde \psi_j =0$. The lemma is proved.

\begin{proposition} \label{8.9.02.6}  Let $Y$ be an arbitrary 
 subvariety of a regular
complex variety $X$ and   $f_1 \wedge ...\wedge  f_n \in \Lambda^n
\C(Y)^*$. Then
$$
d r_{n-1}(f_1\wedge ... \wedge f_n) = 
$$
\begin{equation} \label{p106}
\pi_n \Bigl(d \log f_1 
\wedge ... \wedge d \log f_n\Bigr) 
+ 2 \pi i \cdot (r_{n-2}\circ {\rm Res})
(f_1\wedge ... \wedge f_n)
\end{equation}
\end{proposition}

{\bf Proof}. Let us resolve singularities as in Lemma \ref{Lemma 2.7}. 
Since $\pi$ is a birational isomorphism,
 and the distribution 
$r_{n-1}(f_1 \wedge ... \wedge f_n)$ is determined by its 
restriction to the generic point, 
one has 
$$
\pi_*r_{n-1}( \pi^* f_1 \wedge ... \wedge \pi^*  f_n) = 
r_{n-1}( f_1 \wedge ... \wedge f_n)
$$
So this and 
 Lemma  \ref{Lemma 2.7} imply that we may assume that $\cup_i {\rm div}f_i$ 
is a normal crossing divisor. The proposition 
follows immediately from the Poincar\'e-Lelong formula  (\ref{6.16.04.4})
 $$
d(d i 
\arg f) = 2 \pi i \delta(f) := 2 \pi i \delta_{{\rm div}(f)}
$$
The Proposition \ref{8.9.02.6} is proved. 

{\bf Remark}. To prove the  Poincar\'e-Lelong formula one may 
resolve the singularities and argue just as in the proof of Lemma 
\ref{Lemma 2.7} that 
it is sufficient to prove this formula on a blow up. When the divisor of $f$  
is a normal 
crossing divisor the formula follows from 
$$
d(d i 
\arg z) = 2 \pi i \delta(z):= 2 \pi i (\delta_0 - \delta_{\infty}
)$$

{\bf 7. The distribution $r_{m-1}(L;H)$}. Let $\Omega_L$ be the canonical
$m$-form in ${\Bbb P}^m - L$ with logarithmic singularities
at $L$. It represents a generator of $H^m_{{\rm DR}}({\Bbb P}^m - L)$ defined over $\Z$.  Let us  give its coordinate description. Choose homogeneous coordinates
$(z_0:...:z_m)$ in $\PP^m$ such that 
$L_i $ is given by equation $\{z_i = 0\}$. Then 
$$
\Omega_L =  d\log z_1/z_0 \wedge ... \wedge d\log
z_m/z_0
$$
The form $\Omega_L$ has periods in $\Z(m)$.
So
$\pi_m(\Omega_L)$ is exact. However
there is no 
canonical choice of a primitive $(m-1)$-form for it: the group 
$(\C^{\ast})^m$ acting on $\C \PP^m - L$ leaves the form 
invariant and acts non trivially on the
primitives. 
But if we consider a simplex $L$ in the affine complex space $\AAA^m$
(or, what is the same, choose an additional hyperplane $H$ in  $\C \PP^m$, which should be thought of as the infinite hyperplane) then there
is a {\it canonical}  primitive.

Choose a coordinate system $(z_0:...:z_m)$ in $\PP^m$ as above 
such that $H$ is given by  $\{\sum_{i=1}^{m}
z_i =z_0\}$. Set
\begin{equation} \label{mint4}
r_{m-1}(L;H) := r_{m-1}(z_1/z_0 \wedge ... \wedge  z_m/z_0)
\end{equation}
Here is a more invariant definition. Choose one of the faces 
of the simplex $L$, say 
 $L_0$. Consider the simplex $(H, L_1, ... , L_m)$.
Let $f_i$ be the rational function on $\C \PP^m$ such that 
$(f_i) = L_i - L_0$ normalized by 
$f_i(l_i)= 1$, where $l_i$ is the vertex of the simplex 
$(H, L_1, ... , L_m)$  opposite to the face
$L_i$.  Then $f_i = \frac{z_i}{z_0}$ and 
$$
r_{m-1}(L;H) := r_{m-1}(f_1 \wedge ... \wedge  f_m)
$$

\begin{figure}[ht]
\centerline{\epsfbox{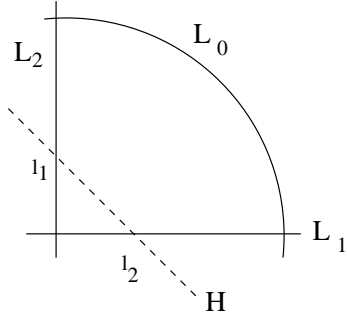}}
\caption{A simplex $L$ and an infinite hyperplane $H$}
\label{gpol5}
\end{figure}
This  form
 is skewsymmetric with respect to the permutation of the 
hyperplane faces of the simplex $L$. 
One has  
$$
d r_{m-1}(L;H) = \pi_m (\Omega_L) \quad \mbox{in} \quad \C \PP^m - L
$$
So a choice of an ``infinite''  hyperplane $H \subset \C 
\PP^{m}$ provides the form 
$r_{m-1}(L;H)$.

{\bf Example}. If $m=1$ then 
$$
\AAA^1 = \PP^1 - \{1\}, \quad L = \{0\} \cup \{\infty \}, \quad  \Omega_L = d\log z, \quad \pi_1 (d\log z) = d\log |z|, \quad 
$$
$$
r_0(\{0\} \cup \{\infty \}; \{1\}) =  \log |z|
$$

The $(n-1)$-form  $ r_{n-1}(L;H)$ provides 
an $(n-1)$-distribution on $\C \PP^n$.
Recall   the simplex $\widehat  L_i$ which is cuted out
 by $L$ in the hyperplane $L_i$, and put 
$H_i := L_i \cap H$. Consider the $(n-2)$-form 
$r_{n-2}(\widehat  L_i; H_i )$ on the hyperplane  
$L_i$  as $n$-distribution in $\C \PP^n$. We denote 
it as $r_{n-2}(\widehat  L_i; H_i ) \cdot \delta_{L_i}$.

\begin {corollary} \label{dif}
One has
\begin{equation} \label{1p1}
d r_{n-1}(L; H) =  \pi_n ( \Omega_L ) + 2 \pi i \cdot 
\sum_{i=0}^n (-1)^ir_{n-2}(\widehat  L_i; H_i ) \delta_{L_i}
\end{equation}
\end {corollary} 
 
{\bf Proof}. Follows immediately from Proposition \ref{8.9.02.6}.

{\bf 8.   A coordinate free description of  the form $r_{m-1}(L;H)$}. 
Let $V_m$ be an $m$-dimensional vector space over a field $F$. 
Choose a volume form ${\rm vol}_m \in {\rm det}V_m^{\ast}$. Set
$\Delta(v_1,...,v_m):= \langle{\rm vol}_m, v_1\wedge ...\wedge v_m\rangle \in F^{\ast}$.

\begin {lemma} For a configuration $(l_0,...,l_m)$ of $m+1$ vectors in
generic position 
\begin{equation} \label{4p}
f_m(l_0,...,l_m):= \sum_{i=0}^m(-1)^i \Lambda_{j\not = i}\Delta(l_0,...,\widehat  l_j,..., 
l_i,...,l_m) \in \Lambda^m F^{\ast}
\end{equation}
does not depend on the choice of the volume form ${\rm vol}_m$.
\end {lemma}

{\bf Proof}. See proof of Lemma 3.1 in [G3].

\begin{figure}[ht]
\centerline{\epsfbox{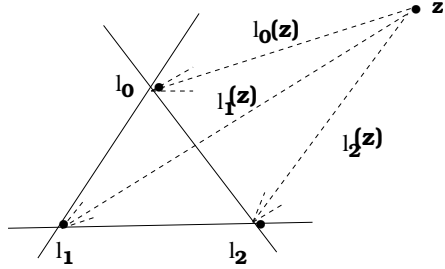}}
\caption{The vectors $l_i(z)$}
\label{gpol4}
\end{figure}

For a point $z\in \AAA^m- L$ let $l_i(z)$ be the
vector from 
$z$ to the vertex $l_i$, see Figure \ref{gpol4}.

We get a canonical element 
$$
f_m(l_0(z),...,l_m(z)) \in \Lambda^m\Q(\AAA^m -  L)^{\ast}
$$
If $F=\C$, applying the homomorphism $r_m$ to this 
element we get a canonical $(m-1)$-form in $\C \PP^{m} - L$.
It  coincides with  $ r_{m-1}(L;H)$. 

{\bf Example}. If $m=1$ and $(L_0, L_1, H)= (0, \infty, 1)$ then $t = \frac{z}{z-1}$ 
is an affine coordinate on $\PP^1 - \{1\}$ and $$
l_0(t) = \frac{z}{z-1}, 
l_1(t) = \frac{1}{z-1}, \quad \mbox{so} \quad  
f_1(l_0(t), l_1(t)) = \frac{l_0(t)}{l_1(t)} = z
$$ 

{\bf Remark}. The map $z\in \AAA^m -  L \lms f_m(l_0(z),...,l_m(z))$ 
provides an isomorphism 
$F_m: CH^m({\rm Spec}(F),0) \lra K_m^M(F)$; 
see [NS] where the isomorphism $F_m$ was presented in a bit different way.

{\bf 9. The main construction}. We have to construct a morphism of complexes
$$
\begin{array}{ccccccccc}
... &\longrightarrow&{\cal Z}^{1}(X; n)&\longrightarrow  &...&\longrightarrow&{\cal Z}^{2n-1}(X; n)&\longrightarrow&{\cal Z}^{2n}(X; n)\\
&&&&&&&&\\
&&\downarrow {\cal P}^{1}(n)&&... &&\downarrow {\cal P}^{2n-1}(n)&&\downarrow 
{\cal P}^{2n}(n)\\
&&&&&&&&\\
0&\lra &{\cal D}^{0,0}_{\R}(n-1)&  \stackrel{}{\longrightarrow}&...&\stackrel{}{\longrightarrow}&{\cal D}_{\R}^{n-1,n-1}(n-1)&\stackrel{2 \overline \partial \partial}{\longrightarrow}&{\cal D}_{\R}^{n,n}(n)
\end{array}
$$
Let $Y \in {\cal Z}^{2n}(X; n)$ 
be a codimension $n$ cycle in $X$. By definition 
$$
{\cal P}^{2n}(n)(Y):= (2\pi i)^n\delta_Y
$$

Let us construct homomorphisms
$$
{\cal P}^{2n-i}(n):  {\cal Z}^{2n-i}(X; n) 
\longrightarrow {\cal D}^{2n-i-1}_{X(\C)}(n-1), \quad i>0 
$$
 Denote by  $\pi_{{\AAA}^i}$ (resp. $\pi_{X}$) the projection of 
$X \times \AAA^i$ to $\AAA^i$ (resp. $X$), and by 
 $\overline \pi_{{\AAA}^i}$ (resp. $\overline \pi_{X}$) the projection of 
$X \times \C{\Bbb P}^i$ to $\C{\Bbb P}^i$ (resp. $X$). 

Recall the element 
\begin{equation} \label{mint3}
\frac{z_1}{z_0} \wedge... \wedge \frac{z_i }{z_0}
\in \Lambda^{i}\C(\AAA^i)^*
\end{equation}
defining the form $r_{i-1}(L;H)$, see (\ref{mint4}). Let 
\begin{equation} \label{mint2}
g_1 \wedge... \wedge g_i \in \Lambda^{i}\C(Y)^*
\end{equation} 
 be the restriction
to $Y$ of the inverse image of element 
(\ref{mint3}) by the projection $\pi^*_{\AAA^i}$. 
The element (\ref{mint2}) provides, by 
Theorem \ref{8.6.02.2}, a distribution on $X(\C) \times \C {\Bbb P}^i$. 
Pushing this distribution down by $(2\pi i)^{n-i}\cdot\pi_{\overline X}$ we get 
the distribution $ {\cal P}^{2n-i}(n)(Y)$:

\begin{definition} \label{mint6}
$$
{\cal P}^{2n-i}(n)(Y) := (2\pi i)^{n-i}\cdot
\pi_{\overline X *}r_{i-1}(g_1 \wedge... \wedge g_i) :=
$$
$$
(2\pi i)^{n-i}\cdot \pi_{\overline X *}r_{i-1}\Bigl( i_Y^*\pi^*_{\AAA^i}
(\frac{z_1}{z_0} \wedge... \wedge \frac{z_i }{z_0})\Bigr)
$$
\end{definition}
Here $i_{Y}: Y 
\hookrightarrow X \times  {\Bbb P}^i$

{\bf Remark}. This definition works if and only if 
the cycle $Y$ has proper intersection with all codimension one faces 
of $X \times L$. Indeed, if $Y$ does not have proper intersection 
with one of the faces, then the equation of this face restricts to zero 
to $Y$, and so (\ref{mint2}) does not make sence. 
As soon as all equations of the 
codimension one faces restrict to non zero functions on $Y$, 
(\ref{mint2}) makes sense, and we can apply  Theorem \ref{8.6.02.2}.

{\bf Remark}. We just proved that the product of distributions $\delta_Y \wedge \pi_{\AAA}^{\ast}r_{i-1}(L;H)$ makes sense and 
$$
{\cal P}^{2n-i}(n)(Y)= \quad (2\pi i)^{n-i}
{\pi_X}_*(\delta_Y \wedge \pi_{\AAA}^{\ast}r_{i-1}(L;H))
$$

It is handy to rewrite Definition \ref{mint6} more explicitly 
as an  integral over $Y(\C)$. Namely,  
let $\omega$ be a smooth form on $X(\C)$ and 
$Y \in  {\cal Z}^{2n-i}(X; n)$. Then 
\begin{equation} \label{mint}
<{\cal P}^{2n-i}(n)(Y),\omega>= \quad (2\pi i)^{n-i}
\int_{Y(\C)}\pi_{\AAA}^{\ast}r_{i-1}(L;H) \wedge \pi_X^{\ast}\omega =  
\end{equation}
\begin{equation} \label{8.9.02.10}
(2\pi i)^{n-i}\int_{Y^0(\C)}  
r_{i-1}(g_1\wedge ... \wedge g_i)\wedge i_{Y(\C)}^*\pi^*_X\omega 
\end{equation}
where $Y^0$ is the nonsingular part of $Y$

Since the form $r_{i-1}(L;H)$ is $\R(i-1)$-valued, for $i>0$ the
 distribution
 ${\cal P}^{2n-i}(n)(Y)$ takes values in $\R(n-1)$. 
Further,  ${\cal P}^{2n}(n)(Y)$ is obviously an $\R(n)$-valued distribution. 

Let us show that for $i>0$ the distribution ${\cal P}^{2n-i}(n)(Y)$ lies precisely 
in the left bottom $(n-1) \times (n-1)$ square of  the Dolbeault
bicomplex. The integral (\ref{mint}) 
is non zero only if 
$\pi^*_{X}\omega \wedge \pi^*_{{\Bbb A}}r_{i-1}(L;H)$ is of type 
$$
({\rm dim}Y, {\rm dim}Y) =  ({\rm dim}X+i-n, {\rm dim}X+i-n) 
$$
Since $r_{i-1}(L;H)$ is an 
$(i-1)$-form we see that the integral vanish if $\omega$ is 
a form of type $(p,q)$ where 
$p$ or $q$ is smaller then ${\rm dim} X+1-n$. This just means that 
the distribution lies in the left bottom $(n-1)\times (n-1)$ square of 
the Dolbeault bicomplex. 
The proposition is proved.

Therefore we have constructed the maps  ${\cal P}^{i}(n)$. 

\begin{theorem} \label{hhmmoo}
${\cal P}^{\bullet}(n)$ is a homomorphism of complexes.
\end{theorem}

{\bf Proof}. One has
\begin{equation} \label{98}
<d {\cal P}^{2n-i}(n)(Y), \omega> = 
%
\int_{\overline Y(\C)}dr_{i-1}(g_1\wedge
... \wedge g_i)\wedge 
i_{\overline Y}^*\pi_X^*\omega 
\end{equation}
where $g_1\wedge
... \wedge g_i$ is as in (\ref{mint2}). 
We use Proposition \ref{8.9.02.6} to calculate 
$dr_{i-1}(g_1\wedge... \wedge
g_i)$.

To handle the first term in (\ref{p106}) observe that 
by the very definition of the complex 
${\cal C}_{\cal D}(X;n)$ we need to investigate integral (\ref{98}) 
only for smooth forms $\omega$ of type $(p,q)$ where $|p-q| \leq i-1$. 
Since $\pi_i( \Omega_L ) = \Omega_L \pm \overline 
\Omega_L$ is a sum of forms of type $(i,0)$
and $(0,i)$ 
the form $\pi_{\AAA}^{\ast}\pi_i(\Omega_L)
\wedge \pi_X^{\ast}\omega$ 
can not be of type $(k,k)$.  
Therefore only the second term in (\ref{p106}) contributes.

Since $r_0(f) = \log|f|$, the commutativity of the last square 
follows from the Poincar\'e-Lelong formula 
$
2 \overline \partial \partial \log |f| = 2 \pi i \cdot \delta_{{\rm div}(f)} 
$.

{\it Commutativity of the $i$-th square of the diagram, 
$i>1$, counting from the right}. 
If $\overline Y$ is  normal this  
follows from Proposition \ref{8.9.02.6}. Indeed, 
since $Y$ meets the codimension two faces properly all equations of the codimension one
faces of  $X \times L$ but  $z_j$ have non zero restriction 
to the generic point of ${\rm div}z_j$. Therefore we may use formula 
(\ref{8.12.02.100}) to calculate ${\rm Res}$, and then the claim is obvious. 
 
In  the case when $\overline Y$ is not normal 
we face the following subtle problem.

 Calculating 
$dr_{i-1}(g_1\wedge... \wedge
g_i)$ and hence $d \circ {\cal P}^{2n-i}(n)(Y)$ when $Y$ is not
normal we need to take 
${\rm Res}_{\widetilde Z}(\widetilde g_1\wedge... \wedge
\widetilde g_i)$ for all irreducible divisors $\widetilde Z$ in the normalization
$\pi: \widetilde Y \to Y$, where $\widetilde g:= \pi^*g$,
and then take
\begin{equation} \label{8.9.02.7}
\sum_{\widetilde Z} \pi_*r_{i-2}\circ {\rm Res}_{\widetilde Z}
(\widetilde g_1\wedge... \wedge
\widetilde g_i)
\end{equation}
Computation of ${\cal P}^{2n-i+1}(n)(Z)\circ d $ 
does not involve the normalization of $Y$: we intersect
$Y$ with all codimension one faces of $X \times L$. 
So we need to compute (\ref{8.9.02.7}) using the intersection data of $Y$ and
$X \times L$.

To handle this  
 we use  the condition that $Y$ meets the codimension two faces of
the simplex $X \times L$ properly. It implies that 
\begin{equation} \label{8.9.02.20}
{\rm Res}_{\widetilde Z}(\widetilde g_1\wedge... \wedge
\widetilde g_i) = \pi^*G_{i-1}; \qquad  
G_{i-1}\in \Lambda^{i-1}\C(\pi(\widetilde Z))^*
\end{equation} 
Indeed, since all equations of the codimension one
faces of  $X \times L$ but $z_j$ have non zero restriction 
to the generic point of ${\rm div}z_j$,  the wedge 
product of these restrictions can be taken as $G_{i-1}$. 

Having  (\ref{8.9.02.20}) the statement is obvious since 
$$
\pi_* r_{i-2} \pi^*G_{i-1} = [\C(\widetilde Z): \C(\pi (\widetilde Z))]
\cdot r_{i-2}(G_{i-1})
$$ 
Indeed, recall (see [F], page 9) that if $\widetilde Y \to Y$ is the
normalization of  $Y$ and $g \in \C(Y)^* = \C(\widetilde Y)^*$ then 
$$
{\rm ord}_Z(g) = \sum_{\widetilde Z} 
{\rm ord}_{\widetilde Z}(g)[\C(\widetilde Z): \C(Z)]
$$
where the sum is over all irreducible 
divisors projecting onto $Z$. 

Theorem \ref{hhmmoo} is proved. Therefore we finished the proof of 
 Theorem-Construction 
\ref{6.11.02.1}.

{\bf 10. The Higher Arakelov Chow groups}. Let $X$ be a regular complex variety. 
Denote by  $\widetilde  C^{\bullet}_{{\cal D}}(n)$ the quotient of the complex $C^{\bullet}_{{\cal D}}(n)$ along the subgroup ${\cal A}_{cl}^{n,n}(n) \subset 
{\cal D}_{cl}^{n,n}(n) $ of closed smooth form of type $(n,n)$ with values in $\R(n)$. 

Consider the cone of the homomorphism ${\cal P}^{\bullet}(n)$ shifted by $-1$:
$$
\widehat  {\cal Z}^\bullet(X; n):= {\rm Cone}
\Bigl( {\cal Z}^\bullet(X; n)
 \stackrel{ }{\longrightarrow} \widetilde  C^{\bullet}_{{\cal D}}(X(\C); n)
\Bigr)[-1]
$$

\begin{definition} \label{arackd} The Higher Arakelov Chow groups are 
\begin{equation} \label{6.11.02.11}
\widehat  {CH}^n(X; i):= H^{2n-i}(\widehat  {\cal Z}^{\bullet}(X; n))
\end{equation}
\end{definition}

 Recall the arithmetic Chow groups  defined by Gillet-Soule [GS] as follows:
$$
\widehat  {CH}^n(X):= 
$$
\begin{equation} \label{6.11.02.10}
\frac{\{(Z,g); \frac { \overline \partial \partial}{\pi i} g+\delta_Z 
\in {\cal A}^{n,n}\}}{\{(0,\partial u + \overline \partial v); ({\rm div} 
f, -\log|f|), f \in \C(Y), 
{\rm codim} (Y) = n-1\}}
\end{equation} 
Here $Z$ is a divisor in $X$, $f$ is a rational function on a divisor $Y$ in $X$, 
$$
g \in {\cal D}_{\R}^{n-1,n-1}(n-1), \quad (u,v) \in ({\cal D}^{n-2,n-1} 
\oplus {\cal D}^{n-1,n-2})_{\R}(n-1)
$$

\begin{proposition} \label{ch}
 $\widehat  {CH}^n(X; 0) = \widehat  {CH}^n(X)$.
\end{proposition}

{\bf Proof}.  Let us look at the very right part of the complex 
$\widehat  {\cal Z}^{\bullet}(X; n)$:
$$
\begin{array}{ccccc}
...&\longrightarrow&{\cal Z}^{2n-1}(X; n)&\longrightarrow&{\cal Z}^{2n}(X; n) \\
&&&&\\
&&\downarrow {\cal P}^{2n-1}(n)&& \downarrow {\cal P}^{2n}(n)\\
&&&&\\
({\cal D}^{n-2,n-1} \oplus {\cal D}^{n-1,n-2})_{\R}(n-1)&\stackrel{(\partial , \overline \partial)}{\longrightarrow}&{\cal D}_{\R}^{n-1,n-1}(n-1)&\stackrel{2\overline \partial \partial}{\longrightarrow}&{\cal D}_{\R}^{n,n}(n)/{\cal A}_{\R}^{n,n}(n)
\end{array}
$$
Consider the very  end of the   Gersten complex on $X$:
$$
\prod_{Y \in X_{n-2}}\Lambda^2\C(Y)^* \stackrel{\partial}{\longrightarrow}  \prod_{Y \in X_{n-1}}\C(Y)^* \longrightarrow {\cal Z}_0(X; n) 
$$
where $\partial$ is the tame symbol. It maps to the complex 
$\widehat  {\cal Z}^{\bullet}(X; n)$, i.e. to the top row  of the bicomplex above, 
as follows. Recall that ${\cal Z}_0(X; n) = {\cal Z}^{2n}(X; n)$, so the very right component of our map is provided by this identification. 
Further, a pair $(Y;f)$ where $Y$ is an irreducible codimension $n-1$ subvariety of $X$ maps to the cycle $\{(y,f(y))| y \in Y\} \subset X \times \AAA^1$.
 Similarly  any element in $\Lambda^2\C(Y)^*$ can be represented as a linear combination
 of elements 
$\sum_i (Y;f_i\wedge g_i)$ where $Y$ is an irreducible codimension $n-2$ subvariety of $X$ and $f_i, g_i$ are rational functions on $Y$ such that ${\rm div} f_i$ and 
${\rm div} g_i$  share no irreducible divisors. Then we send $(Y;f_i\wedge g_i)$ to 
the cycle $(y,f_i(y),g_i(y)) \subset X \times \AAA^2$. 
It is well known that in this way we get  an isomorphism on the 
last two cohomology groups. 
Computing the composition of this map with the homomorphism 
${\cal P}^{\bullet}(n)$ we end up precisely with the denominator 
in (\ref{6.11.02.10}). The proposition is proved.

\section  { The Chow polylogarithms}                                                              

Suppose $X = {\rm Spec}(\C)$. Then ${\cal P}_n(Y): = {\cal P}^{1}(n)(Y)$ is a function on the space of all codimension $n$ cycles in $\PP^{2n-1}$ intersecting properly faces of the simplex $L$. It is called 
the {\it Chow polylogarithm function}. For $i>1$ all the distributions ${\cal P}^{i}(n)(Y)$ are zero. However modifying 
the construction of the previous chapter we get a very interesting object, the Chow polylogarithm,  even when $X$ is a point. The Chow polylogarithm function is the 
 first component of the Chow polylogarithm. One can define the Chow 
polylogarithm for an arbitrary variety $X$, but we spell out the details 
in the most 
interesting case when $X$ is a point.

{\bf 1   Chow polylogarithms [G5]}. Let $L = (L_0, ..., L_{p+q})$ 
be  a simplex in $\PP^{p+q}$,
 $H$ is a hyperplane in generic position to $L$, and  $H_i := H \cap L_i$. 

Let $ {\cal  Z}^{q}_{p}(L)$ be the variety of all codimension $q$
effective   algebraic cycles  in 
$\C \PP^{p+q}$ which intersect properly, i.e. each irreducible component in the right codimension, 
 all faces of the simplex $L$. 
 It  is a union of an infinite number of finite
dimensional complex algebraic varieties. 

{\bf Example}.   
$ \C \PP^n - L = (\C^{\ast})^n$ is an irreducible component of 
$ {\cal  Z}^{n}_{0}(L)$ parametrizing the irreducible subvarieties, i.e. points. 

Let $\widehat L_i$ be the simplex in the projective space $L_i$ cut by the hyperplanes $L_j$, $j \not = i$. 
The intersection of a cycle with a codimension
1 face $L_i$ 
of the simplex $L$ provides a map
$$
a_i: {\cal  Z}^{q}_{p}(L) \longrightarrow  {\cal  Z}^{q}_{p-1}(\widehat L_i), 
\qquad 0 \leq i \leq p+q 
$$

 Let $l_j$ be the vertex opposite to the face $L_j$. 
Consider an open part ${\cal  Z}^{q}_{p}(L)^0$ of ${\cal  Z}^{q}_{p}(L)$ 
parametrising the cycles $C$ such that projection with the center $l_j$ sends $C$ to a codimension $q-1$ cycle. Then 
projection with the center at  the vertex $l_j$ of $L$ 
 defines a  map 
$$
b_j:    {\cal  Z}^{q}_{p}(L)^0 \longrightarrow 
{\cal Z}^{q-1}_{p}(\widehat L_j), \qquad 0 \leq i \leq p+q 
$$

\begin{theorem-construction} \label{chw}
For given  $q \geq 0$ there is an          
explicitly constructed chain of     $(q-p-1)$-distributions  $\omega_{p}^q =
\omega_{p}^q(L; H)$ on ${\cal Z}^{q}_{p}(L)$ such that

\begin{equation} \label{i}
i)\qquad d \omega^q_{0}(L,H) = \pi_q( \Omega_L)
\end{equation}                                                  
\begin{equation} \label{ii}                                    
ii)\quad d\omega_{p}^q(L; H) = \sum_{i=0}^{p+q}(-1)^ia_i^{\ast}\omega_{p-1}^q(L; H_i)
\end{equation}
\begin{equation} \label{iii}                                     
iii)\qquad \sum_{j=0}^{p+q+1}(-1)^j b_j^{\ast}\omega_{p}^q(L; H) = 0
\end{equation}
The restriction of  $\omega_{p}^q$ to the subvariety $\widehat 
{\cal Z}^{q}_{p}(L)$ of smooth cycles
in generic position with respect to the simplex $L$ is a real-analytic
differential $(q-p-1)$ form.   
\end{theorem-construction}

For a given positive integer 
$q$ the collection  $\{\omega_{p}^q\}$ is called the {\it $q$-th Chow
polylogarithm}.   

The varieties $ {\cal  Z}^{q}_{p}(L)$ for $ p\geq 0$ form a truncated simplicial variety 
$ {\cal  Z}^{q}_{\bullet}(L)$. The conditions i) and ii) just mean that 
the  sequence of forms $\omega_{p}^q$ is 
 a $2q$-cocycle in
the  complex computing the Deligne cohomology $H^{2q}(
{\cal Z}^{q}_{\bullet}(L), \R_{\cal D}(q))$.

{\bf Proof}. 
We define 
 $\omega^q_p$  as the Radon transform of the distribution 
$r_{p+q-1}(L;H)$ in $\C \PP^{p+q}$ over the family of cycles $Y_{\xi}$
parametrized by 
$ {\cal Z}^q_p(L)$. This means the following. Consider the incidence variety:
$$
  \Gamma_p := \{(x,\xi) \in \C\PP^{p+q} \times   {\cal Z}^q_p(L)(\C)
\quad  \mbox{such that} \quad x \in Y_{\xi}\}
$$
where $Y_{\xi}$ is the cycle in $\C \PP^{p+q}$ corresponding to $\xi
\in   {\cal  Z}^q_p(L)$. We get a double bundle
$$
\begin{array}{ccccc}
&&  \Gamma_p \subset \C \PP^{p+q} \times {\cal  Z}^q_p(L)(\C)&&\\
&&&&\\
&\pi_1\swarrow&&\searrow \pi_2&\\
&&&&\\
\C \PP^{p+q}&&&&  {\cal  Z}^q_p(L)(\C)
\end{array}
$$
Then
$$
\omega_p^q:= {\pi_2}_{\ast} {\rm Res}_{\Gamma_p} \pi_1^{\ast} (2\pi i)^{-q}r_{p+q}(L;H) 
$$
Observe that $r_{p+q}(L;H)$ is a distribution on $\C \PP^{p+q}$, 
and hence 
$\pi_1^{\ast} r_{p+q}(L;H)$ is a distribution on 
$\C \PP^{p+q} \times {\cal  Z}^q_p(L)(\C)$. The fact that this 
distribution can be restricted to 
$\Gamma_p$ is a version of Theorem \ref{8.6.02.2}, and is proved in the same way.  
The push forward ${\pi_2}_{\ast}$ of this 
distribution is well defined since  $\pi_2$ is a proper map. 

The property i)  is true by the very definition.

\begin {lemma} \label{2.4}
$\sum_{j=0}^{n+1}(-1)^j b_j^{\ast}\omega^n_0 = 0$.
\end {lemma}

{\bf Proof}. Let $s(z_0,...,z_n):= z_1/z_0 \wedge ... 
\wedge z_{n}/z_0$. The lemma follows from  the identity 
$$
\sum_{j=0}^{n+1}(-1)^j s(z_0,...,\widehat  {z_j},..., z_{n+1}) = 0
$$
So we have iii). 
To check  ii) observe that the push forward ${\pi_2}_{\ast}$ of 
distributions commutes with the De Rham differential.  
The theorem is proved.

{\bf 2. Properties of the Chow polylogarithm function}. 
The  function
${\cal P}_q:= \omega_{q-1}^q$  on $ {\cal  Z}^{q}_{q-1}(L)$
 is called the {\it Chow $q$--logarithm function}. 
It satisfies
two functional 
equations:
$$
\sum_{i=0}^{2q}(-1)^ia_i^{\ast}{\cal P}_q =0 \qquad \sum_{j=0}^{2q}(-1)^j
b_j^{\ast}{\cal P}_q = 0 
$$

\begin{theorem} \label{6.10.04.1q}
The Chow polylogarithm function is invariant under the natural 
action of the torus 
$(\C^*)^{p+q}$ on ${\cal Z}^q_p(\C)$. In particular it does not 
depend on the choice of the hyperplane $H$. 
\end{theorem}

{\bf Remark}. The statements of Theorem \ref{6.10.04.1q} are no longer true 
for the  forms $\omega_{p}^q$ for
$p<q-1$. 

Here is a reformulation of Theorem \ref{6.10.04.1q}.
\begin {theorem}  
\label {0.2}
Suppose that $dimX=n$ and $f_1,...,f_{2n+1
}$ are rational functions on $X$. Then
 the integral 
$$
(2\pi i)^{1-n}\int_{X(\C)}r_{2n}(f_1,...,f_{2n+1})
$$
does not change if we multiply one of the
functions $f_i$ by a non zero constant.
\end {theorem}
 
{\bf Proof}. Multiplying, say, $f_1$ by $\lambda$ we see that the difference
between the two integrals is
\begin{equation} \label{a6}
\log|\lambda| \sum_k a_k \int_{X(\C)}  {\rm Alt}_{2n} d\log|f_2|\wedge ... \wedge
d\log|f_{2k-1}| \wedge d\arg f_{2k} \wedge ... \wedge d\arg f_{2n+1}
\end{equation}
where the $a_k$ are some rational constants (easily computable from
(\ref{1wq})). We will prove that for each $k$ the corresponding integral in this sum
is already zero. Using the identity
$$
(d\log|f_2| + i d \arg  f_2) \wedge ... \wedge (d\log|f_{2n+1}| + i d \arg
 f_{2n+1}) = 0
$$
we can rewrite the integral
$$
\int_{X(\C)}   d\arg f_2\wedge ... \wedge d\arg f_{2n+1}
$$
as a sum of similar integrals containing $d\log |f_i|$. Our statement
follows from
\begin {proposition}  \label{5.2}
\label {0.2b} Suppose that $dimX=n$. Then
$$ 
 d {\rm Alt}_{2n}\Bigl(\log|f_2|d\log|f_3|\wedge ... \wedge
d\log|f_{2k-1}| \wedge d\arg f_{2k} \wedge ... \wedge d\arg f_{2n+1}\Bigr) = 
$$
\begin{equation} \label{d6}
  {\rm Alt}_{2n} \Bigl( d\log|f_2| \wedge ... \wedge
d\log|f_{2k-1}| \wedge d\arg f_{2k} \wedge ... \wedge d\arg f_{2n+1}\Bigr)
\end{equation}
in the sense of distributions.
\end {proposition}

{\bf Proof}. Since $d di \arg f = 2\pi i \delta(f) $, the left hand
side is equal to the right hand side {\it 
plus} the following terms concentrated on the divisors $f_{2k+j} = 0$:
$$
4\pi\cdot2(n-k+1) {\rm Alt}_{2n} \Bigl( \delta(f_{2n+1})\log|f_2|d\log|f_3| \wedge ... \wedge
d\log|f_{2k-1}| \wedge d\arg f_{2k} \wedge 
... \wedge d\arg f_{2n}\Bigr)
$$
However all these additional terms  vanish thanks to the following 
proposition.

\begin {proposition}  \label{5.3}
\label {0.2a}
Suppose that $dim X = n$. Then for each $0 \leq j \leq n-1$ one has 
\begin{equation} \label{b6}
 {\rm Alt}_{2n}\Bigl(d \log|f_1|
\wedge ... \wedge d\log|f_{2j+1}| \wedge d\arg f_{2j+2}\wedge ... \wedge
d\arg f_{2n} \Bigr) =0
\end{equation}
\begin{equation} \label{c6}
 {\rm Alt}_{2n}\Bigl(d \log|f_1|
\wedge ... \wedge d\log|f_{2j}| \wedge d\arg f_{2j+1}\wedge ... \wedge
d\arg f_{2n} \Bigr) =
\end{equation}
$$
 \frac{(2n)! {n \choose j}}{{2n \choose 2j}}\cdot d \log|f_1|
\wedge ... \wedge d\log|f_{2n}| 
$$
\end {proposition}

{\bf Proof}. The idea is this. 
One has $n$ equations 
$$
d\log f_{1} \wedge ... \wedge d\log f_{2n} =0
$$
$$
d\log |f_{1}| \wedge d\log f_{2} \wedge ... \wedge d\log f_{2n} =0
$$
\begin{equation} 
............
\end{equation}
$$
d\log |f_{1}|\wedge ... \wedge d\log |f_{n-1}| \wedge d\log f_{n}
\wedge ... \wedge d\log f_{2n} = 0
$$
Taking imaginary part of each of them and alternating $f_1,...,f_{2n}$ we
get $n$ linear equations. Solving them we get the proposition. 
See the details in the Appendix in [GZ].

\section{The Grassmannian polylogarithms}

{\bf 1. Configurations of vectors and Grassmannians: a dictionary}. 
Let $G$ be a group. Let $X$ be a $G$-set. 
We define {\it configurations of $m$ points of} $X$ as $G$-orbits in $X^m$. 

{\bf Example 1}. If $X:= V$ is a vector space and $G:= GL(V)$ we get 
configurations  of vectors
in $V$. 
A configuration of vectors $(l_1,...,l_m)$ is  in generic position 
if each $k\leq {\rm dim}V$ of
the vectors are linearly independent.

{\bf Example 2}. If $X:= \PP(V)$ is a projective space and $G:= PGL(V)$ we get 
configurations  $(x_1,...,x_m)$ of $m$ points in $\PP(V)$. 
A  configuration of points is in generic position if each $k\leq {\rm dim}V$ of
them generate a plane of dimension $k-1$. 

Let $T_{p+q}$ be the quotient of the torus ${\Bbb G}_m^{p+q+1}$ 
by the diagonal subgroup ${\Bbb G}_m = (t, ..., t)$. Below $V_n$ denotes a vector space of dimension $n$. 

\begin{lemma-construction} \label{6.12.02.5}
i) There are canonical isomorphisms  between the following sets of geometric
objects: 

(a) Configurations of $p+q+1$ vectors in generic position in $V_q$.

(b) Isomorphism classes of triples 
$\{$a projective space ${\PP}^{p+q}$ together with a  
simplex $L$, an ``infinite'' hyperplane $H$ in generic position to $L$, and a $p$-dimensional plane in generic position to $L$ (but not necessarily to $H$) $\}$. 

(c) Isomorphism classes of triples 
$\{$a vector space $V_{p+q+1}$, a basis $(e_0,...,e_{p+q})$ of 
$V_{p+q+1}$, and a $p+1$-dimensional subspace  of $V_{p+q+1}$ in generic position with respect to the coordinate hyperplanes$\}$.

ii) The torus $T_{p+q}$ acts naturally, and without fixed points,  
on each of the objects 
a), b), c), and the isomorphisms above are compatible with this action. 
\end{lemma-construction}

{\bf Proof}. i) (a) $ \rightarrow$ (c). 
For the $(p+1)$-dimensional  subspace, take  the kernel of the linear map from
$V_{p+q+1}$ to $V_q$ sending $e_i$ to $l_i$. 

(c) $ \rightarrow$ (a). Take the quotient of $V_{p+q+1}/h$ along the given subspace $h$ and consider the images of the vectors $(e_0,...,e_{p+q})$ there. 

(c) $ \rightarrow$ (b). Let $\PP^{p+q}:= \PP(V_{p+q+1})$. Let 
 ${\AAA}^{p+q}$ be the affine hyperplane in $V_{p+q+1}$
passing through the ends of the basis vectors $e_i$. 
Then ${\AAA}^{p+q} \subset \PP^{p+q}$.
 The coordinate hyperplanes in $V_{p+q+1}$ provide a simplex 
$L_{p+q} \subset {\PP}^{p+q}$.
The projectivization of a generic $(p+1)$-dimensional 
subspace $h$ in $V_{p+q+1}$ gives a
$p$-plane $\overline h$ in generic position with respect to this simplex. 
(Notice that we do not impose any condition on the mutual 
location of $H$ and $\overline h$. For instance $\overline h$ may be inside of $H$.)

(b) $\rightarrow$ (c). The triple $(\PP^{p+q}, H, L)$ provides a unique up to an isomorphism data 
$(V_{p+q+1}, (e_0,...,e_{p+q}))$. Namely, the partial data $(\PP^{p+q}, H)$ 
provides us with $(V_{p+q+1}, \widetilde  H)$ where $\widetilde  H$ is the subspace of 
$V_{p+q+1}$ projecting to $H$. Now the vertices $l_i$ of the simplex $L$ provide coordinate lines $\widetilde  l_i$ in 
$V_{p+q+1}$. Intersecting these coordinate lines with a parallel shift of the subspace $\widetilde  H$ we get a point  on each of the coordinate lines. By definition the endpoints of the basis vectors $e_i$ are these points. Taking the  subspace $\widetilde  h$ in $V_{p+q+1}$ projecting to a given $p$-plane $h$ in 
$\PP^{p+q}$ we get the desired correspondence.

 ii) The torus $T_{p+q}$ acts on the 
configurations of vectors in a) 
as $$ 
(t_1, ..., t_{p+q+1}):  
(l_1, ..., l_{p+q+1}) \lms (t_1l_1, ..., t_{p+q+1} l_{p+q+1})
$$
The torus $T_{p+q}$ is identified with $\PP^{p+q} - L$ in b), 
and so acts naturally on the data in b). The action on the data in c) is similar. 
 The lemma is proved. 

If we use the description c)  for  the Grassmannians then 
$b_j$ is obtained by factorization along the coordinate
axis $\langle e_j \rangle$.

{\bf 2. The Grassmannian and bi-Grassmannian polylogarithms}. 
Let us fix a positive integer $q$. 
The operations $a_i$ and $b_j$ from Section 3.1 transform planes to planes. So we get
the following diagram of varieties called the 
 bi-Grassmannian $\widehat  G(q)$:
$$
\begin{array}{cccccccccccc}
&&&&&&\downarrow...\downarrow&&\downarrow...\downarrow\\
&&&&&&&&&\\
\widehat G(q):=&&&&&\stackrel{\rightarrow}{\stackrel{...}{\rightarrow}}&\widehat G_1^{q+1}&\stackrel{\rightarrow}{\stackrel{...}{\rightarrow}}&\widehat  G_0^{q+1}\\
&&&&&&&&&\\
&&&&\downarrow...\downarrow&&\downarrow...\downarrow&&\downarrow...\downarrow\\
&&&&&&&&&\\
&&...&\stackrel{\rightarrow}{\stackrel{...}{\rightarrow}}&\widehat  G_2^{q}&\stackrel{\rightarrow}{\stackrel{...}{\rightarrow}}&\widehat  G_1^{q}&\stackrel{\rightarrow}{\stackrel{...}{\rightarrow}}&\widehat  G_0^{q}\\
\end{array}
$$
Here the  horizontal arrows are the maps $a_i$ and the vertical ones
are $b_j$.

{\bf Remark}. The bi-Grassmannian $\widehat  G(n)$ is not a (semi)bisimplicial
scheme. (It is a truncated semi{\it hypersimplicial} scheme. 
See s.2.6 in [G4]).

{\it Configurations of  hyperplanes  and torus quotients of Grassmannians}. 
Let $\widehat  G_p^q$ be the 
Grassmannian of $p$-planes in ${\Bbb  P}^{p+q}$ 
in generic position with respect to a given simplex $L$.

Taking the $T_{p+q}$--orbits of the objects 
a) and b) in  the lemma  we arrive at

\begin{corollary} \label{6.13.02.1}
There is a 
 bijective correspondence
$$
\widehat G^q_p/T_{p+q} 
<->  \left \{ \mbox{Configurations of $ p+q+1$ generic hyperplanes in } \PP^{p}\right \}
$$
sending a $p$--plane $h$ to the configuration  $
(h \cap L_0,...,h \cap L_{p+q})$      in $h$.
  \end{corollary} 

Let $\psi^q_{p}(q)$ be the  restriction of the differential form
$\omega_{p}^q$ to $\widehat  G^{q}_p$. The properties i), ii) from Theorem \ref{chw} 
are exactly the defining conditions for the single--valued 
Grassmannian polylogarithm whose existence was conjectured in 
 [HM], [BMS], see also [GGL]. 

Let us extend these forms by zero 
to the other rows of the bi-Grassmannian $\widehat G(q)$, i.e.   
set $\psi^{q+i}_{p}(q) =0$ if $i>0$.   Then the 
property iii) from Theorem \ref{chw} guarantees that the forms $\psi^{q+i}_{p}(q)$ 
form  a  $2q$-cocycle in
the  bicomplex computing the Deligne cohomology $H^{2q}(\widehat 
G(q)_{\bullet}, \R(q)_{\cal D})$. It is called {\it the
bi-Grassmannian $q$-logarithm}. 
( [G5]).

A sequence of {\it  
multivalued analytic } forms on Grassmannians satisfying conditions
 similar to  i),  ii) was defined  in [HaM1], [HaM2].
Another construction of the multivalued analytic Grassmannian 
polylogarithms was suggested in [G5] in the more general setting 
of the multivalued Chow polylogarithm.

{\bf 3. The Grassmannian $n$-logarithm function}. 
By  Theorem \ref{0.2} the Chow polylogarithm function is invariant under the action 
of the torus $(\C^*)^{2n-1}$. So restricting it  to the open Grassmannian $\widehat  G^n_{n-1} \subset 
\widehat  {\cal Z}^n_{n-1}$ and using  the bijection
$$
\mbox{ $\{ (n-1)$--planes in 
 $\PP^{2n-1}$ in generic position with respect to a simplex $L$\}}/({\Bbb G}_m^*)^{2n-1} \quad 
$$
$$  <--> \quad \left \{ \mbox{Configurations of} \quad 2n \quad\mbox {generic 
hyperplanes in } \PP^{n-1}\right \}
$$
we get a function on the configurations of $2n$ hyperplanes in $\C \PP^{n-1}$, called 
the Grassmannian polylogarithm function ${\cal L}^G_n$. 

The Grassmannian polylogarithm function has the following  simple description on 
the language of configurations of hyperplanes. It is intresting that in this description we can  
work with {\it any}  configuration of $2n$ 
hyperplanes,   assuming nothing about their mutual   location. 
 
 Let $h_1,...,h_{2n}$ be $2n$ arbitrary hyperplanes in $\C \PP^{n-1}$. 
Choose an additional hyperplane $h_0$. Let $f_i$ be a
rational function on $\C \PP^{n-1}$ with divisor $(h_i) - (h_0)$. It is 
defined up to a scalar
factor. 
Set
$$
{\cal L}^G_n(h_1,...,h_{2n}):= (2\pi i)^{1-n}\int_{\C \PP^{n-1}}r_{2n-2}
(\sum_{j=1}^{2n}(-1)^j f_1\wedge ...
\wedge \widehat  f_j \wedge ... \wedge f_{2n})
$$
It is skewsymmetric by definition.
Notice that 
$$
\sum_{j=1}^{2n}(-1)^j f_1\wedge ...
\wedge \widehat  f_j \wedge ... \wedge f_{2n} = \frac{f_1}{f_{2n}} \wedge
\frac{f_2}{f_{2n}} \wedge ... \wedge \frac{f_{2n-1}}{f_{2n}} 
$$
So we can define ${\cal L}^G_n(h_1,...,h_{2n})$ as follows: choose rational
functions $g_1,...,g_{2n-1}$ such that ${\rm div} g_i = (h_i) - (h_{2n})$ and put
$$
{\cal L}^G_n(h_1,...,h_{2n}) = (2\pi i)^{1-n}
\int_{\C \PP^{n-1}}r_{2n-2}(g_1,...,g_{2n-1})
$$

{\bf Remark}. The function ${\cal L}^G_n$ is defined on the set of all configurations of $2n$ hyperplanes in $\C\PP^{n-1}$. However it is 
not even continuous on this set. It is real analytic on the submanifold of generic 
configurations. Since we put no restrictions on the hyperplanes $h_i$ 
the following theorem is stronger than Theorem \ref{chw} in the case of 
linear subvarieties. 

\begin {theorem}  
\label {0.3a} The function ${\cal L}^G_n$ has  the following properties:

a) It does not depend on the choice of hyperplane $h_0$.

b) For any $2n+1$ hyperplanes $(h_1, ..., h_{2n+1})$ in $\C \PP^{n}$ one has
\begin{equation} \label{4}
 \sum_{j=1}^{2n+1}(-1)^j {\cal L}^G_n(h_j \cap h_1,..., \widehat {h_j \cap h_j}, ..., h_j \cap h_{2n+1}) =0
\end{equation}

c) For any $2n+1$ hyperplanes $(h_1, ..., h_{2n+1})$ in $\C \PP^{n-1}$ one has
\begin{equation} \label{5}
 \sum_{j=1}^{2n+1}(-1)^j {\cal L}^G_n(h_1,...,\widehat  h_j,...,h_{2n+1})=0
\end{equation}
\end {theorem}

{\bf Proof}. a) Choose another hyperplane $ h_0'$. Take a rational
function $f_0$ with divisor $(h_0') - (h_0)$. Set $f_i' = \frac{f_i}{f_0}$. Then
$$
\sum_{j=1}^{2n+1}(-1)^j f_1\wedge ...
\wedge \widehat  f_j \wedge ... \wedge f_{2n+1} - \sum_{j=1}^{2n+1}(-1)^j
 f_1' \wedge ...
\wedge \widehat  { f_j'} \wedge ..., \wedge  f_{2n+1}' = 0
$$

Indeed, substituting $f_i' = \frac{f_i}{f_0}$ in this formula we find
that the only possible 
nontrivial term $f_0 \wedge f_1 \wedge ... \wedge \widehat  f_i \wedge ...
\wedge \widehat  f_j \wedge ... \wedge f_{2n}$ vanishes because it is symmetric
in $i,j$.

b) Let $g_1,...,g_{2n+1}$ be rational functions on $\C \PP^{n}$ with
${\rm div} g_i = (h_i) - (h_0)$. 
Then
\begin{equation} \label{8}
 dr_{2n-1}\Bigl(\sum_{j=1}^{2n+1}(-1)^j g_1 \wedge ...
\wedge \widehat  g_j \wedge ... \wedge g_{2n+1}\Bigr) = 
\end{equation}
$$
\sum_{j \not = i}(-1)^{j+i-1} 2 \pi i \delta(f_i) \wedge r_{2n-2}\Bigl(
g_1\wedge ... \wedge \widehat  g_i \wedge ...
\wedge \widehat  g_j \wedge ... \wedge g_{2n+1}\Bigr)
$$
(Notice that $d\log g_1 \wedge ... \wedge \widehat  {d\log g_j} \wedge ...
\wedge d\log g_{2n+1} =0$ on $\C \PP^{n}$). 
Integrating (\ref{8}) over $\C \PP^{n}$ we see that the left hand side equals
 zero, while the right hand side equals to the sum of the expressions staying 
on the left of (\ref{4}). So we get b).

c) is obvious: we apply $r_{2n-2}$ to the zero element.
Theorem is proved.

{\bf 4. ${\Bbb P}^1 - \{0, \infty\}$ as a 
special stratum in the configuration space of $2n$ points in $\PP^{n-1}$}. 
A {\it special configuration}  is 
a configuration of $2n$  points  
\begin{equation} \label{--00}
(l_0,...,l_{n-1}, m_0,...,m_{n-1}) 
\end{equation} 
in $\PP^{n-1}$ such that  $l_0,...,l_{n-1}$ are vertices of a 
simplex in $\PP^{n-1}$ and $ m_i$ is a point 
on the edge $l_il_{i+1}$ of the simplex different from $l_i$ and 
$l_{i+1}$,  as on the picture. 

\begin{figure}[ht]
\centerline{\epsfbox{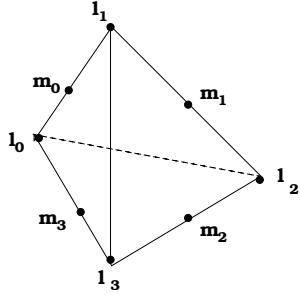}}
\caption{A special configuration of $8$ points in $\PP^{3}$.}
\label{fig1coc}
\end{figure}

\begin{proposition} \label{6.12.02.1}
The set of special configurations of $2n$ points in $\PP^{n-1}$
is canonically identified with $\PP^1 \backslash \{0, \infty \}$. 
\end{proposition}

\begin{figure}[ht]
\centerline{\epsfbox{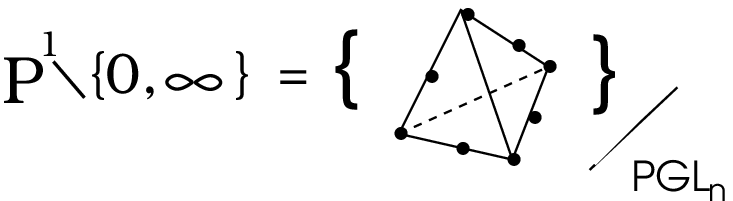}}
\caption{$\PP^1 - \{0, \infty\}$ is a stratum in the configuration space of $2n$ points in $\PP^{n-1}$.}
\label{gpol12}
\end{figure}

{\bf Proof}. 
We define {\it generalized cross-ratio} 
$$
r(l_0,...,l_{n-1},m_0,...,m_{n-1})  \in  F^*,
$$ 
where $F$ is the common field of definition of the points $l_i, m_j$, as follows. 
Consider the one-dimensional subspaces $L_i, M_j$ in the $n$-dimensional vector space $V$ projecting  to the points $l_i,m_j$ in $\PP^{n-1}$ respectively. The 
subspaces $L_i, M_i, L_{i+1}$ generate a two dimensional subspace. Its  quotient
 along   $M_i$ can be identified with $L_i$ as well as with $L_{i+1}$. 
So we get a canonical linear map 
$
\overline M_i: L_i   \longrightarrow L_{i+1}
$. 
The composition of these maps (the ``linear monodromy'')
$$
\overline M_0 \circ    ... \circ \overline M_{n-1} : L_0 \longrightarrow L_0
$$
 is   multiplication by  an element of $F^*$ 
called the generalized cross-ratio of the special configuration (\ref{--00}).

\begin{figure}[ht]
\centerline{\epsfbox{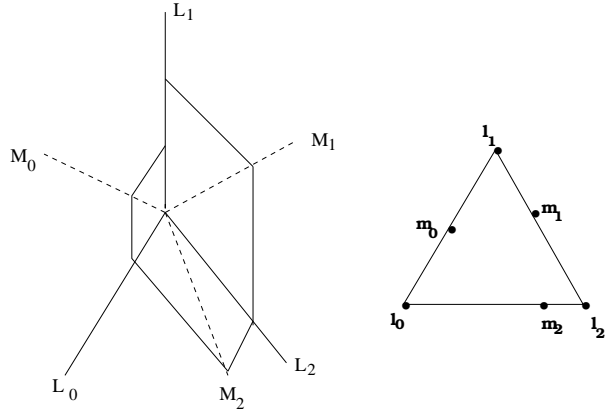}}
\caption{The generalized cross-ratio of a special configuration}
\label{gpol7}
\end{figure}

It is clearly invariant under the cyclic permutation  
$$
l_0 \to l_0 \to ... \to l_{n-1}\to l_0; \quad m_0 \to m_1 \to ... \to m_{n-1}\to m_0
$$
  Notice that $r(l_0,...,l_{n-1},m_0,...,m_{n-1}) =1$
if and only if the points $m_0,...,m_{n-1}$ belong to a hyperplane.

Let $\widehat  m_i$
be the point of intersection of the line $l_il_{i+1}$ with the
hyperplane passing through all the points $m_j$ except $m_i$. 
Let $r(x_1, ..., x_4)$ be the cross-ratio of the four points $x_i$ on $\PP^1$. Then 
$$
r(l_0,...,l_{n-1},m_0,...,m_{n-1}) =
r(l_i,l_{i+1},m_i,\widehat  m_{i+1})
$$

{\it The special configurations and classical polylogarithms}. 
Consider the configuration of $2n$  hyperplanes in $\PP^{n-1}$ given by the following 
equations in homogeneous coordinates $z_0: ... :z_{n-1}$
$$
z_0 = 0, \quad ..., \quad z_{n-1} = 0, \quad z_0 = z_1, \quad z_1 + z_2 = z_0,
$$
\begin{equation} \label{ewa99}
 z_2-z_3=0,\quad  ..., \quad z_{n-2} - z_{n-1} = 0, \quad z_{n-1} = az_0
\end{equation}
It admits the following interpretation. Recall that the classical
polylogarithm function 
$Li_{n-1}(z)$ can be defined by an iterated integral:
$$
Li_{n-1}(a) = 
\int_0^a\frac{dt}{1-t} \circ \frac{dt}{t} \circ ... \circ \frac{dt}{t}  = 
\int_{\Delta_a}\frac{dz_1}{z_1} \wedge ... \wedge \frac{dz_{n-1}}{z_{n-1}}
$$
If $a \in (0,1]$, then  the simplex $\Delta_a$ is defined by the equations
$$
\Delta_a:= \quad \{(z_1,...,z_{n-1}) \in \R^{n-1}| \quad 
0 \leq 1-z_1 \leq z_2 \leq z_3 \leq ... \leq z_{n-1} \leq a\} 
$$
The faces of the simplex $\Delta_a$ can be defined for arbitrary $a$. 
Then the codimension one faces $\{z_i =0\}$ of the coordinate 
simplex 
and the codimension one faces of the simplex $\Delta_a$ 
form the configuration (\ref{ewa99}).

We can  
reorder hyperplanes of this configuration as follows:
$$
z_0 =0, \quad z_1 =0, \quad z_1 =z_0, \quad z_1+z_2 = z_0, \quad z_2 =0, 
$$
$$
z_2 = z_3, \quad z_3 =0, \quad ...\quad , \quad z_{n-2} =0, \quad 
z_{n-2} = z_{n-1}, \quad z_{n-1} = a z_0
$$
 Applying the projective duality to this configuration of hyperplanes  
we get the special configuration of $2n$ points in $\PP^{n-1}$ 
with the generalized cross ratio $a$. 

\begin{figure}[ht]
\centerline{\epsfbox{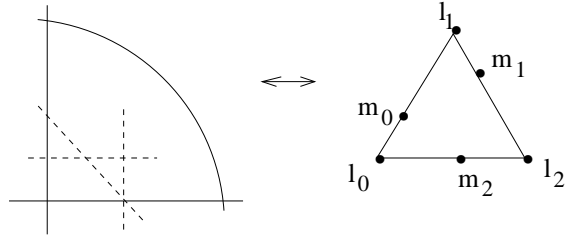}}
\caption{Classical polylogarithm configurations and special configurations}
\label{gpol11.fig}
\end{figure}

The correspondence between the configuration (\ref{ewa99}) and the 
special configuration of points is illustrated in the case $n=3$ 
on Figure \ref{gpol11.fig}.

{\bf Remark}. It is amusing that the special configuration of 
$2n$ points in $\PP^{n-1}$, 
which is related to the classical $n$-logarithm  by Theorem  \ref{clpoly} below,
 is constructed using the geometry of the mixed motive corresponding to $Li_{n-1}(a)$.

{\bf 5. Restriction of the Grassmannian $n$--logarithm to the special stratum}. 
The function $Li_n(z)$ has a  remarkable  single-valued
version  ([Z1], [BD]): 
\begin{eqnarray*} 
{\cal L}_{n}(z) &:=& \begin{array}{ll} 
{\rm Re} & (n:\ {\rm odd}) \\ 
{\rm Im} & (n: \ {\rm even}) \end{array} 
\left( \sum^{n-1}_{k=0} \beta_k
\log^{k}\vert z\vert \cdot Li_{n-k}(z)\right)\; , \quad n\geq 
2 \\ 
\end{eqnarray*}           
It is continuous on $\C \PP^1$. Here   $\frac{2x}{e^{2x} -1}  =
\sum_{k=0}^{\infty}\beta_k x^k $, so $\beta_k = \frac{2^kB_k}{k!}$ 
where the $B_k$ are the Bernoulli numbers. For example ${\cal L}_2(z)$
is the Bloch - Wigner function.

Let us consider the following modification of the 
 function ${\cal L}_n(z)$ proposed by A. M. Levin in [Le]:
$$
\widetilde {\cal L}_{n}(x):= 
$$
$$
\frac{(2n-3)}{(2n-2)}
\sum_{\mbox{$k$ even;  $0 \leq k \leq n-2$}}
\frac{2^k (n-2)!(2n-k-3)!}{(2n-3)!(k+1)!(n-k-2)!} {\cal L}_{n-k}(x)\log^k|x|
$$
For example $\widetilde {\cal L}_{n}(x) =  {\cal L}_{n}(x)$ for $n \leq 3$, 
but already $\widetilde {\cal L}_{4}(x)$ is different from ${\cal L}_{4}(x)$. 
A direct integration carried out in Proposition 4.4.1 of [Le] 
shows that  
$$
-(2\pi i)^{n-1}(-1)^{(n-1)(n-2)/2}\widetilde {\cal L}_{n}(x) =
$$
$$
\int_{\C \PP^{n-1}}\log|1-z_1| \prod_{i=1}^{n-1}d\log|z_i| \wedge \prod_{i=1}^{n-2}d\log|z_{i} - z_{i+1}| 
\wedge d\log|z_{n-1} -a| 
$$
This combined with Proposition \ref{0.3} below implies 

\begin {theorem}  \label{clpoly}
The value of the function  ${\cal L}^G_n$ on the special configuration
(\ref{--00}) is equal to $$
-(-1)^{n(n-1)/2}4^{n-1}{2n-2\choose n-1}^{-1}\widetilde {\cal L}_{n}(a)
$$ where $a = r(l_0,...,l_{n-1},m_0,...,m_{n-1})  $.
\end {theorem}

Another proof in the case $n=2$ is given in Proposition 6.8.

\begin {conjecture} \label{conj}
The Chow $n$-logarithm function can be expressed by the
Grassmannian $n$-logarithm function.
\end {conjecture}

 {\bf Remark}. Suppose that an element $\sum_k\{f^{(k)}_1,...,f^{(k)}_{2n+1}\}
\in K^M_{2n+1}(\C(X))$ has zero residues at all the divisors on an
$n$-dimensional variety $X$ over $\C$. Then it defines an element 
$$
\alpha \in gr^{\gamma}_{2n+1}K_{2n+1}(X) = Ext^{2n+1}_{{\cal M}}(\Q(0)_X,\Q(2n+1)_X) 
$$
 Its direct image to the point is
an element 
$$
\pi_{\ast}(\alpha) \in gr^{\gamma}_{n}K_{2n+1}(Spec \C) = Ext^{1}_{{\cal M}}(\Q(0),\Q(n))
$$
 Applying the regulators we see that the
integral
$\sum_k \int_{X(\C)}r_{2n}(f^{(k)}_1,...,f^{(k)}_{2n+1})$ coincides, up to a factor, 
 with
the value of the Borel regulator map on $\pi_{\ast}(\alpha)$ and so  by
results of the next chapter is 
expressible by the Grassmannian $n$-logarithms. Conjecture  \ref{conj}
tells us that this should be true {\it for any} element in
$K^M_{2n+1}(\C(X))$.

\section{Grassmannian polylogarithms, symmetric 
spaces and Borel regulators}
 
{\bf 1.  The function $\psi_n$.} Let $V_n$ be a complex 
vector space of dimension $n$. 
Let
$$
\HH_n: = \left \{ \mbox{ positive definite Hermitian forms in  
   $V_n$ }\right \}/\R ^*_+ = SL_n(\C)/SU(n) 
$$
$$
= \left \{ \mbox{ positive definite Hermitian forms in
   $V_n$ with determinant } =1 \right \} 
$$
It is a symmetric space    of rank $n-1$.  For example 
$\HH_2 = {\cal H}_3$ is the hyperbolic 3-space. 
Replacing positive definite by non-negative definite Hermitian forms 
we get a compactification $\overline \HH_n$ of the symmetric space $\HH_n$.

Let $G_x$ be the subgroup of $SL_N(\C)$ stabilizing the point  $x \in \HH_n$ . 
A point $x$ defines a one-dimensional vector space $M_x$:
$$
x \in \HH_n \longmapsto M_x:= \left \{ \mbox{measures on    } \C \PP^{n-1} 
\mbox{ invariant under }  G_x\right \}
$$
   Namely, a point $x$ corresponds to a hermitian metric in $V_n$.  
This metric provides the Fubini-Studi metric on $\C \PP^{n-1} =P(V_n)$.  
Moreover there is the Fubini-Studi K\"ahler form on $\C \PP^{n-1} = P(V_n)$; 
its imaginary part is a symplectic form. Raising it to 
$(n-1)$-th power we get the Fubini-Studi volume form.
The elements of   $M_x$ are the multiples of  the Fubini-Studi volume form.

So $\HH_n$  embeds
 into the projectivization of the space of all measures in $\C \PP^{n-1}$. 
 Taking its closure we get a compactification of $\HH_n$.

 Let us choose for any point $x \in \HH_n$ an invariant measure $\mu_x \in M_x$. 
  Then, for any $y \in {\Bbb H}_n$, the ratio     $ \mu_x/\mu_y$ is a real function on $\C \PP^{n-1}$.

Let $x_0,...,x_{2n-1}$ be points  of the symmetric space $SL_n(\C)/SU(n)$.
Consider the following function
  \begin{equation} \label{1221q} 
\psi_n(x_0,...,x_{2n-1}) := \int_{\C \PP^{n-1}}
 \log | \frac{\mu_{x_1}}{\mu_{x_0}}| d\log|\frac{\mu_{x_2}}{\mu_{x_0}} | \wedge ... \wedge d\log| \frac{\mu_{x_{2n-1}}}{\mu_{x_0}} |
\end{equation}

  {\bf 2. General properties of the function $\psi_n$}.   Let us study the properties of 
integral (\ref{1221q}) in a more general situation. Let $X$ be an $m$-dimensional manifold. For any $m+2$ measures $\mu_0,...,\mu_{m+1}$
 on $ X$ such that $\frac{\mu_{i}}{\mu_j}$ are smooth functions we can construct a differential $m$-form  on $X$:
$$
\overline r_{ m}(\mu_0:...:\mu_{ m+1}) := \log | \frac{\mu_1}{\mu_0}| d\log|\frac{\mu_2}{\mu_0} | \wedge ... \wedge d\log| \frac{\mu_{ m+1}}{\mu_0} |
$$

\begin{proposition} \label{1.2.}
The integral 
\begin{equation} \label{i2}
\int_{ X}\overline r_{ m}(\mu_0:...:\mu_{ m+1}) 
\end{equation}
satisfies the following properties:

1) Skew symmetry with respect to the permutations of $\mu_i$.

2) Homogeneity:
$$
\int_{ X}\overline r_{ m}(\lambda_0  \mu_0: ... :\lambda_{ m+1}   \mu_{ m+1}) = \int_{X }\overline r_{ m}(\mu_0:...:\mu_{ m+1}) 
$$

3)Additivity: for any $ m+3$ measures $\mu_i$ on $X $ one has
$$
  \sum_{i=0}^{m+2} (-1)^i \int_{X }\overline r_{ m}(\mu_0: ... :\widehat \mu_i: ... :\mu_{ m+2}) =0
$$

4)  Let $g$ be a diffeomorphism of $X$.  Then   
$$
\int_{X }\overline r_{ m}(g^*\mu_0: ... :g^*\mu_{ m+1})  = \int_{X }\overline r_{ m}( \mu_0: ... : \mu_{ m+1})
$$
\end{proposition}

{\bf Proof}. 1). Follows from $\log f \cdot d\log g + \log g \cdot d\log f = 
d (\log f \cdot \log g)$.

2) Using 1) we may assume $\lambda_i =1$ for $i>0$. Then 
$$
\int_{ X} (\overline r_{ m}(\lambda_0\mu_0:\mu_1: ...:\mu_{ m+1}) - \overline r_{ m}(\mu_0:\mu_1 :...:\mu_{ m+1})) = 
$$
$$
-  \log |\lambda| \cdot \int_{ X} 
d(\log| \frac{\mu_2}{\mu_0}| d \log| \frac{\mu_3}{\mu_0}| \wedge ... \wedge d\log| \frac{\mu_{m+1}}{\mu_0}|) =0
$$

3) Taking into account the skewsymmetry of the integral we  have  to prove that 
\begin{equation} \label{12321}
{\rm Alt}_{(0,..., m+2)} \left \{
\log | \frac{\mu_2}{\mu_1}| d\log| \frac{\mu_3}{\mu_1}| \wedge ... \wedge d\log| \frac{\mu_{m+2}}{\mu_1} | \right \} = 0
\end{equation}
    Let us write $ \frac{\mu_i}{\mu_j} = \frac{\mu_i}{\mu_0}/ \frac{\mu_j}{\mu_0}$  and substitute it into  (\ref{12321}). Then the    terms in (\ref{12321})  where
$$
\log | \frac{\mu_2}{\mu_0}| d\log| \frac{\mu_3}{\mu_0}| \wedge ... \wedge d\log| \frac{\mu_{m+2}}{\mu_0} | 
$$
 will appear look as follows:
$$
\log | \frac{\mu_2}{\mu_1}| d\log| \frac{\mu_3}{\mu_1}| \wedge ... \wedge d\log| \frac{\mu_{m+2}}{\mu_1} | - \log | \frac{\mu_1}{\mu_2}| d\log| \frac{\mu_3}{\mu_2}| \wedge ... \wedge d\log| \frac{\mu_{m+2}}{\mu_2} | 
$$
$$
-\log | \frac{\mu_2}{\mu_0}| d\log| \frac{\mu_3}{\mu_0}| \wedge ... \wedge d\log| \frac{\mu_{m+2}}{\mu_0} | + \log | \frac{\mu_0}{\mu_2}| d\log| \frac{\mu_3}{\mu_2}| \wedge ... \wedge d\log| \frac{\mu_{m+2}}{\mu_2} | 
$$
(The first two terms comes from ${\rm Alt}_{(1,..., m+2)}\overline r_m(\mu_1:  ... :\mu_{m+2})$ and the second two from ${\rm Alt}_{(0,2,..., m+2)}\overline r_m(\mu_0: \mu_2: ... :\mu_{m+2})$. The expression ${\rm Alt}_{ m+2}\overline r_m(\mu_0: ... :\widehat  \mu_i: ... :\mu_{m+2})$ provides no such terms if $i>1$). 

4) Clear. The proposition is proved.

Recall the following general construction. Let $G$ be a group. 
Let $X$ be a $G$-set and $f$ a function on $X^n$ satisfying 
$$
\sum_{i=1}^{n+1}(-1)^{i}f(x_1, ..., \widehat x_i, ..., x_{n+1}) =0
$$
 Choose a point $x \in X$. Then there is an $(n-1)$-cocycle of the group $G$:
$$
f_x(g_1, ..., g_n):= f(g_1 x, ..., g_n x)
$$  

\begin{lemma} \label{point}
The cohomology class of the cocycle $f_x$ does not depend on $x$.
\end{lemma}

{\bf Proof}. The difference $f_y - f_x$ is the coboundary of the $(n-2)$-chain 
\begin{equation} \label{1.22.1}
h_{x,y}(g_1, ..., g_{n-1}) = 
\sum_{i=1}^{n-1}(-1)^{k-1} f(g_1 x, g_2 x, ..., g_kx, g_ky, g_{k+1}y, 
..., g_{n-1} y)
\end{equation}
Here is the geometric picture leading to this formula. Consider the prism 
$\Delta_{g_1, ..., g_{n}}^{(n-1)}\times \Delta_{x,y}^{(1)}$ given by product 
of the $(n-1)$-simplex with vertices 
$g_1, ..., g_n$ by the  $1$-simplex with vertices $(x,y)$. Decomposing its side face  
$\Delta_{g_1, ...,  g_{n-1}}^{(n-2)}\times \Delta_{x,y}^{(1)}$ 
into simplices 
 we come to the right hand side of (\ref{1.22.1}).  Then the terms of the 
formula $f_y - f_x - \delta h_{x,y}$ 
correspond to the boundary faces of the prism. Cutting the prism into
 simplices we see that 
the sum of the terms corresponding to the prism boundary is zero thanks 
to the cocycle relation. 
The lemma is proved. 

So, for any $x \in {\Bbb H}_n$, 
$$
(\psi_n)_x(g_0,...,g_{2n-1}):= \psi_n(g_0x,...,g_{2n-1}x) 
$$ 
is a smooth $(2n-1)$-cocycle of $GL_n(\C)$. 

{\bf Remark}. This   cocycle   is the restriction to $GL_n(\C)$ 
of  the Bott cocycle for the group of diffeomorphisms of $\C \PP^{n-1}$.

Let $h_0, \dots, h_{2n-1}$ be any hyperplanes in ${\Bbb CP}^{n-1}$. 
Recall that the Grassmannian $n$-logarithm is defined by
$$
{\cal L}_n^{G}(h_0,\dots,h_{2n-1})=(2\pi i)^{1-n}
\int_{{\Bbb CP}^{n-1}} r_{2n-1}(f_1,\dots, f_{2n-1})$$
where $f_i$ is a rational function on ${\Bbb CP}^{n-1}$ with the divisor $(h_i) -(h_0)$.

\begin{proposition} \label{0.3}
One has 
$$
{\cal L}_n^{\rm G}(h_0,\dots,h_{2n-1}) = -\frac {(-4)^{n-1}(n-1)!^2} {(2\pi i)^{n-1} (2n-2)!}
\int_{{\Bbb CP}^{n-1}}\log|f_1|\bigwedge_{j=2}^{2n-1} d\log|f_j| 
$$
\end{proposition}

{\bf Proof}. See Proposition 6.2 in [GZ].

   {\bf 3. The Grassmannian polylogarithm ${\cal L}_n^G$ as the boundary value
 of the function $\psi_n$.}  
We start from an explicit formula for the Fubini-Studi form. Let $\widehat \PP^{n-1}$ 
be the variety of all hyperplanes in $\PP^{n-1}$. Consider the incidence divisor
$$
D \subset  \widehat \PP^{n-1} \times  \PP^{n-1}  \qquad D := \{(h, x)| x \in h\}
$$
 where $h$ is a hyperplane and $x$ is a point in $\PP^{n-1}$. 

Let $(x_0:...:x_{n-1})$ be  the 
 homogeneous coordinates of a point $x$ in $\PP^{n-1}$. 
Let
$$
\sigma_n(x,dx):= \sum_{i=0}^{n-1} (-1)^i x_i dx_0 \wedge ... 
\wedge \widehat dx_i \wedge ... \wedge dx_{n-1}  = i_E {\rm vol}_x
$$
be the Leray from. Here ${\rm vol}_x = dx_0 \wedge ... 
\wedge dx_{n-1}$ and $E =\sum x_i \partial_{x_i}$. 

 There is a canonical differential $(n-1,n-1)$-form
$\omega_D$ on $ \widehat \PP^{n-1} \times  \PP^{n-1} - D$ 
with a polar singularity at the divisor $D$. Namely, 
let $x \in V_n$ and $\xi \in V_n^*$. Then 
$$
\omega_D:= \frac{1}{(2\pi i)^{n-1}}\frac{\sigma_n(\xi,d\xi) 
\wedge \sigma_n(x,dx)}{<\xi,x>^{n}}
$$
It is $PGL_{n}$-invariant. 
A hermitian metric $H$ in $V_n$ provides an isomorphism 
$H: V_n \lra \overline V_n^*$, and hence an isomorphism 
${\Bbb C} \PP^{n-1} \lra \widehat {\overline {{\Bbb C}P}}^{n-1}$. 
The graph $\Gamma_H$ of this map does not intersect the incidence divisor $D$. Thus 
restricting the form $\omega_D$ to $\Gamma_H$ we 
get a volume form on ${{\Bbb C}P}^{n-1}$:
\begin{equation} \label{FUBS}
\omega_{FS}(H) := \frac{1}{(2\pi i)^{n-1}}\frac{\sigma_n(z,dz) \wedge \sigma_n(\overline z, d \overline {z})}{H(z, \overline z)^{n}}
\end{equation} 
It is clearly invariant under the group preserving the Hermitian form $H$. Moreover, 
it is the Fubini-Studi volume form:  
a proof can be obtained by using the explicit formula for the Fubini-Studi
 K\"ahler form given  in [Ar], complement 3.

One can realize 
$\C \PP^{n-1}$ as the smallest stratum  of the boundary of $\HH_n$ . Namely, for 
a hyperplane $h$ in an $n$-dimensional complex vector space $V_n$ let 
$$
F_h:= \left \{\mbox{nonnegative definite hermitian forms in 
  $V_n$ with  kernel}  h \right \}/ \R _+^*
$$
The set of 
 hermitian forms in $V_n$ with the 
kernel  $h$  is isomorphic to $\R _+^*$,  so  $F_h$ defines a point on 
the boundary of $\overline \HH_n$. 

 For any nonzero 
nonnegative definite hermitian form $H$ one can define the corresponding 
Fubini-Studi form by formula (\ref{FUBS}). It is 
a differential form with  singularities along the projectivization of 
the kernel of $H$. In particular if $h$ is a hyperplane then the degenerate hermitian 
form $F_{h}$ provides the Lebesgue measure on the affine space $\C\PP^{n-1} -  h$. 
Indeed, if 
$h_0 = \{z_0 =0\}$ then (\ref{FUBS}) specializes to 
    $$
\frac{1}{(2\pi i)^{n-1}}  d\frac{z_1}{z_0} \wedge ... \wedge d\frac{z_{n-1}}{z_0}  \wedge  d\frac{ \overline z_1}{\overline z_0} \wedge ... \wedge d\frac{ \overline z_{n-1}}{\overline  z_0}  
$$
Denote by $M_{h}$ the one dimensional real vector space generated by this form. 
For any hyperplane $h$ in $\C \PP^{n-1}$ 
let us choose a measure $\mu_{h} \in M_{h}$.

\begin{proposition} \label{1/18.1}
For any $2n$ hyperplanes $h_0,...,h_{2n-1}$ in $\C \PP^{n-1}$ the integral 
\begin{equation} \label{1221N} 
\psi_n(h_0,...,h_{2n-1}) := \int_{\C \PP^{n-1}}
 \log | \frac{\mu_{h_1}}{\mu_{h_0}}| d\log|\frac{\mu_{h_2}}{\mu_{h_0}} | \wedge ... \wedge d\log| \frac{\mu_{h_{2n-1}}}{\mu_{h_0}} |
\end{equation} 
 is convergent and equals to 
$$
(-4)^{-n} \cdot (2\pi i)^{n-1} (2n)^{2n-1}{2n-2\choose n-1} \cdot {\cal L}_n^G(h_0, ..., h_{2n-1})
$$
\end{proposition} 

{\bf Proof}. Let  $h_1,h_2$  be hyperplanes in $\C \PP^{n-1}$ and $f$ be a rational function   
such that   $(f) = (h_1) -(h_2)$. 
From the explicit description of $M_h$ given above we immediately see that 
\begin{equation} \label{AG362}
\mu_{h_1}/\mu_{h_2} =   \lambda \cdot |f|^{2n}
\end{equation}
Using this and Theorem 2.4 we see that  
integral (\ref{1221N}) is convergent. The second statement follows from 
Proposition \ref{0.3} and (\ref{AG362}). The proposition is proved.

More generally, take any $2n$ Hermitian forms 
$H_0, ..., H_{2n-1}$, possibly degenerate. For each of the forms $H_i$ 
 consider the corresponding measure $\mu_{H_i}$ ( a multiple of the Fubini-Studi form related to $H_i$). Using the convergence of the integral (\ref{1221N}) 
we can deduce that the integral 
\begin{equation} \label{1221M} 
\psi_n(H_0,...,H_{2n-1}) := \int_{\C \PP^{n-1}}
 \log | \frac{\mu_{H_1}}{\mu_{H_0}}| 
d\log|\frac{\mu_{H_2}}{\mu_{H_0}} | \wedge ... \wedge 
d\log| \frac{\mu_{H_{2n-1}}}{\mu_{H_0}} | = 
\end{equation}
$$
-n^{2n-1}\cdot \int_{\C \PP^{n-1}}
 \log | \frac{H_1(z, \overline z)}{H_0(z, \overline z)}| d\log|\frac{H_2(z, \overline z)}{H_0(z, \overline z)} | \wedge ... \wedge d\log| \frac{H_{2n-1}(z, \overline z)}{H_{0}(z, \overline z)} |
$$
is also convergent. This enables us to extend  $\psi_n$ to the function 
$\overline \psi_n(x_0,...,x_{2n-1})$ 
on the configuration space of $2n$ points in $\overline \HH_{n-1}$. The function $\overline \psi_n$ is discontinuous. For instance it is discontinuous at the point 
$x_1 = ... = x_{2n-1} = F_h$ for a given hyperplane $h$ in $\C \PP^{n-1}$. 
It is however a smooth function on an open part of any given strata. We will keep the notation
$$
\psi_n(h_0,..., h_{2n-1}) = \overline \psi_n(F_{h_0},...,F_{h_{2n-1}})
$$

Applying Lemma \ref{point} to the case when $X$ is $\overline {\Bbb H}_n$ and using only the fact that 
the function $\overline \psi_n(x_0,...,x_{2n-1})$ is well defined for {\it any}
 $2n$ points in $\overline {\Bbb H}_n$ and satisfies the cocycle condition for any $2n+1$ of them we get 

\begin{corollary} \label{pointx1}
Let $x \in {\Bbb H}_n$ and let $h$ be a hyperplane in $\C\PP^{n-1}$. 
Then the cohomology classes of the following cocycles coincide:
$$
\psi_n(g_0 x,...,g_{2n-1}x) \quad \mbox{and} \quad \psi_n(g_0 h,...,g_{2n-1}h)
$$ 
\end{corollary}

{\bf 4. A normalization of the Borel class $b_n$}. Choose a hermitian metric 
in $V_n$. 
Let $e$ be the corresponding point of the symmetric space ${\Bbb H}_n$; its 
stabilizer is the subgroup $SU(n)$. One has 
$$
\Bigl(\Lambda^{\bullet}T^*_e{\Bbb H}_n\Bigr)^{SU(n)} = 
{\cal A}^{\bullet}(SL_n(\C)/SU(n))^{SL_n(\C)}  
$$
There are well known  
canonical ring isomorphisms (see [B2] and references there):
$$
\Bigl(\Lambda^{\bullet}T^*_e{\Bbb H}_n\Bigr)^{SU(n)}\otimes_{\R} \C = 
\Lambda^{\bullet}(sl_n(\C))^{sl_n(\C)} = 
$$
\begin{equation} \label{1.20.2}
H^{\bullet} (sl_n(\C), \C) \stackrel{\alpha}{=} H_{\rm top}^{\bullet} (SU(n), \C) 
\stackrel{\beta}{=}
H^{\bullet}_m(SL_n(\C), \C)
\end{equation}
where $H^{\bullet}$ is the Lie algebra cohomology, $H^{\bullet}_{\rm top}$ is the topological cohomology, and $H_m(G)$ denotes the measurable cohomology of a Lie group $G$. 
The first isomorphism  is obvious: $T_e{\Bbb H}_n\otimes_{\R} \C = sl_n(\C)$. 
The map 
$$
\alpha_{\rm DR}: \Lambda_{\Q}^{\bullet}(sl_n)^{sl_n} \stackrel{\sim}{\lra} 
H_{\rm DR}^{\bullet} (SL_n(\C), \Q) 
$$ 
sends an $sl_n$--invariant exterior form 
on $sl_n$  to the right--invariant one, and hence 
biinvariant differential form on $SL_n(\C)$. 
Let us describe the map 
$$
\beta_{\rm DR}: H_{\rm DR}^{\bullet} (SL_n(\C), \Q) 
\stackrel{}{\lra}
H^{\bullet}_m(SL_n(\C), \C)
$$ 

Let $C$ be a biinvariant, and hence closed, 
 differential $(2n-1)$--form on $SL_n(\C)$. Let us restrict it first 
to the Lie algebra, and then  
to the orthogonal 
complement $su(n)^{\perp}$ to the Lie subalgebra $su(n) \subset  sl_n(\C)$. 
We identify the $\R$--vector spaces $T_e{\Bbb H}_n$ and  $su(n)^{\perp}$. 
The obtained exterior form on $T_e{\Bbb H}_n$ is the restriction of 
 an invariant differential form, denoted $\omega_{C}$,  
on the symmetric space 
${\Bbb H}_n$. It is a closed differential form.

For any ordered $2n$ points $x_1,...,x_{2n}$ 
in ${\Bbb H}_n$ there is 
a  geodesic 
simplex $I(x_1,...,x_{2n})$ in ${\Bbb H}_n$. It is constructed inductively 
as follows. Let 
$I(x_1,x_2)$ be the geodesic from $x_1$ to $x_2$. The  
geodesics from $x_3$ to 
the points of  $I(x_1,x_2)$ form a geodesic triangle 
$I(x_1,x_2,x_3)$. All the geodesics from $x_4$ to the points 
of the geodesic triangle $I(x_1,x_2,x_3)$ form a geodesic simplex 
$I(x_1,x_2,x_3,x_4)$, and so on. When the 
rank of the symmetric space is greater than 1 (i.e. $n>2$) 
the geodesic simplex $I(x_1,...,x_k)$ depends on the 
ordering of the vertices $x_1,...,x_k$. 

The 
 differential $(2n-1)$-form $\omega_{C}$ on 
$SL_n(\C)/SU(n)$ provides a volume of the geodesic simplex: 
$$
{\rm vol}_C I(x_1,...,x_{2n}):= \int_{I(x_1,...,x_{2n})}\omega_{C}
$$
For every $2n+1$ points $x_1,...,x_{2n+1}$ the boundary of the simplex 
$I(x_1,...,x_{2n+1})$ is the alternating sum of the simplices 
$I(x_1,..., \widehat x_i,..., x_{2n+1})$. Since the form $\omega_{C}$ is 
closed, the Stokes theorem yields 
\begin{equation} \label{6.15.02.1}
\sum_{i=1}^{2n+1}(-1)^i \int_{I(x_1,..., \widehat x_i,..., x_{2n+1})}\omega_{C} = 
\int_{I(x_1,..., x_{2n+1})}d\omega_{C} = 0
\end{equation}
This just means that for a given point $x$ the function
$ {\rm vol_C} I(g_1 x,...,g_{2n} x)$ is a smooth $(2n-1)$-cocycle 
of the Lie group 
$SL_n(\C)$. It was considered by J.Dupont [D]. 
By Lemma \ref{point} cocycles corresponding to different points $x$ 
are canonically cohomologous. 
The obtained cohomology class is the class  $\beta_{\rm DR}([C])$.

{\bf Remark}. 
${\rm vol}_CI(x_1,...,x_{2n})$ is independent up to a sign 
of the ordering of its vertices. 
Indeed, 
consider $2n+1$ points $(x_1, x_2, x_1, x_3,...,x_{2n})$ 
and apply relation (\ref{6.15.02.1}).

{\it The Betti cohomology of $SL_n(\C)$}.  Recall  that $SU(n)$ is a retract of 
$SL_n(\C)$. 
It is well known that 
$$
H_{\rm top}^{\bullet} (SU(n), \Z) = 
H_{\rm top}^{\bullet} (S^3 \times S^5 \times ... \times S^{2n-1}, \Z) = 
\Lambda^*(B_3, 
B_5, ..., B_{2n-1})
$$
The restriction from $SU(n)$ 
to $SU(m)$ kills the classes $B_{2k-1}$ for $k>m$. If $k\leq m$ 
it identifies the class 
$B_{2k-1}$ for $SU(n)$ with the one for $SU(m)$. 
The class $B_{2n-1}$ for $SU(n)$ 
is provided by the fundamental class of the sphere $S^{2n-1} \subset \C^n$. 
Namely, it is the pull back of the 
fundamental class under the map $SU(n) \lra S^{2n-1}$ provided by a choice of a point on 
$S^{2n-1}$.  This sphere has the orientation 
induced by the one of $\C^n$. 
Thus 
\begin{equation} \label{1.20.1}
\Z \cdot B_{n} = {\rm Ker}\Bigl(H_{\rm top}^{2n-1}(SU(n), \Z) \lra 
H_{\rm top}^{2n-1}(SU(n-1), \Z)  \Bigr)
\end{equation} 
The transgression in the Leray spectral sequence for the universal $SU(n)$-bundle 
$EU(n) \lra BU(n)$ provides an isomorphism 
$$
\Z\cdot B_n \lra \frac{H^{2n}(BSU(n), \Z)}{\oplus_{0 < i < 2n}H^i \cdot H^{2n-i}}
$$
and identifies $B_{n}$ with the Chern class 
$c_n \in H^{2n}(BSU(n), \Z)$ of the associated vector bundle. 

{\it The De Rham cohomology of $SL_n(\C)$}. Consider the differential form
\begin{equation} \label{7.1.02.2}
\omega_{D_n}:= {\rm tr}(g^{-1}dg)^{2n-1} \in \Omega^{2n-1}(SL_N)
\end{equation} 
Its restriction to the subgroup $SL_m$ is zero for $m<n$. It follows that 
the cohomology class 
$$
[\omega_{D_n}] \in H^{2n-1}_{\rm DR}(SL_n, \C)
$$
is a multiple of $B_n$. 
The Hodge considerations 
show that $[\omega_{D_n}] \in  (2\pi i)^n \Q \cdot 
B_n$. 

\begin{lemma} \label{7.1.02.3} The differential form $\omega_{D_n}$ is an 
$\R(n-1)$-valued form. In particular it provides 
a cohomology  class 
$$
 b_n := \beta_{\rm DR}(\omega_{D_n}) \in H^{2n-1}_m(SL_n(\C), \R(n-1)) 
$$
\end{lemma}

{\bf Proof}. An easy calculation
shows that the value of the exterior form ${\omega_{D_n}}_{|T_e{\Bbb H}_n}$ 
on 
$$
(e_{1,n} + e_{n,1}) \wedge  i(e_{1,n} - e_{n,1}) \wedge ... \wedge  
(e_{n-1, n} + e_{n, n-1}) \wedge  i(e_{n-1, n} - e_{n, n-1}) \wedge e_{n,n} 
$$
is non zero, and obviously lies in $\Q(n-1)$. 

On the other hand the values of the form $\omega_{D_n}$ 
lie in a one dimensional $\R$--vector space. Indeed, the space of 
$su(n)$--invariant real exterior 
$(2n-1)$--forms on the space of all 
hermitian $n \times n$ matrices,
 which have  zero restriction to the subspace of hermitian 
$(n-1) \times (n-1)$ matrices,  
is one-dimensional.  The exterior form ${\omega_{D_n}}_{|T_e{\Bbb H}_n}$ 
belongs to the complexification of this space. The lemma follows from this.

We call the cohomology class provided by this lemma the Borel class, and 
use it below to construct the Borel regulator.

{\bf 5. Comparison of the Grassmannian and Borel cohomology classes of $GL_n(\C)$}. 
Let $[C_n^G]$ be the cohomology class of the $(2n-1)$--cocycle 
of $GL_n(\C)$ provided by the Grassmannian $n$--logarithm (see Corollary 5.5). 
We want to compare it with the Borel class. 

Let us consider the following integral 
$$
\widetilde C_n(H_1, ..., H_{2n-1}):=
$$
\begin{equation} \label{A362} 
-n^{2n-1}\cdot \int_{\C \PP^{n-1}}
 \frac{H_1(z, \overline z)}{(z, \overline z)} d\frac{H_2(z, \overline z)}{(z, \overline z)}  \wedge ... \wedge d\frac{H_{2n-1}(z, \overline z)}{(z, \overline z)} 
\end{equation}
where  the $H_i$ are arbitrary complex matrices and  $H_i(z, \overline z)$ 
are the bilinear form in $z, \overline z$ given by  the matrix $H_i$. 
We claim that it is a $(2n-1)$-cocycle of the Lie algebra $gl_n(\C)$, and it is 
obtained by differentiating  the 
group cocycle provided by the function (\ref{1221M}). 
We put these  facts in the following framework. 

If we restrict to the case when $H_i$ are hermitian matrices, 
integral (\ref{A362}) admits the following interpretation. 
Let us construct a  map
$$
{\Bbb M}_e: \C \PP^{n-1} \lra T^*_e{\Bbb H}_n
$$
which is a version of the moment map. 
For a point $z \in \C \PP^{n-1}$ the value of the $\langle M_e(z), v\rangle$ 
of the functional ${\Bbb M}_e(z)$ on a vector $v \in T_e{\Bbb H}_n$ is defined as follows. 
Let $e(t)$ be a path in 
${\Bbb H}_n$ such that $e(0) =e$ and $\stackrel{{\bf .}}{e}(0) =v$. 
recall the measure $\mu_x$ defined in Section 5.1. Then 
$$
<{\Bbb M}_e(z), v>:= \frac{d}{dt}\log \frac{\mu_{e(t)}(z)}{\mu_{e}(z)}|_{t=0} 
$$ 
Choose coordinates $z_1, ..., z_{n}$ in $V_n$ such that 
$(z, \overline z):= |z_1|^2 + ... + |z_{n}|^2$ corresponds to the point $e$. 
 Then $T_e{\Bbb H}_n$ is identified with the space of hermitian $(n\times n)$ 
matrices $H$. It follows from 
(\ref{FUBS}) that 
\begin{equation} \label{1/18.2} 
<{\Bbb M}_e(z), H>:= n \frac{H(z,\overline z)}{(z,\overline z)}
\end{equation}
The map ${\Bbb M}_e$ is clearly $SU(n)$-invariant. Its image is an  $SU(n)$-orbit 
in $T^*_e{\Bbb H}_n$ isomorphic to $\C \PP^{n-1}$. 

We need the following general construction. Let $V$ be a real vector space and $M$  a compact subset of  $V^*$ which is the closure of a $k$-dimensional submanifold. Any element $\omega \in \Lambda^kV$ can be viewed as  a  $k$-form 
$\omega$ on $V^*$. Integrating it over $M$ we get 
an exterior form $C_M \in \Lambda^kV^*$. If $M$ is a cone over $M'$ with the vertex at the origin then $\int_M\omega = \int_{M'}i_E\omega$ where $E$ is the Euler vector field on $V$.

 Applying this construction to the cone over the orbit $M_e$ 
based at  the origin we get an  
 $SU(n)$-invariant element 
$\widetilde C_n \in \Lambda^{2n-1}T^*_e{\Bbb H}_n$. 
It follows from (\ref{1/18.2}) that 
it is given by formula (\ref{A362}) multiplied by $2n$.

Another invariant $(2n-1)$-cocycle 
 $C_n$ of the Lie algebra $gl_n$, considered  by Dynkin [Dy],  is given by
\begin{equation} \label{7.1.02.1}
C_n(X_1, ..., X_{2n-1}) = \frac{1}{n! }{\rm Alt}_{2n-1}{\rm Tr}(X_1 X_2 ... X_{2n-1})
\end{equation}
Let $[C_n]$ be the cohomology class 
of $GL_n(\C)$ 
corresponding to the cocycle $C_n$.  

\begin{theorem} \label{1.20.10} One has 
$$
\widetilde C_n \quad = -(-1)^{\frac{(n-1)n}{2}}\frac{(2\pi i)^{n-1}n^{2n-1}(n-1)!}{(2n-1)!} 
\cdot  C_n
$$
and the class $[C_n^G]$ is a non zero rational multiple of  $[C_n]$:
$$
[C_n^G] = - (-1)^{\frac{n(n+1)}{2}}\frac{(n-1)!^3}{(2n-2)!(2n-1)!} [C_n]
$$
\end{theorem}

{\bf Proof}. The second claim follows from the first 
using Proposition \ref{1/18.1}.

Let us prove the first claim. 
The restriction of the cocycle $C_n$ to the Lie subalgebra $gl_{n-1}(\C)$ equals to zero. 
This follows, for instance,  from the Amitsur-Levitsky theorem: for any $n \times n$ 
matrices $A_1, ..., A_{2n}$ one has 
${\rm Alt}_{2n} (A_1, ..., A_{2n})=0$.

On the other hand the restriction of the cocycle $\widetilde C_n$ 
 to the Lie subalgebra of matrices 
$(a_{ij})$ where $a_{1j} = a_{j1} =0$ 
is zero. 
Indeed, in this case the form we integrate 
in (\ref{A362}) is a differential $(2n-2)$-form in $dz_2, ..., dz_{n-1}, 
d \overline z_2, ..., d\overline z_{n-1}$ and thus it is zero.  
So thanks to (\ref{1.20.2}) and (\ref{1.20.1}) we conclude that  the 
cocycle $C_n $ is proportional to $\widetilde C_n$. 
To determine the proportionality coefficient 
we compute the values of the both cocycle on a special element 
$E_n \in \Lambda^{2n-1}gl_n$. To write it down denote by $e_{i,j}$ the 
elementary $n\times n$ matrix whose only non-zero entry is $1$ on the $(i,j)$ place. Then 
\begin{equation} \label{1/16.2}
E_n:= \bigwedge_{j=1}^{n-1}(e_{j,n}\wedge  e_{n,j} )\wedge e_{n,n} 
\end{equation}
 A direct computation shows that  
$$<C_n, E_n> = 1
$$
Indeed, to get a non zero trace we have to multiply $(n-1)$ blocks 
$e_{n,j} e_{j,n}$, as well as $e_{n,n}$, which can be inserted 
anywhere between these blocks. 
So there are $(n-1)! n = n!$ possibilities. 

Let us compute the value of the cocycle $\widetilde C_n$ on $E_n$.  
\begin{lemma}
Integral (\ref{A362}) equals 
\begin{equation} \label{1/18.4} 
\frac{-n^{2n-1}}{(2n-1)!}\cdot {\rm Alt}_{2n-1}\int_{\C \PP^{n-1}}
 \frac{H_1(z, \overline z)d H_2(z, \overline z)  \wedge ... \wedge dH_{2n-1}(z, \overline z)}{(z, \overline z)^{2n-1}} 
\end{equation}
\end{lemma}

{\bf Proof}. By Proposition \ref{1.2.} integral (\ref{A362}) equals to 
\begin{equation} \label{1?1} 
\frac{-n^{2n-1}}{(2n-1)!}{\rm Alt}_{2n-1}\cdot \int_{\C \PP^{n-1}}
 \frac{H_1(z, \overline z)}{(z, \overline z)} d\frac{H_2(z, \overline z)}{(z, \overline z)}  \wedge ... \wedge d\frac{H_{2n-1}(z, \overline z)}{(z, \overline z)} 
\end{equation}
One has, for $i=2, ..., 2n-1$, that  
$$
d\frac{H_{i}(z, \overline z)}{(z, \overline z)}  = 
\frac{(z, \overline z) d H_{i}(z, \overline z) - H_{i}(z, \overline z) d (z, \overline z)}{(z, \overline z)^2}
$$
Substituting  
$$
\frac{- H_{i}(z, \overline z) d (z, \overline z)}{(z, \overline z)^2}\quad \mbox{instead of} \quad d\frac{H_i(z, \overline z)}{(z, \overline z)} \quad \mbox{in (\ref{1?1})} 
$$
we get zero since $H_1$ and $H_i$  
appear in a symmetric way and thus disappear after the alternation. 
The lemma follows.

Let us calculate integral (\ref{1?1}) in the special case 
$$
H_{2n-1}(z, \overline z) = |z_n|^2, \quad H_{2k-1}(z, \overline z) = z_k \overline z_n, \quad H_{2k}(z, \overline z) = z_n \overline z_k 
$$  
so that   $H_1 \wedge ... \wedge H_{2n-1} =E_n$. 
We will  restrict the integrand to the affine part $\{z_n=1\}$ 
and then perform the integration.  
Since $dH_{2n-1}(z, \overline z) =0$ on $\{z_n=1\}$ and 
$d z_k \wedge d\overline z_k  = -2i dx_k \wedge dy_k$ 
we get
$$
-(-1)^{\frac{(n-1)(n-2)}{2}}\frac{(-2i)^{n-1}n^{2n-1}}{2n-1}\int_{\C ^{n-1}}\frac{
 d^{n-1}x d^{n-1} y}
{(1+ |z_1|^2 + ... + |z_{n-1}|^2)^{2n-1}} =
$$
$$
-(-1)^{\frac{(n-1)n}{2}}\frac{(2i)^{n-1}n^{2n-1}}{2n-1} {\rm vol}(S^{2n-3})
\int_0^{\infty}\frac{r^{2n-3}dr}{(1+r^2)^{2n-1}} =
$$
$$
 -(-1)^{\frac{(n-1)n}{2}}\frac{(2\pi i)^{n-1}n^{2n-1}}{(n-1)!(2n-1)} 
{2n-2 \choose n-1}^{-1}
$$
since the volume of the sphere $S^{2n-3}$ is 
$\frac{2\pi^{n-1}}{(n-2)!}$ and 
\begin{equation} \label{1/18.3}
\int_0^{\infty}\frac{r^{2n-3}dr}{(1+r^2)^{2n-1}} \quad = \quad \frac{1}{2}\int_0^{\infty}\frac{r^{n-2}dr}{(1+r)^{2n-1}} \quad = \quad \frac{1}{2n-2} {2n-2 \choose n-1}^{-1}
\end{equation}
To get the last equality we integrate by parts: 
$$
\int_0^{\infty}\frac{r^{n-2}dr}{(1+r)^{2n-1}} = - \frac{1}{2n-2}\int_0^{\infty}r^{n-2} \left(
\frac{1}{(1+r)^{2n-2}}\right)'dr = 
$$
$$
\frac{n-2}{2n-2}\int_0^{\infty}r^{n-3}\frac{dr}{(1+r)^{2n-2}} = ... = 
\frac{(n-2)!n!}{(2n-2)!}\int_0^{\infty}\frac{dr}{(1+r)^{n+1}} = 
\frac{1}{n-1} {2n-2 \choose n-1}^{-1}
$$
Theorem \ref{1.20.10} is proved.

{\bf 6. Construction of the Borel regulator 
via Grassmannian polylogarithms}. 
Let $G$ be a group. The diagonal map $\Delta: G \lra G \times G$ provides a homomorphism $\Delta_*: H_n(G) \lra H_n(G \times G)$. 
Recall that 
$$
{\rm Prim}H_n(G) := \{x \in H_n(G)| \Delta_*(x) = x \otimes 1 + 1 \otimes x\}
$$
Set $A_{\Q}:= A \otimes \Q$. One has 
$$
K_n(F)_{\Q} = {\rm Prim}H_n(GL(F))_{\Q} = {\rm Prim}H_nGL_n(F)_{\Q}
$$
where the second isomorphism is provided by Suslin's stabilization theorem. 

The Borel regulator is a map 
$$
r_{n}^{\rm Bo}: K_{2n-1}(\C)_{\Q} \lra \R(n-1)
$$
provided by pairing the class $b_n \in H^{2n-1}(GL_{2n-1}(\C), \R(n-1))$ with 
the subspace $K_{2n-1}(\C)_{\Q} \subset H_{2n-1}(GL_{2n-1}(\C), {\Q})$. 

Recall the Grassmannian complex $C_*(n)$
$$
... \stackrel{d}{\lra} C_{2n-1}(n) \stackrel{d}{\lra}  C_{2n-2}(n) 
\stackrel{d}{\lra} ... \stackrel{d}{\lra} C_{0}(n)
$$
where $C_k(n)$ is the free abelian group generated by configurations 
of $k+1$ 
vectors $(l_0, ..., l_{k})$ in generic position in an $n$--dimensional 
vector space over a field $F$, and $d$ is given by the standard 
formula 
(see s. 3.1 in [G2]). The group $C_k(n)$ is in degree $k$. 
Since it is a homological resolution of the trivial $GL_n(F)$--module $\Z$  
(see Lemma 3.1 in [G2]), 
there is a canonical homomorphism 
$$
\varphi_{2n-1}^n: H_{2n-1}(GL_n(F)) \lra H_{2n-1}(C_*(n)) 
$$
Thanks to Lemma \ref{7.1.02.3} 
the Grassmannian $n$--logarithm function provides a homomorphism 
\begin{equation} \label{6.24.02.3}
{\cal L}_n^G: C_{2n-1}(n) \lra \R(n-1); \quad (l_0, ..., l_{2n-1}) \lms 
{\cal L}^G_n(l_0, ..., l_{2n-1})
\end{equation}
Thanks to the first $(2n+1)$--term functional equation for ${\cal L}_n^G$, 
see (\ref{4}),  
it is zero on the subgroup $dC_{2n}(n)$. So it induces a homomorphism
$$
{\cal L}_n^G: H_{2n-1}(C_{*}(n)) \lra \R(n-1);
$$ 
\begin{lemma} \label{6.24.02.1} The composition ${\cal L}_n^G\circ 
\varphi_{2n-1}^n$ coincides with the class $[C_n^G]$.
\end{lemma}

{\bf Proof}. Standard, see [G4].

To construct the Borel regulator we extend, as in s. 3.10 of [G2], 
the class $[C_n^G]$ 
to a class of $GL_{2n-1}(\C)$. Let us recall the key steps. 

Let $\Z[S]$be the free abelian group generated 
by a set $S$. Let $F$ be a field. 
Applying the covariant functor $\Z \lms \Z[X(F)]$ 
to the bi-Grassmannian 
$\widehat G(n)$ (see s. 4.2), and taking the alternating sum of 
the obtained homomorphisms, 
 we get a bicomplex. 
Using Lemma \ref{6.12.02.5} we see that 
it looks as follows ([G2], s. 3.7):
$$
\begin{array}{cccccccc}
&&&&&... & \stackrel{d}{\lra}&C_{2n-1}(2n-1)\\
&&& &&&&\downarrow \\
&&&...&&...&&...\\
&&& \downarrow &&&&\downarrow \\
& ... & \stackrel{d}{\lra} & C_{2n-1}(n+1) & \stackrel{d}{\lra} & ... & 
\stackrel{d}{\lra} & C_{n+1}(n)\\
&\downarrow && \downarrow &&&&\downarrow \\
... \stackrel{d}{\lra} & C_{2n-1}(n) & \stackrel{d}{\lra} & C_{2n-2}(n) & 
\stackrel{d}{\lra}& ... &\stackrel{d}{\lra} & C_n(n) 
\end{array}
$$
In particular the bottom row is the stupid truncation of the 
Grassmannian complex at
the group $C_n(n)$. 
The total complex of this bicomplex is a homological complex, 
 called the weight $n$ bi--Grassmannian complex $BC_*(n)$.  
In particular there is a homomorphism 
\begin{equation} \label{6.24.02.2}
H_{2n-1}(C_*(n)) \lra H_{2n-1}(BC_*(n)) 
\end{equation} 

In [G1-2] we proved that there are homomorphisms 
$$
\varphi_{2n-1}^m:  H_{2n-1}(GL_m(F)) \lra H_{2n-1}(BC_*(n)), 
\quad m \geq n
$$ 
whose restriction to the subgroup $GL_n(F)$ coincides 
with the composition 
$$
H_{2n-1}(GL_n(F)) \stackrel{{\varphi_{2n-1}^n}}{\lra} H_{2n-1}(C_*(n)) 
\stackrel{(\ref{6.24.02.2})}{\lra} H_{2n-1}(BC_*(n)), 
$$ 
Let us  extend homomorphism 
(\ref{6.24.02.3}) to a homomorphism
$$
{\cal L}_n^G: BC_{2n-1}(n) \lra \R(n-1)
$$
by setting it zero on the groups $C_{2n-1}(n+i)$ for $i>0$. 
The second $(2n-1)$--term functional equation for 
the Grassmannian $n$--logarithm function, see (\ref{5}),  
just means that the composition 
$$
C_{2n}(n+1) \lra C_{2n-1}(n) \stackrel{{\cal L}_n^G}{\lra} \R(n-1),
$$
where the first map is a vertical arrow in $BC_*(n)$, is zero. 
Therefore we get a homomorphism
$$
{\cal L}_n^G: H_{2n-1}(BC_{*}(n)) \lra \R(n-1)
$$ 
 \begin{corollary} \label{1.20.11as}
One has 
$$
[C^G_n] = -(-1)^{n(n+1)/2}  \frac{(n-1)!^2}{(2n-2)!(2n-1)!} 
\cdot \frac{b_n}{n}
$$
\end{corollary}

{\bf Proof}. Indeed, 
$$
{\rm Alt}_{2n-1}{\rm Tr}(X_1 \cdot ... \cdot X_{2n-1}) = 
<{\rm tr}(g^{-1}dg)^{2n-1}_{|sl_n}, X_1 \wedge ... \wedge X_{2n-1}>
$$
So the claim  follows from Theorem \ref{1.20.10} since, as it clear from  comparison of formulas 
(\ref{7.1.02.1}) and (\ref{7.1.02.2}), 
$
b_n = n![C_n]
$. 
The corollary is proved.

\begin{theorem}  \label{6.24.02.7}
The composition 
\begin{equation} \label{6.24.02.6}
 K_{2n-1}(\C) \stackrel{\sim}{\lra} 
{\rm Prim} H_{2n-1}(GL_{2n-1}(\C), \Q) \stackrel{\varphi_{2n-1}^{2n-1}}{\lra} 
\end{equation} 
$$
H_{2n-1}(BC_{*}(n)_{\Q}) \stackrel{{\cal L}_n^G}{\lra} \R(n-1)
$$
equals to 
\begin{equation} \label{7.02.1.4}
-(-1)^{n(n+1)/2} \cdot \frac{(n-1)!^2}{n(2n-2)!(2n-1)!}r_{n}^{\rm Bo}
\end{equation}
\end{theorem} 

{\bf Proof}. Recall that 
restriction to  $GL_n$ of the map $\varphi_{2n-1}^{2n-1}$ 
coincides with the map $\varphi_{2n-1}^n$. Therefore  Lemma \ref{6.24.02.1} 
guarantees that restriction to $GL_n(\C)$ 
of the composition of the last two arrows 
 coincides with the map given by the class $[C_{2n-1}^G]$. 
So Corollary \ref{1.20.11as} implies the theorem.

{\bf 7. Comparing $[D_n]$ and $B_n$}. 
The following result is not used below. 
 
\begin{theorem} \label{1.20.11}
$$
[D_n] \quad = (2\pi i)^{n}(2n-1)B_n
$$
\end{theorem} 

{\bf Proof}. The transgression identifies the class $B_n$ with 
the Chern class of the universal bundle over $BG$, where $G = GL_n(\C)$. 
We will compute explicitly the transgression of the $n$-th component of the 
Chern character of the 
universal vector bundle $p: E \lra BG$. 
Let ${\cal A}$ be a connection 
on $E$. Then the $n$-th Chern class is 
represented by the $2n$-form 
$$
c_n(A):= \frac{tr F_{\cal A}^n}{{(2\pi i)}^n}  
$$
where $F_A := d{\cal A} + {\cal A} \wedge {\cal A}$ is the curvature form. 

Let $q: EG \lra BG$ be the principal fibration associated with $E$. 
Then the form $q^*c_n({\cal A})$ is exact. If 
$d\omega = q^*c_n({\cal A})$ and 
$F$ is a fiber of $q$ then $\omega|_F$ 
is closed, its cohomology class is transgressive, and goes to 
$[c_n({\cal A})]$. To do 
the computation we choose a connection ${\cal A}_0$ on $BG$ which is 
flat in a neighbourhood $U$
of a point $x \in BG$. It provides a trivialization of the bundle $E$ 
over $U$ as well as a trivialization  $\varphi: EG|_U \lra G \times U$. 

The bundle 
$q^*E$ has  a canonical trivialization. It provides a connection ${\cal A}_1$ 
on $q^*E$. So there are two connections, $q^*{\cal A}_0$ and 
${\cal A}_1$ on $q^*E$. One has ${\cal A}_1 = q^*{\cal A}_0 + g^{-1}dg$, where 
$(g,u) = \varphi(x)$. Let 
$$
{\cal A}(t):= t {\cal A}_1  +(1-t) q^*{\cal A}_0 = 
tg^{-1}dg + q^*{\cal A}_0
$$ 
It can be thought of as a connection on the lifting of the 
bundle $q^*E$ 
to $EG \times [0,1]$; here $t \in [0,1]$. The curvature $F(t)$ of this 
connection is 
$$
F(t)=  g^{-1}dg dt + t^2g^{-1}dg\wedge  g^{-1}dg
$$
The push forward of the form $trF(t)^n$ down to $EG$ is a primitive for the form 
$trF_{{\cal A}(1)}^n - trF_{{\cal A}(0)}^n$. It is given (in  $q^{-1}U$)  by
$$
\int_0^1tr F(t)^n = \frac{1}{2n-1}tr (g^{-1}dg )^{2n-1}
$$
Theorem \ref{1.20.11} is proved.

{\bf 8. On the 
motivic nature of the Grassmannian $n$--logarithm functions}. 
According to the results of the previous section understanding of the Borel regulator, 
and hence special values of the Dedekind $\zeta$--functions, is reduced to 
study of properties of the Grassmannian $n$--logarithm function 
${\cal L}_n^{G}$. 

Recall that a framed mixed Hodge-Tate structure has a natural $\R$-valued invariant ([BD]), called the Lie period. Thus a variation of Hodge-Tate structures ${\Bbb L}$ 
over a base $X$ 
provides  a period function ${\Bbb L}^{\cal M}$ on $X(\C)$.

\begin{conjecture} \label{6.24.02.10} There exists a variation ${\Bbb L}^{\cal M}_n$ 
of framed 
mixed Tate motives over $\widehat G_{n-1}^n$ such that:
$$
\sum_{i=0}^{2n}(-1)^i a_i^*{\Bbb L}_n^{\cal M} = 0; \quad 
\sum_{j=0}^{2n}(-1)^j b_j^*{\Bbb L}_n^{\cal M} = 0;
$$
and the  Lie period ${\cal L}^{\cal M}_n$ 
of its Hodge realization satisfies 
\begin{equation} \label{6.24.02.11}
{\cal L}^{\cal M}_n - {\cal L}^{G}_n = \sum_{i=0}^{2n-1}(-1)^i 
a_i^*{F}_n
\end{equation}
where $F_n$ is a function on $\widehat G_{n-2}^n(\C)$. 

b) The functional equations satisfied by  ${\cal L}^{\cal M}_n$ 
 essentially determine it: 
the space of all smooth/measurable functions satisfying these functional equations is finite dimensional. 
\end{conjecture}

{\bf Remark}. The function $F_n$ is obviously 
not determined by 
(\ref{6.24.02.11}) -- add a function coming from $\widehat G^n_{n-3}(\C)$. 
Nevertheless we expect that there exists  a {\it canonical} explicit 
choice for  $F_n$. Then formula (\ref{6.24.02.11}) can be considered 
as an explicit formula for 
${\cal L}^{\cal M}_n$ in terms of ${\cal L}^{G}$'s. 

Moreover we expect that there exists a canonical homotopy 
between the Grassmannian $n$-logarithm (understood as 
a cocycle in the Deligne cohomology of the bi-Grassmannian) 
and its ``motivic'' 
bi-Grassmannian counterpart. 
Observe that the motivic bi-Grassmannian $n$-logarithm should
have  non-trivial components outside of the bottom line of 
the bi-Grassmannian, while the defined above (or in [G5]) 
Grassmannian $n$--logarithm 
 is concentrated entirely 
at the bottom line.

A variation of mixed Tate motives over $\widehat G_{n-1}^n$ 
was constructed in [HaM]. However it is not clear how to relate it to 
the function ${\cal L}^G_n$.

Conjecture \ref{6.24.02.10} is known for $n=2$ and $n=3$. 

The $n=2$ case 
follows from (\ref{6.13.02.101}), the well known 
motivic realization of the dilogarithm, and 
Bloch's theorem characterizing the Bloch-Wigner function by 
Abel's $5$--term equation it satisfies. 

The $n=3$ case of conjecture \ref{6.24.02.10} follows from the results of 
[G1-2], [GZ] and the motivic realization of the trilogarithm. 
In particular the part b) is given by  
 Theorem 1.10 in [G1]. 

{\bf Examples}. 1. $n=2$. Then  ${\cal L}^{\cal M}_2 = {\cal L}^{G}_2$. 

2. $n= 3$. The motivic Grassmannian trilogarithm function 
has been constructed in [G1-2] in terms of the 
classical trilogarithm function. Namely, one has 
$$
{\cal L}^{\cal M}_3(l_0, ..., l_5) = \frac{1}{90} {\rm Alt}_6
{\cal L}_3\Bigl( \frac{\Delta(l_0, l_1, l_3)\Delta(l_1, l_2, l_4)
\Delta( l_2, l_0, l_5)}{\Delta(l_0, l_1, l_4)
\Delta(l_1, l_2, l_5)\Delta(l_2, l_0, l_3)}
\Bigr) 
$$

According to Theorem 1.3 of [GZ]  ${\cal L}^{\cal M}_3$ 
is different from ${\cal L}^{G}_3$, and  
$$
F_3(l_0, ..., l_4) = \frac{1}{9}{\rm Alt}_5\Bigl(\log 
|\Delta(l_0, l_1, l_2)| \log |\Delta(l_1, l_2, l_3)| \log |\Delta(l_2, l_3, l_4)| \Bigr) 
$$

3. If $n>3$  then ${\cal L}^{\cal M}_n$ 
is different from ${\cal L}^{G}_n$ since it is already so 
for  the restriction to the special configuration, 
see Theorem \ref{clpoly}.

The space of the functional equations 
for the function ${\cal L}^G_3$  is smaller than the one for 
${\cal L}^{\cal M}_3$, see Chapter 1.5 of [G3].   
A similar situation is expected  
for all $n>3$.

The space of the functional equations for the motivic 
$n$--logarithm function ${\cal L}^{\cal M}_n$ should provide an 
explicit construction of the weight $n$ part of the motivic Lie coalgebra of 
an arbitrary field $F$, as explained in s. 4.1 in [G6], taking into 
account the following correction. 

{\it Correction}. In s. 4.2 of [G6] the subgroup of the functional equations ${\cal R}_n^G$ is supposed to be defined as the subgroup of all functional equations for the function ${\cal L}^{\cal M}_n$, not ${\cal L}^{G}_n$.

\section{The Chow dilogarithm and a reciprocity law}

The Chow dilogarithm provides a homomorphism $\Lambda^3 \C(X)^* \to \R$ 
given by 
\begin{equation} \label{chowhh}
  f_1\wedge f_2\wedge f_3 \lms {\cal P}_2(X;f_1,f_2,f_3):= \frac{1}{2\pi i}
\int_{X(\C)}r_2(f_1,f_2,f_3)
\end{equation}

In this section we show that the Chow dilogarithm can be 
expressed by the function ${\cal L}_2(z)$.  
The precise versions of this claim are discussed below.
 

{\bf 1. The set up ([G1-2])}. 
For any field $F$ we defined in [G1] the  groups 
$$
{\cal B}_n(F):= \frac{\Z[F^*]}{{\cal R}_n(F)}, \quad n \geq 2
$$
and homomorphisms
$$
\delta_n : {\cal B}_n(F) \lra {\cal B}_{n-1}(F)\otimes F^*; \quad 
\{x\}_n \lms \{x\}_{n-1} \otimes x, \quad n\geq 3
$$
$$
\delta_2: 
{\cal B}_2(F) \lra \Lambda^2F^*, \quad \{x\}_2 \lms (1-x) \wedge x
$$
There is a complex $\Gamma(F;n)$
$$
{\cal B}_{n}(F) \stackrel{\delta_n}{\longrightarrow} {\cal B}_{n- 
1}(F)\otimes F^{\ast}  \stackrel{\delta_n}{\rightarrow}  \ldots \stackrel{\delta_n}{\rightarrow} {\cal B}_{2}(F)\otimes  \Lambda^{n-2}F^{\ast}
\stackrel{\delta_n}{\rightarrow}
 \Lambda^{n}F^{\ast}
$$
where $\delta_n(\{x\}_k \otimes Y) := \{x\}_{k-1} \otimes x \wedge Y$ for $k>2$, and $(1-x)\wedge x \wedge y$ for $k=2$, 
called the weight $n$ polylogarithmic complex. 

If $K$ is a  field with 
 a discrete valuation $v$ and  the residue field $
k_v$, then there is a homomorphism of complexes 
${\rm res}_v: \Gamma(K, n) \lra \Gamma(k_v, n-1)[-1]$ (see s.\ 14 of \S
1 in [G1]). For example for $n=3$ we have 
\begin{equation} 
\begin{array}{ccccc}
{\cal B}_3(K)&\stackrel{\delta_3}{\longrightarrow}&{\cal B}_2(K)\otimes
K^{\ast}&
\stackrel{\delta_3}{\longrightarrow}&\Lambda^3K^{\ast}\\
  &&\downarrow {\rm res}_v &&\downarrow  {\rm res}_v \\
&&{\cal B}_2(k_v)&\stackrel{\delta_2}{\longrightarrow}&\Lambda^2k_v^{\ast}
\end{array}
\end{equation}
Here ${\rm res}_v(\{x\}_2 \otimes y)$ is zero unless $v(x) =0$. In the latter case 
it is ${\rm res}_v(\{x\}_2 \otimes y) = v(y) \{\overline x\}_2$, where $\overline x$ 
denotes  projection of $x$ to the residue field of $K$.

Let $X$ be a regular  curve over an algebraically closed field
$k$ and $F:= k(X)^{\ast}$.   
 Set ${\rm Res}:= \sum_x {\rm res}_x$ where ${\rm res}_x$ is
the residue homomorphism for the valuation on $F$
 corresponding to a
point $x$ of $X$. For instance for $n=3$ we get a morphism of complexes
$$
\begin{array}{ccccc}
{\cal B}_3(F)&\stackrel{\delta_3}{\longrightarrow}&{\cal B}_2(F )\otimes
F^{\ast}&
\stackrel{\delta_3}{\longrightarrow}&\Lambda^3F^{\ast}\\
  &&\downarrow {\rm Res}&&\downarrow {\rm Res}\\
&& {\cal B}_2(k )&\stackrel{\delta_2}{\longrightarrow}& \Lambda^2k^{\ast}
\end{array}
$$
  
We will also need a more explicit version $B_2(F)$ 
of the group ${\cal B}_2(F)$. 
Denote by $R_2(F)$  the subgroup of $\Z[{\Bbb P}^1(F)]$ 
generated by the elements $$
\{0\}, \{\infty\}\quad \mbox{and} \quad \sum_{i=1}^{5}(-1)^i
\{r(x_1,...,\widehat   x_i,...,x_5)\}
$$ 
when $(x_1,...,x_5)$ runs through all $5$-tuples of distinct points in ${\Bbb P}^1(F)$,
and the cross-ratio $r(...)$ is normalized by $r(\infty, 0, 1, x) =x$.  Define the Bloch group $B_2(F)$ as
$$
B_2(F):= \frac{\Z[{\Bbb P}^1(F)]}{R_2(F)}
$$

One can show that $R_2(F)\subset {\cal R}_2(F)$. 
So there is a map 
\begin{equation} \label{6.18.02.10}
i: B_2(k)\lra {\cal B}_2(k)
\end{equation}
induced by the identity map on the generators. 

\begin{proposition} \label{mot}
Let $k$ be a number field. Then (\ref{6.18.02.10})  is an isomorphism 
modulo torsion.
\end{proposition}

{\bf Proof}. The map $i$ is clearly surjective. The diagram
$$
\begin{array}{ccc}
B_2(k) &\lra&{\cal B}_2(k) \\
&&\\
\downarrow \delta_2&&\downarrow \delta_2\\
&&\\
\Lambda^2k^*&=&\Lambda^2k^*
\end{array}
$$
is commutative. So we need only to show that if $0 \not = x \in B_2(k)_{\Q}$ and $\delta_2(x)=0$, then $i(x)\not =0$. This  follows from the injectivity of the regulator map on $K^{ind}_3(k)_{\Q}$.  
Indeed, by Suslin's theorem  for a field $F$ one has $K^{ind}_3(F)_{\Q} = {\rm Ker} \delta_2\otimes \Q $. Let us identify $K^{ind}_3(\C)_{\Q}$ with this subgroup of $B_2(\C)_{\Q} $. The restriction of the dilogarithm map
$$
B_2(k) \lra {\cal B}_2(k)\lra  (\Z[Hom(k,\C)]\otimes 2 \pi i \R)^+, \qquad \{z\}_2 \lms \{2 \pi i{\cal L}_2(\sigma_i(z))\}
$$
 to the subgroup ${\rm Ker} \delta_2\otimes \Q $  gives the  Borel regulator   $K^{ind}_3(k)_{\Q} \lra \R^{r_2}$ ([G1]) and  thus is injective   by Borel's theorem.  

{\bf Remark}. For any field  $k$ the rigidity conjecture for $K_3^{ind}$ 
implies that the map $i$ should be an isomorphism, see [G1].

{\bf 2. The strong reciprocity law}. 
 \begin{conjecture} \label{reh}
Let $X$ be a regular projective curve over an algebraically closed field
$k$ and $F:= k(X)^{\ast}$.
Then there exists a canonical homomorphism of groups  
$
h:\Lambda^3F^* \rightarrow
{\cal B}_2(k)
$ 
satisfying the following two conditions:

a) $h(k^* \wedge \Lambda^2F^*) =0$ and the diagram 
\begin{equation} \label{6.17.02.3}
\begin{array}{ccccc} \label{reh1}
{\cal B}_3(F)& \stackrel{\delta_3}{\longrightarrow}&{\cal B}_2(F)\otimes F^{\ast} 
& \stackrel{\delta_3}{\longrightarrow} & \Lambda^3F^{\ast}\\
&&&&\\
&&{\rm Res}\downarrow&h \swarrow&\downarrow {\rm Res}\\
&&&&\\
&&{\cal B}_2(k)&\stackrel{\delta_2}{\longrightarrow}&\Lambda^2 k^{\ast} 
\end{array}
\end{equation}
is commutative.

b) If $X$ is a  curve over $\C$
then 
\begin{equation} \label{homotq}
\frac{1}{2\pi i}\int_{X(\C)}r_2(f_1\wedge f_2 \wedge f_3) =   {\cal L}_2 \Bigl(h(f_1\wedge
f_2 \wedge f_3)\Bigr) 
\end{equation}
\end{conjecture}

{\bf Remarks}. 1. b) follows easily if we have a functorial map $h$ such that 
${\rm Res} = \delta_2 \circ h $, see Theorem \ref{55} below.  

2. According to  Suslin's 
 reciprocity law for the Milnor group $K^M_3(F)$ the projection of 
${\rm Res}(\Lambda^3F^{\ast}) \subset \Lambda^2k^*$ to $K_2(k)$ is zero. 
Since by Matsumoto's 
theorem  $K_2(k)={\rm Coker}(\delta_2)$,  one has 
${\rm Res}(\Lambda^3F^{\ast}) \subset {\rm Im}(\delta_2)$. However  
${\rm Ker}(\delta_2)$ is nontrivial, so  it is {\it a priory} 
unclear that we can lift naturally  the map ${\rm Res}$  to a map $h$.
 One  of the reasons   why we can do   this  is provided by (\ref{homotq}).

We prove  this conjecture  in the following cases: 

a)  $X={\Bbb P}^1$; we construct {\it explicitly} a reciprocity homomorphism  
$h: \Lambda^3F^*\to  B_2(k)$ in Theorem \ref{p4}.

b)  $X$ is an elliptic curve over an  algebraically closed  field; we construct {\it explicitly} a reciprocity homomorphism  $h: \Lambda^3F^* \to  {\cal B}_2(k)$ in Theorem \ref{ecur}.

c) $k = {\overline \Q}$, $X$ is any curve; see Theorem \ref{rl}.

In the cases a)  and b)  the homomorphism $h$  satisfies the following additional property. Let $F = k(X)$ and $k$ is not necessarily algebraically closed. Let $k'$ be the field of definition of the divisors $(f_1), (f_2), (f_3)$. Then $h(f_1\wedge f_2\wedge f_3) \in {\cal B}_2(k')$. 

\begin{conjecture}
Let $X$ be a projective regular curve over an algebraically closed field
$k$ and $F:= k(X)$. Then the  homomorphism 
$$
{\rm Res}: \Gamma(F;n) \longrightarrow \Gamma(k;n-1)[-1]
$$  
is homotopic to zero.
\end{conjecture}

\begin{lemma} \label{mot1}
Assume that we have a map $h$ such that $h(k^* \wedge \Lambda^2F^*) =0$ and 
${\rm Res} = \delta_2 \circ h $. Then $ h \circ \delta_3 = {\rm Res}$.
\end{lemma}

{\bf Proof}. The image of the group ${\cal B}_2(F)\otimes F^{\ast}$ under the map $h \circ \delta_3$ belongs to the subgroup ${\rm Ker}\delta_2$. Since $h(k^* \wedge \Lambda^2F^*) =0$ one has $h \circ \delta_3 = {\rm Res}$ on ${\cal B}_2(F)\otimes k^{\ast}$. Any element of $k(X)^*$ can be connected via a curve to a constant. This together with the rigidity of ${\rm Ker}\delta_2$ (built into the definition of the group ${\cal B}_2(k)$) implies the result.

{\bf 3. The $X= \PP^1$ case}. 
Recall that  $v_x(f)$ is the order of zero of $f\in k(X)$ at $x$. 
Choose a point $\infty$ on ${\Bbb P}^1$.

\begin{theorem} \label{p4}
 Assume that $k = \overline k$. Then the map $
h:\Lambda^3 k({\Bbb P}^1)^* \rightarrow
 B_2(k)
$   given by the formula 
$$
h(f_1\wedge f_2\wedge f_3):=  \sum_{x_i \in {\Bbb P}^1(k)} 
v_{x_1}(f_1)v_{x_2}(f_2)v_{x_3}(f_3)\{r(x_1, x_2, x_3, \infty)\}_2
$$
satisfies all the conditions of conjecture \ref{reh} modulo $6$-torsion.
\end{theorem}

{\bf Proof}. Let us show that  $h$ is independent of the choice of  $\infty$, i.e. 
$$
\sum_{x_i \in {\Bbb P}^1(k)}
v_{x_1}(f_1)v_{x_2}(f_2)v_{x_3}(f_3)\{r(x_1, x_2, x_3, a)\}_2
\in B_2(k)
$$
does not depend on  $a$. Indeed, the 5-term
relation  for the 5-tuple of points $(x_1, x_2, x_3, a, b)$ gives
$$
\sum_{x_i \in {\Bbb P}^1(k)}
v_{x_1}(f_1)v_{x_2}(f_2)v_{x_3}(f_3)\Bigl(\{r(x_1, x_2, x_3, a)\}_2 -
\{r(x_1, x_2, x_3, b)\}_2\Bigr) =
$$
$$
-\sum_{x_i \in {\Bbb P}^1(k)}
v_{x_1}(f_1)v_{x_2}(f_2)v_{x_3}(f_3)\Bigr(\{r(x_1, x_2, a, b )\}_2
-\{r( x_1, x_3, a, b)\}_2 + \{r(x_2, x_3, a, b)\}_2\Bigr)
$$
Each of these 3 terms vanishes because $\sum_{x \in {\Bbb P}^1(k)}v_{x}(f) =0$.
for any $f \in k(\PP^1)^*$. 

 \begin{proposition} \label{not}
Let $k= \overline k$. Then modulo  $6$-torsion
$$
h((1-f) \wedge f \wedge g) = \sum_{x \in {\Bbb P}^1(k)}v_{x}(g)\{f(x)\}_2
$$
\end{proposition}

{\bf Proof}. Using linearity with respect to $g$ and projective invariance of the cross--ratio we see that it is sufficient to prove the identity for $g=t$. 
Then it boils down to 
\begin{equation} \label{7.3.02.2}
\sum_{x_i \in {\Bbb P}^1(k)}
v_{x_1}(1-f)v_{x_2}(f)\{r(x_1, x_2, 0, \infty)\}_2  = 
\{f(0)\}_2 - \{f(\infty)\}_2
\end{equation}

\begin{lemma} \label{notxzc}
Applying $\delta_2$ to both parts of 
(\ref{7.3.02.2}) we get the same result  modulo  $6$-torsion.
\end{lemma} 

{\bf Proof}. 
Choose a coordinate $t$ on $\PP^1$ such that $f(\infty) =1$. 
Then 
\begin{equation} \label{7.3.02.4}
f(t) = \frac{\prod_i (a_i -t)^{\alpha_i}}{\prod_k (c_k -t)^{\gamma_k}}; \quad 
1-f(t)= \frac{B \prod_j (b_j -t)^{\beta_j}}{\prod_k (c_k -t)^{\gamma_k}}
\end{equation}
Observe that  $\{f(\infty)\}_2 =0$ modulo $6$-torsion. 
The left hand side equals 
\begin{equation} \label{7.3.02.1}
\sum_{x_i \in {\Bbb P}^1(k)}
v_{x_1}(1-f)v_{x_2}(f)\{x_1/x_2\}_2  =
\end{equation}
$$
\sum \beta_j \alpha_i \{b_j/a_i\}_2 - 
\sum \gamma_k \alpha_i \{c_k/a_i \}_2
- \sum \beta_j\gamma_k \{b_j/c_k\}_2  
$$
Applying $\delta_2$ to it we get 
$$
\sum \beta_j \alpha_i \cdot \frac{a_i-b_j}{a_i} \wedge \frac{b_j}{a_i} - 
\sum \gamma_k \alpha_i \cdot \frac{a_i- c_k}{a_i} \wedge \frac{c_k}{a_i} 
- \sum \gamma_k \beta_j \cdot \frac{c_k - b_j}{c_k} \wedge \frac{b_j}{c_k} = 
$$
$$
\sum \beta_j \alpha_i \cdot b_j \wedge a_i + 
\sum \alpha_i\gamma_k  \cdot a_i \wedge c_k  
+ \sum \gamma_k \beta_j \cdot c_k \wedge b_j \quad + 
$$
$$
- \sum_i\frac{\prod_j (a_i-b_j)^{\beta_j }}{\prod_k (a_i-c_k)^{\gamma_k }} 
\wedge   a_i^{\alpha_i}  
+   \sum_j\frac{\prod_i (a_i - b_j)^{\alpha_i}}
{\prod_k (c_k - b_j)^{\gamma_k}} \wedge  b_j^{\beta_j}
-  \sum_k \frac{\prod (a_i- c_k)^{\alpha_i}}{\prod (c_k  - b_j)^{\beta_j}} 
\wedge  c_k^{\gamma_k} 
$$
Using (\ref{7.3.02.4}) we see that the second line equals modulo $2$--torsion 
to
$$ 
- \sum_i (1-f(a_i)) \wedge a_i^{\alpha_i} + \sum_j f(b_j) 
\wedge b_j^{\beta_i}  - \sum_k
 \frac{\prod (a_i- c_k)^{\alpha_i}}{\prod (b_j - c_k)^{\beta_j}} 
\wedge   c_k^{\gamma_k} + 
B \wedge \prod a_i^{\alpha_i}
$$
The first two terms are zero since $f(a_i)=0$ 
and $f(b_j) =1$. The third term equals to 
$-(B \wedge \prod c_k^{\gamma_k})$ since, as it follows from (\ref{7.3.02.4}), 
$f(c_k) = \infty$ and thus 
$$
\frac{\prod (a_i- c_k)^{\alpha_i}}{\prod (b_j - c_k)^{\beta_j}}  = -B 
$$
On the other hand, we have 
$$
\delta_2(\{f(0)\}_2) = (1-f(0)) \wedge f(0) =  
\frac{B\prod_j b_j^{\beta_j}}{\prod_k c_k^{\gamma_k}} \wedge 
\frac{\prod_i a_i^{\alpha_i}}{\prod_k c_k^{\gamma_k}} 
$$
which matches the  expression we got for the left hand side. The lemma is proved.

To prove the proposition it remains to use a rigidity argument.
Namely, we need to show that the identity is valid 
for some particular $f$, which is easy, or use proposition \ref{p44} 
plus injectivity of the regulator on $K_3(\overline \Q)_{\Q}$. 
The proposition is proved. 

Now let us prove the key fact that  ${\rm Res} = \delta_2\circ h$. We need to show that 
 for any 3 rational functions $f_1, f_2, f_3$ on
${\Bbb P}^1$  
$$
\sum_{x \in {\Bbb P}^1(k)} {\rm res}_x (f_1\wedge f_2\wedge f_3) = 
$$
\begin{equation} \label {1}
\delta_2  \Bigl( \sum_{x_i \in {\Bbb P}^1(k)}
v_{x_1}(f_1)v_{x_2}(f_2)v_{x_3}(f_3)\{r(x_1, x_2, x_3, \infty)\}_2 \Bigr)
\end{equation}

Both sides are obviously homomorphisms
from $\Lambda^3 F^{\ast}$ to $\Lambda^2 k^{\ast}$ which are zero on
$k^{\ast}  \wedge \Lambda^2 F^{\ast}$. (The last property for the map
$\sum {\rm res}_x$ 
is provided by the Weil reciprocity law). We normalize the cross ratio of four points on the projective 
line by  $r(\infty,0,1,z) =z$. 
 So it suffices to
check the formula on elements
$$
\frac{z-a_2}{z-a_1} \wedge \frac{z-b_2}{z-b_1} \wedge \frac{z-c_2}{z-c_1}
$$
In this case it follows from 
$$
\delta_2 \{r(a_2, b_2, c_2, \infty)\}_2 = \delta_2\left\{\frac{a_2-c_2}{b_2-c_2}\right\}_2 = \frac{b_2-a_2}{b_2-c_2}\wedge \frac{a_2-c_2}{b_2-c_2} 
$$

It remains to prove the following proposition.

\begin{proposition} \label{p44}
$$
{\cal P}_2({\Bbb P}^1;f_1,f_2,f_3)=   \sum_{x_i \in {\Bbb P}^1(\C)} v_{x_1}(f_1)v_{x_2}(f_2)v_{x_3}(f_3){\cal L}_2(r(x_1, x_2, x_3, \infty))
$$
\end{proposition}

{\bf Proof}.  We  immediately reduce the statement  to the situation when $f_1 = 1-z, f_2 = z, f_3 = z-a$ 
which is a particular case of the following lemma  

\begin{lemma} \label{l454} Let $X$ be an arbitrary curve over $\C$. Then
$$
\int_{X(\C)}r_2((1-f)\wedge f\wedge g)= - \sum_{x \in X(\C)} v_{x}(g) \cdot
{\cal L}_2(f(x))
$$
\end{lemma}

{\bf Proof}. For  functions $f(z)$ and $g(z)$ on $X(\C)$ set
$$
\alpha(f,g): = \log|f|d\log|g| - \log|g|d\log|f|
$$
 Consider the following 1-form on $X(\C)$
\begin{equation} \label{99}
  {\cal L}_{2}(f) d\arg g  -\frac{1}{3} \alpha(1-f,f)\log|g|
\end{equation}
It defines a current on $X(\C)$. We claim that its derivative is equal to:
\begin{equation} \label{11}
  2 \pi \cdot {\cal L}_{2}(f)\delta(g)  + r_2((1-
f) \wedge f \wedge g) 
\end{equation}
Using $d (d \arg g) = 2\pi  \cdot \delta(g)$
and 
\begin{equation} \label{dl2}
d{\cal L}_2(z) = -\log|1-z| d\arg z  + \log|z| d\arg(1-z)
\end{equation}
we see that the  differential of the current (\ref{99}) equals to
$$
2 \pi  {\cal L}_{2}(f)\delta(g) \quad + \quad \Bigl( -\log|1-f| d\arg f  + \log|f| d\arg(1-f)\Bigr)\wedge d\arg g +
$$
$$
\frac{1}{3}\Bigl(\log|1-f|d\log|f| - \log|f|d\log|1-f|\Bigr)\wedge d\log|g| - \frac{2}{3}\log|g| \cdot d\log|1-f|\wedge d\log|f|
$$
Since $d\log(1-f)\wedge d\log f  = 0$ we have $$
d\log|1-f|\wedge d\log|f| = d\arg(1-f)\wedge d\arg f$$
 Using this and writing  $r_2(f_1 \wedge f_2 \wedge f_3)$ as
$$
\frac{1}{3}(\log|f_1|d\log|f_2|\wedge d\log|f_3| + \quad \mbox{cyclic permutations})
$$
$$
 -  (\log|f_1|d\arg(f_2)\wedge d\arg f_3 + \quad \mbox{cyclic permutations})
$$
we come to (\ref{11}).
Integrating we get the lemma. The theorem is proved. 

{\bf 4. Expressing the Chow dilogarithm via  the classical one}. Let $\pi: 
Y \rightarrow S$ be a family of curves over a base $S$ over $\C$ and 
$f_1,f_2,f_3 \in \C(Y)^*$. 
We get a function at the generic point  of $S$.  Its value at 
 $s \in S$ is given by the Chow dilogarithm 
${\cal P}_2(Y^s;f^s_1,f^s_2,f^s_3)$, where 
$Y^s$ is the fiber of 
$\pi$ at $s$. Denote it by 
${\cal P}_2(Y \rightarrow S ;f_1,f_2,f_3)$.

\begin{theorem} \label{55}
a) Let $\pi: 
Y \rightarrow S$ be a family of curves over a  base $S$ over $\C$. Then 
there are rational functions $\varphi_i$ on $S$ such that  
$$
{\cal P}_2(Y \rightarrow S;f_1,f_2,f_3) = \sum_i {\cal L}_2(\varphi_i(s))
$$

b) Let $k = \C(S)$, $X$ is the generic fiber of $\pi$, and $F=k(X)$. Suppose that there exists a map $h:\Lambda^3F^* \rightarrow
{\cal B}_2(k)$ such that ${\rm Res} = \delta_2 \circ h $. Then 
\begin{equation} \label{3211}
 d{\cal P}_2(Y \rightarrow S;f_1,f_2,f_3) \quad = \quad d{\cal L}_2(h(f_1,f_2,f_3))
\end{equation}
\end{theorem}

{\bf Proof}. 
a) 
We use  the existence of the transfer map on  $K^M_3$  
to reduce the statement  to the case  $X= {\Bbb P}^1$. 

Choose a projection $p: X \rightarrow {\Bbb P}^1$. 
We may suppose without loss of generality that $p$ is a (ramified) 
Galois covering 
with the Galois group $G$. 
 Indeed, let  $p_1:Y \rightarrow X$ 
be such a covering that its composition with $p$ 
 is a Galois covering. Indeed, 
$$
{\cal P}_2(Y\lra S; p_1^*f_1,p_1^*f_2, p_1^*f_3) = \frac{1}{{\rm deg} p_1}{\cal P}_2(X\lra S; f_1, f_2, f_3)
$$

Then $\sum_{g \in G} g^*\{f_1, f_2, f_3\} \in p^*K_3^M(k({\Bbb P}^1))$. 
It coincides with $p^*$ of the transfer of the element $\{f_1, f_2, f_3\} \in K_3^M(F)$. 
This means that there exist  $s^{(i)}_1,s^{(i)}_2,s^{(i)}_3  \in k({\Bbb P}^1)$ 
and $g_j,h_j  \in k(X)$ such that
\begin{equation} \label{6.18.02.1}
\sum_{g \in G} g^*(f_1\wedge f_2 \wedge f_3) - p^*\sum_i 
s^{(i)}_1\wedge  s^{(i)}_2\wedge s^{(i)}_3 =  
\sum_j(1-g_j)\wedge g_j \wedge h_j
\end{equation} 
Therefore
$$
 {\cal P}_2(Y \rightarrow S;f_1,f_2,f_3)= 
\frac{1}{|G|}\sum_i {\cal P}_2({\Bbb P}^1 \times S \rightarrow S;
 s^{(i)}_1, s^{(i)}_2, s^{(i)}_3) + 
$$
$$
\sum_j {\cal P}_2({\Bbb P}^1 \times S \rightarrow S; (1-g_j), g_j, h_j)
$$
It remains to use Lemma  \ref{l454}  and Proposition  \ref{p44}. 
The part a) of the theorem is proved. 

b) We need the following lemma.
\begin{lemma} \label{2.14}
\begin{equation} \label{dhomot}
d {\cal P}_2(Y \rightarrow S;f_1,f_2,f_3) =  (2 \pi)^{-1} \cdot {\rm Alt}_3 \Bigl(v_x(f_1) \log|f_2(x)| d_s \arg f_3(x) \Bigr)
\end{equation}
\end{lemma}

{\bf Proof}. Using $dd \log f  = 2\pi i \delta(f)$  we get  an equality of $3$-currents on $Y$
$$
dr_2(f_1,f_2,f_3) = 
$$
\begin{equation} \label{2233}
 \pi_3\Bigl(\frac{df_1}{f_1}\wedge\frac{df_2}{f_2}\wedge\frac{df_3}{f_3}\Bigr) +  2 \pi  \cdot {\rm Alt}_3 \Bigl(\delta(f_1) \log|f_2(x)| d \arg f_3(x) \Bigr)
\end{equation}
The second term in (\ref{2233}) is a $1$-form on the divisor 
$D:= \cup_{i=1}^3{\rm div}(f_i)$ considered as  a $3$-current on $Y$.  
This $1$-form is the composition of the residue map
$$
{\rm res}: \Lambda^3 \C(Y)^* \longrightarrow \coprod_{X \in Y_1}
\Lambda^2\C(X)^*
$$
with the map
$$
r_1: \Lambda^2 \C(X)^* \longrightarrow  {\cal A}^1({\rm Spec }(\C(X))), \qquad f \wedge g \longmapsto -2\pi (\log|f| d\arg g - \log|g| d\arg f)
$$

The push forward of the first term in (\ref{2233}) vanishes  
(since the fibers are complex curves the push down  of any $(3,0)$-form 
to $S$ is zero). Integrating the second $3$-current  in (\ref{2233}) 
along the fibers of $Y$ we get (\ref{dhomot}). The lemma is  proved.

 According to  Lemma \ref{2.14}  and formula (\ref{dl2}) for $d {\cal L}_2$, and using ${\rm Res} = \delta_2 \circ h $ we get the proof of the 
part b) of the theorem.

{\bf Remark}  The  function ${\cal L}_2(z)$ is continuous on $\C {\Bbb P}^1$. Therefore part a) of the theorem 
implies that the function ${\cal P}_2(Y \rightarrow S ;f_1,f_2,f_3)$ can be extended to a continuous 
function on $S$.

{\bf 5. Conjecture \ref{reh} for $\overline k = \overline \Q$}. 
\begin{theorem} \label{rl}
Let $X$ be a   regular  projective curve over  
${\overline \Q}$ and $F:= {\overline \Q}(X)$. 
Then there exists a homomorphism   
$
h:\Lambda^3F^* \rightarrow
B_2({\overline \Q})\otimes \Q
$ as in conjecture \ref{reh}  such that for any embedding $\sigma: {\overline \Q} 
\hookrightarrow  \C$ one has 
\begin{equation} \label{homot}
\frac{1}{2\pi i}\int_{X(\C)}r_2(\sigma(f_1\wedge f_2 \wedge f_3)) =  {\cal L}_2 \Bigl(\sigma(h(f_1\wedge
f_2 \wedge f_3))\Bigr) 
\end{equation}
\end{theorem}

{\bf Proof}. It is similar to the proof of theorem \ref{55}. 
Choose a projection $p:X \rightarrow {\Bbb P}^1$. 
We may suppose that $p$ is a Galois covering 
with the Galois group $G$. 
 Indeed, let  $p_1:Y \rightarrow X$ 
be such a covering that its composition with $p$ 
 is a Galois covering. 
Setting  
$$
h(f_1\wedge f_2 \wedge f_3):= h(\frac{1}{{\rm deg}p_1} \cdot p_1^*(f_1\wedge f_2 \wedge f_3))
$$ 
we may suppose that $p$ is Galois.  

Then $\sum_{g \in G} g^*\{f_1, f_2, f_3\} \in p^*K_3^M(k({\Bbb P}^1))$. 
So there exist  $s^{(i)}_1,s^{(i)}_2,s^{(i)}_3,  \in k({\Bbb P}^1)$ 
and $g_j,h_j  \in k(X)$ such that (\ref{6.18.02.1}) holds. 
Set
$$
|G| \cdot h(f_1\wedge f_2 \wedge f_3) := \sum_i h(s^{(i)}_1\wedge  s^{(i)}_2\wedge s^{(i)}_3) + \sum_j \sum_{x \in X(\overline \Q)}\{g_j(x)\}_2 \cdot v_x(h_j)
$$
   \begin{lemma} \label{1.8}
Suppose $\sum_i (1-f_i) \wedge f_i \wedge g_i = 0$ in  $\Lambda^3 {\overline \Q}(X)^*$. 
Then $$
\sum_i \sum_{x \in X(\overline \Q)} v_x(g_i) \cdot \{f_i(x)\}_2 =0 \quad \mbox{in the group $ B_2({\overline \Q})$.}
$$ 
\end{lemma}

 This lemma implies that $h$ is well defined.
Indeed, suppose that we have a different presentation  
$$
\sum_{g \in G} g^*\{f_1, f_2, f_3\} = p^*\sum_k 
\widetilde  s^{(k)}_1\wedge  \widetilde  s^{(k)}_2\wedge \widetilde  
s^{(k)}_3 +  \sum_j(1-\widetilde  g_j)\wedge \widetilde  g_j \wedge \widetilde  h_j
$$
We need to show that 
\begin{equation} \label{**}
h\Bigl( \sum_k \widetilde  s^{(k)}_1\wedge  \widetilde  s^{k)}_2\wedge 
\widetilde  s^{(k)}_3 - 
 \sum_i  s^{(i)}_1\wedge   s^{(i)}_2\wedge  s^{(i)}_3\Bigr) + 
\end{equation} 
$$
\sum v_x(h_j)\{g_j(x)\}_2 - \sum v_x(\widetilde  h_j)\{\widetilde  g_j(x)\}_2 = 0
$$
There exist $a_j, b_j \in k({\Bbb P}^1)$ such that modulo $k^* \wedge 
\Lambda^2k({\Bbb P}^1)^*$  one has   
$$
 \sum_k \widetilde  s^{(k)}_1\wedge  \widetilde  s^{(k)}_2\wedge \widetilde  s^{(k)}_3 - 
 \sum_i  s^{(i)}_1\wedge   s^{(i)}_2\wedge  s^{(i)}_3 - \sum_j (1-a_j)\wedge a_j 
\wedge b_j =0
$$
According to  Theorem \ref{p4} the homomorphism $h$ for ${\Bbb P}^1$ annihilates 
the left hand side. On the other hand 
$$
p^*\Bigl( \sum_k \widetilde  s^{(k)}_1\wedge  \widetilde  s^{(k)}_2\wedge \widetilde  
s^{(k)}_3 - 
 \sum_i  s^{(i)}_1\wedge   s^{(i)}_2\wedge  s^{(i)}_3\Bigr)  - 
$$
$$
\sum_j(1-\widetilde  g_j)\wedge \widetilde  g_j \wedge \widetilde   h_j +\sum_j(1- g_j)\wedge  g_j \wedge h_j = 0
$$
Using Lemma \ref{1.8} we get (\ref{**}). 
To get (\ref{homot}) we use Theorem \ref{55} and notice that
\begin{equation} \label{***}
{\cal P}_2(Y \rightarrow S;f_1,f_2,f_3) = 1/m \cdot {\cal P}_2(Z \rightarrow S;
p_1^*f_1, p_1^*f_2, p_1^*f_3)
\end{equation}

{\bf Proof of Lemma \ref{1.8}}.  
For a regular curve $X$ over an algebraically closed field $k$ there is a commutative diagram  
$$
\begin{array}{ccc} \label{reh1a}
B_2(F)\otimes F^{\ast} 
& \stackrel{\delta_3}{\longrightarrow} & \Lambda^3F^{\ast}\\
&&\\
{\rm Res}\downarrow&&\downarrow {\rm Res}\\
&&\\
B_2(k)\otimes \Z[X(k)]&\stackrel{\delta_2}{\longrightarrow}&\Lambda^2 k^{\ast}\otimes \Z[X(k)] 
\end{array}
$$
Thus for any point $x$ of the curve $X$ the element $\sum_{i} v_x(g_i) \cdot 
\{f_i(x)\}_2 $ lies in ${\rm Ker }\delta_2$. Therefore it defines an 
element $\gamma_x \in K_3({\overline \Q})_{\Q}$.

 For any  embedding $\sigma: {\overline \Q} 
\hookrightarrow  \C$ the value of the 
Borel regulator on $\sigma(\gamma_x)$  is equal to $\sum_{x \in X} v_x(g_i) \cdot {\cal L}_2(\sigma(f_i(x)))$. 
So by Lemma \ref{l454} the value of the 
Borel regulator on $\sum_x \sigma(\gamma_x)$ is equal to $2 \pi  \cdot\int_{\C {\Bbb P}^1} r_2(\sum_i (1-f_i) \wedge f_i \wedge g_i)$
and hence it is zero by our assumption. So Borel's theorem implies that the element is also zero.

A similar argument using Lemma \ref{1.8} shows that the homomorphism $h$ does not depend on the choice of the (finite) Galois extension of $\Q({\Bbb P}^1)$ containing the field $\Q(X)$.

{\bf 6. Explicit formulas for the reciprocity homomorphism $h$ and the Chow dilogarithm in the case of an elliptic curve}. 
Let $E$ be an elliptic curve. We want to calculate the integral 
$\int_{E(\C)} r_2(f_1 \wedge f_2 \wedge f_3)$. 
Let us suppose that $E$ is realized as a
plane curve. Then any rational function $f$ on $E$ can be written as a ratio
of products of {\it linear} homogeneous functions: 
$$
f = \frac{l_1\cdot ... \cdot l_k}{l_{k+1}\cdot ... \cdot l_{2k}}
$$
So it is enough to calculate the integral
$\int_{E(\C)} r_2(l_1/l_0 \wedge l_2/l_0 \wedge l_3/l_0)$  where $l_i$ are 
  linear functions in homogeneous coordinates. 
We will do this in a more general setting. 

\begin{figure}[ht]
\centerline{\epsfbox{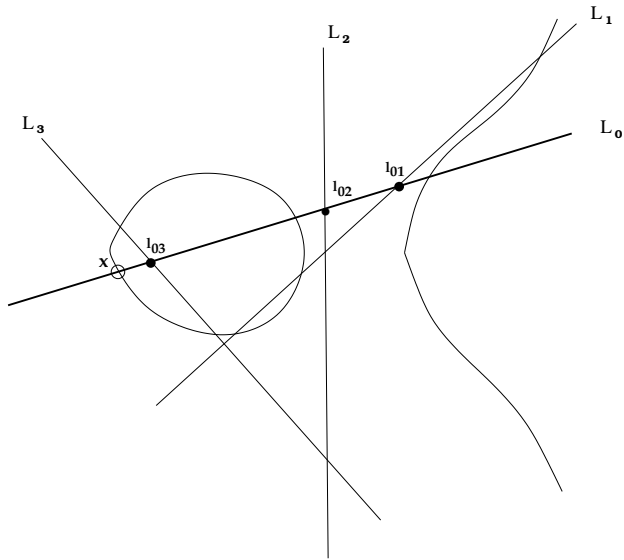}}
\caption{Defining $\{r(a,b,c,D)\}_2$ for a plane algebraic curve $X$.}
\label{kj}
\end{figure}

{\it Notations}. Let $X$ be a plane algebraic curve
 and $l_i$  linear functions in homogeneous coordinates. 
Denote by  $L_i$   the line $l_i =0$ in 
$\PP^2$. 
Let $D_i$ be the divisor $L_i \cap X$. 
Set $l_{ij}:= L_i \cap L_j$. 
For three points $a,b,c$ and a 
divisor $D =\sum n_i (x_i)$ on a 
line we will use  the following notation (see Figure \ref{kj}):
$$
\{r(a,b,c,D)\}_2:= \sum_i n_i \{r(a,b,c,x_i)\}_2
$$

\begin{theorem} \label{ecur}
 Let $E$ be an elliptic curve over an algebraically closed field $k$. Then there exists a   homomorphism of groups $h:\Lambda^3F^* \rightarrow
{\cal B}_2(k)$ such that for any linear homogeneous functions $l_0,...,l_3$ one has  
\begin{equation} \label{hrule}
h(l_1/l_0 \wedge l_2/l_0 \wedge l_3/l_0) = -\sum_{i=0}^3(-1)^i\{r(l_{i0},..., \widehat l_{ii}, ..., \l_{i3}, D_i)\}_2
\end{equation}
and which satisfies  all the properties of conjecture \ref{reh}. In particular, if $k =\C$ then
\begin{equation} \label{homot1}
{\cal P}_2(E; f_1\wedge f_2 \wedge f_3) =  {\cal L}_2 \Bigl( h(f_1\wedge
f_2 \wedge f_3 )\Bigr) 
\end{equation}
\end{theorem}

{\bf Proof}. Suppose we have four generic lines $L_0,L_1,L_2,L_3$ in $\PP^2$. Any two of them, say   $L_0$ and $L_1$, provide  a {\it canonical} rational function $(l_0/l_1)$ on $\PP^2$ with the divisor $L_0 - L_1$ normalized by the condition that its value at the point $l_{23}$  is equal to $1$.

\begin{lemma} \label{3.9}
a) On the line $L_3$ one has $(l_1/l_0) + (l_2/l_0) =1$.

b) $\frac{(l_1/l_0)}{(l_2/l_0)} = -(l_1/l_2)$
\end{lemma} 

{\bf Proof}.  Let $m$ be a point on the line $L_3$. Then
\begin{equation} \label{88}
(l_1/l_0)(m) = r(l_{03},l_{13},l_{23},m); \qquad (l_2/l_0)(m) = r(l_{03},l_{23},l_{13},m)
\end{equation}
This gives a).
It follows from this that if the point $m$ approaches the point $l_{03}$ then $\frac{(l_1/l_0)}{(l_2/l_0)}$ tends to $-1$. This implies b).

\begin{lemma} \label{lrez}
For any plane curve $X$ one has
$$
\sum_{x \in X} {\rm res}_x \Bigl((l_1/l_0) \wedge (l_2/l_0) \wedge ( l_3/l_0)\Bigr) =
- \delta_2\Bigl(\sum_{i=0}^3(-1)^i    \{r(l_{i0},..., \widehat l_{ii}, ..., \l_{i3}, D_i)\}_2\Bigr)
$$
\end{lemma}

{\bf Proof}. Let us compute first the   residues at the divisors $D_1,D_2,D_3$  using part a) of Lemma \ref{3.9}. For example 
the residue at $x \in D_1$  is equal to
$$
v_x((l_1/l_0))\cdot (l_2/l_0)(x) \wedge (l_3/l_0)(x) = v_x((l_1/l_0))\cdot (1-(l_3/l_0)(x) ) \wedge (l_3/l_0)(x) 
$$
$$
\stackrel{(\ref{88})}{=} v_x((l_1/l_0))\cdot \{r(l_{13}, l_{10}, l_{12}, x)\}_2 = v_x((l_1/l_0))\cdot \{r(l_{10}, l_{12}, l_{13}, x)\}_2
$$
It remains to compute the residues on the line $L_0$. According to part b) of  Lemma \ref{3.9} one has
$$
(l_1/l_0) \wedge (l_2/l_0) \wedge ( l_3/l_0) = - (l_0/l_1) \wedge (l_2/l_1) \wedge 
(l_3/l_1) 
$$
Using this we reduce the calculation to the previous case. 

 \begin{proposition} \label{hrez} Let $E$ be an elliptic curve over an algebraically closed field $k$. Then formula (\ref{hrule}) provides a well defined homomorphism of groups $h:\Lambda^3F^* \rightarrow
{\cal B}_2(k)$.
\end{proposition}

{\bf Proof}. Let $D:= \sum_i n_i(x_i)$ be the divisor of a rational function $f$ on $E$. 
To decompose it into a fraction of products of linear homogeneous functions $l_i$ we proceed as follows. Let $l_{x,y}$ (resp. $l_{x }$)  be a linear homogeneous equation of the line in $\PP^2$ through the points $x$ and $y$ on $E$  (resp. $x $ and $-x $).  The divisor of the function $l_{x,y}/l_{x}$ is $(x) +(y) - (x+y) - (0)$. 
If  $D = (x)+(y) + D_1$, we write $f =
l_{x,y}/l_{x+y}\cdot f'$, so $(f') = (0) + (x+y) + D_1$. 
 After a finite number of such steps we get the desired decomposition.

\begin{figure}[ht]
\centerline{\epsfbox{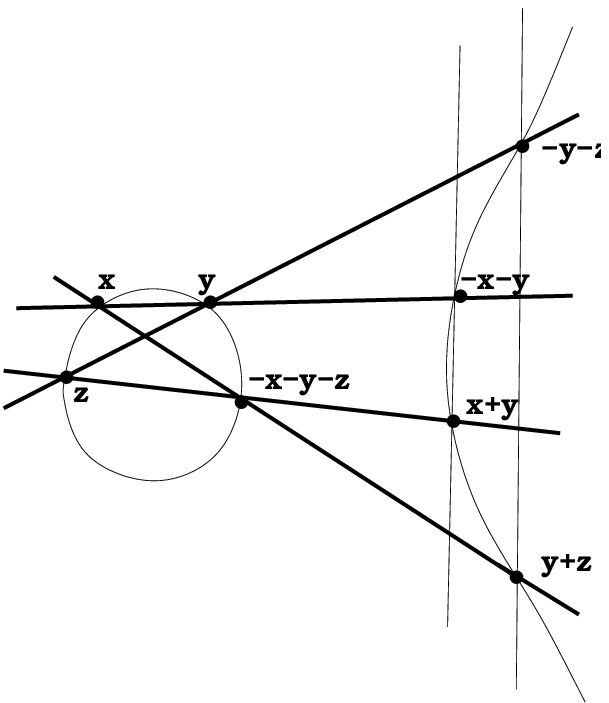}}
\caption{}
\label{kj3}
\end{figure}

There are the following relations
$$
\frac{l_{x,y}}{l_{x+y}} \cdot \frac{l_{x+y,z}}{l_{x+y+z}} \Big /\frac{l_{y,z}}{l_{y+z}} \cdot \frac{l_{x,y+z}}{l_{x+y+z}} \quad =  \quad \mbox{constant}
$$

One can prove that they generate all the relations between the functions $l_{x,y}/l_{x+y}$. 
 So $h$ is well defined if it annihilates the following expression:
$$
 F(x,y,z; l_0,l_2,l_3):= \Bigl(\frac{l_{x,y}}{l_{x+y}} \cdot \frac{l_{x+y,z}}{l_{x+y+z}} \Big / \frac{l_{y,z}}{l_{y+z}} \cdot \frac{l_{x,y+z}}{l_{x+y+z}} \Bigr) \wedge (l_2/l_0) \wedge ( l_3/l_0)  
$$
It follows from   Lemma \ref{lrez} that $\delta_2(F(x,y,z; l_0,l_2,l_3)) =0$. Thus according to the definition of the group ${\cal B}_2(k)$ it is enough to check that $h(F(x',y',z'; l_0,l_2,l_3))=0$ for a certain triple of points $(x',y',z')$. It is easy to see that $h(F(a,a,a; l_0,l_2,l_3))=0$ since then the first factor in $F$ is a constant. The proposition is proved.

\begin{proposition} \label{hrezza} Let $X$ be an algebraic curve in $\PP^2$ over $\C$ and $l_0,...,l_3$ linear homogeneous functions on $\C^3$. Then one has (using the notations defined above)
\begin{equation} \label{elfo}
\int_{X(\C)} r_2(l_1/l_0 \wedge l_2/l_0 \wedge l_3/l_0) = 2 \pi \cdot \sum_{i=0}^3(-1)^i
{\cal L}_2(r(l_{i0},..., \widehat l_{ii}, ..., \l_{i3}, D_i)
\end{equation}
\end{proposition}

{\bf Proof}. It follows from Lemma \ref{lrez} and Theorem \ref{55} that the differentials of the both sides coincide. So their difference is a constant. 
To show that this constant is zero 
we deform $X$ to a   union of lines in $\PP^2$. Using Proposition  \ref{p44} one sees that   formula (\ref{elfo}) is valid when $X$ is a 
 line in $\PP^2$.

\section{Appendix: on volumes of simplices in symmetric spaces}

  {\bf 1. Volumes of hyperbolic geodesic simplices as boundary integrals}.
 A point $y$ of  the $n$-dimensional hyperbolic space ${\cal H}_n$ defines a one 
dimensional space $M_y$ of volume forms  on the absolute $\partial{\cal H}_n$. It  
consists of the volume forms  invariant under the 
action of the isotropy group of $y$. We write them as follows. 
Let $x_0, ..., x_{n}$ be the coordinates in a vector space $V_{n+1}$ of dimension $n+1$, 
and 
$
Q(x):= x_0^2 + ...  +  x_{n-1}^2 - x_{n}^2 
$. 
Then ${\cal H}_n$ can be realized as the  
projectivization of the cone $Q(x) < 0$, and its boundary is the projectivization of the cone
 $Q(x) = 0$. Choose a point $y \in {\cal H}_n$. Lifting $y$ to a vector $y' \in V_{n+1}$ 
we have the following 
volume form on the boundary $\partial{\cal H}_n$:
$$
\frac{\delta(Q(x)) \sigma_{n+1}(x,dx)}{(x,y')^{n-1}}
$$
If  $y$ belongs to the boundary $\partial{\cal H}_n$, this formula provides a 
space $M_y$ of singular volume forms   on the absolute; they are invariant under the 
  isotropy group of  $y$.

Let us choose for any point $y$ such a volume form $\mu_y$. 
For two points $x,y$ the 
ratio $  \mu_x/\mu_y$ is a nonzero  function on the absolute.

Let $I(y_0,...,y_{n })$ be the geodesic simplex with vertices at 
$y_0,...,y_{n }$ where the points $y_i$ could be on the absolute. 
Denote by ${\rm vol}I(y_0,...,y_{n})$ the volume of this simplex with respect to the invariant volume form  in ${\cal H}_n$ normalized by the following condition: if we realize the hyperbolic space as the interior of the unit ball $y_1^2 + ... + y_n^2 \leq 1$ then the volume form restricted to the tangent space at the origin $(0, ..., 0)$ is $dy_1 \wedge ... \wedge dy_n$.  
 
\begin{theorem} \label{2112}
For any hyperbolic geodesic simplex 
$I(y_0,...,y_{n})$ one has 
\begin{equation} \label{1221} 
\frac{(n-1)^n{\rm vol}(S^{n-1})}{n} \cdot {\rm vol}I(y_0,...,y_{n}) =   
\int_{\partial {\cal H}_n}
\log | \frac{\mu_{y_1}}{\mu_{y_0}}| d\log|\frac{\mu_{y_2}}{\mu_{y_0}} | \wedge ... \wedge d\log| \frac{\mu_{y_{n}}}{\mu_{y_0}} | 
\end{equation} 
 \end{theorem}

 Let    $\varphi(y_0, ..., y_{n})$ be the function defined by the right hand side of (\ref{1221}). Thanks to property 2) of Proposition \ref{1.2.} it does not depend on the choice of invariant volume forms $\mu_y$.

\begin{proposition} \label{1.21.1}
The function $\varphi(y_0,...,y_{n})$   has  the following properties: It is

1)  A smooth function on the vertices $y_i$.

2) Equal to zero if three of the vertices  belong to the same geodesic. 

3) Additive with respect to cutting of a simplex, i.e. if $y_0,...,y_{n}$ are points 
such  that $y_0,y_1,y_2 $ are on the same geodesic, $y_1$ between $y_0$ and $y_2$, then
$$
\varphi(y_0,y_2,...,y_{n+1}) = \varphi(y_0,y_1,y_3, ...,y_{n+1})  + \varphi(y_1,y_2,...,y_{n+1}) 
$$

4) Invariant under the action of the symmetry group $SO(n,1)$.
\end{proposition}

{\bf Proof}. 1) This is clear. 

2) Let us  realize the  hyperbolic space as the interior part of  
 the unit ball in $\R ^{n}$. Consider the geodesic $l$ passing through  the center of the ball  in the vertical  direction. The subgroup $SO(n-1) \subset SO(n,1)$ preserves pointwise this geodesic. So for any point $y$ on the geodesic the invariant volume form $\mu_y$  is invariant under the action of the group $SO(n-1)$. 
 The quotient of $\partial {\cal H}_n$  under the action of $SO(n-1)$ is given by the projection $p: \partial {\cal H}_n \longrightarrow \R $ onto the vertical axis. 
Take three points $y_1,y_2,y_3$ on the geodesic.  Then $\frac{\mu_{y_2}}{\mu_{y_1}}$ and $\frac{\mu_{y_3}}{\mu_{y_1}}$ are  lifted  from the line $\R $. Therefore $d\log|\frac{\mu_{y_2}}{\mu_{y_1}}| \wedge d\log|\frac{\mu_{y_3}}{\mu_{y_1}}|$ = 0.
 So for a degenerate simplex $(y_0,...,y_n)$ one has 
$$ 
\log|\frac{\mu_{y_0}}{\mu_{y_1}}| d\log|\frac{\mu_{y_2}}{\mu_{y_1}}| 
\wedge d\log|\frac{\mu_{y_3}}{\mu_{y_1}}|\wedge ... \wedge d\log|
\frac{\mu_{y_n}}{\mu_{y_1}}| =0
$$
It remains to mention the skewsymmetry  of the integral (\ref{1221}). 
The property 2) is proved.

3) Follows from 2) and the additivity property from  Proposition  
\ref{1.2.}.

4) This is clear from 4) of Proposition  
\ref{1.2.}. The  proposition is proved.

The leading term of the Taylor expansion of the function $\varphi(y_0, y_1, ..., y_n)$ 
when $y_0$ is fixed and $y_1, ..., y_n$ are near $y_0$ provides an exterior $n$-form in 
$T_{y_0}{\cal H}_n$ denoted $\varphi_{y_0}(Y_1, ..., Y_n)$,  $Y_i \in T_{y_0}{\cal H}_n$. 
Let us compare it with the volume form 
${\rm Vol}_{y_0}(Y_1, ..., Y_n)$ in $T_{y_0}{\cal H}_n$ normalized as 
before Theorem 7.1. 

\begin{lemma} \label{1.21.1}
$\varphi_{y_0}(Y_1, ..., Y_n)  = \frac{(n-1)^n}{n}{\rm vol}(S^{n-1}) \cdot
{\rm Vol}_{y_0}(Y_1, ..., Y_n)$. 
\end{lemma}

{\bf Proof}.   Below we abuse notations by writing $y$  for $y'$.  
One has $\mu_{y_1}/\mu_{y_2} = \frac{(y_2,x)^{n-1}}{(y_1,x)^{n-1}}$. So 
$$
\varphi(y_0, ..., y_n) = (n-1)^{n}\int_{\partial {\cal H}_n} \log\vert 
\frac{(y_1,x)}{(y_0,x)}\vert  d  \log\vert 
\frac{(y_2,x)}{(y_0,x)}\vert \wedge ... \wedge d  \log\vert 
\frac{(y_n,x)}{(y_0,x)}\vert 
$$
Thus 
$$
\varphi_{y_0}(Y_1, ..., Y_n) = (n-1)^{n}\int_{\partial {\cal H}_n} (Y_1,x) d( Y_2, x) \wedge ... \wedge 
d(Y_n, x)
$$
To do the computation of this integral we may suppose that $y_0 = (0, ..., 0, 1)$, $Y_i = 
\frac{\partial}{\partial y_i}$, so $(Y_i, x) = x_i$. Then the last integral equals to 
$$
(n-1)^{n}\int_{S^{n-1}} x_1 dx_2 \wedge ... \wedge dx_n =  \frac{(n-1)^{n}}{n}{\rm vol}(S^{n-1})
$$
where $S^{n-1}$ is the sphere $x_1^2 + ... + x_n^2 =1$. 
The lemma follows.

{\bf Proof of Theorem \ref{2112}}. Let us suppose first that the points $x_i$ 
are inside of the hyperbolic space. 

The  function   $\varphi(x_0,...,x_{n})$  defines an 
$n$-density $\widetilde \varphi$ on
  $ {\cal H}_n$. 
Namely,   to define the integral $\widetilde \varphi$  over a simplex $M$ one has to 
subdivide it into small simplices and take the sum of the functions 
$\varphi$ corresponding to their vertices.   
When the simplices are getting smaller the limit exists and is by definition  
$\int_M \widetilde \varphi$. Here we used  the properties 1)  - 3).
More precisely 1) and 2)  imply that $\varphi$ defines an additive volume form on 
${\cal H}_n$, and 3) (together with 1)) guarantee that this volume form is $\sigma$-additive.

The skewsymmetry property implies that  $\widetilde \varphi$ is actually a differential $n$-form. It is invariant under the action of 
the group $SO(n,1)$. Therefore it is proportional to the standard volume form. 

Now suppose that the vertices $x_i$ can be  on the absolute. Then 
it is easy to see that the corresponding integral 
(\ref{i2})  is still convergent. Moreover, 
if the vertices of the geodesic simplex are in general position 
then it is a continuous function of the vertices. This  implies that  the volume of an 
ideal geodesic simplex is finite (which is, of course, an elementary fact) and 
 coincides with the corresponding integral (\ref{1221}).

  Completely similar results are valid for the complex $n$-dimensional   hyperbolic 
space ${{\cal H}_n}^{\C}:= \{|z_1| + ... +|z_n|^2 < 1\}$, $z_i \in \C$, and the  
quaternionic hyperbolic space ${{\cal H}_n}^{\HH}:= \{|q_1| + ... +|q_n|^2 < 1\}$ ($q_i$ are quaternions). Indeed, a point $x$ in each of these spaces defines an invariant volume form $\mu_x$ on the boundary.

{\bf 2. Calculation of the volume of a three dimensional ideal geodesic simplex}. 
If $n=3$ the absolute can be identified with $\C \PP^1$, 
and for the ideal simplex with vertices at 
the points $\infty, 0,1,a$ on the absolute we get 
$$
{\rm vol} (I(\infty, 0,1,a)) = 3 c_3 \cdot 
\int_{\C \PP^1} (\log|z| d\log|1-z| - \log|1-z| d\log|z|) \wedge d\log|z-a| =
$$
$$
3 c_3 \cdot \int_{\C \PP^1} (\log|z| d\arg (1-z) - \log|1-z| d\arg (z))\wedge d\arg (z-a)
$$
because $d \log (z-1) \wedge d \log (z-a) = 0$. 
Here $\log f = \log |f| + i \arg(f)$. 
Using 
$$
d {\cal L}_2(z) =  \log|z| d\arg (1-z)  - \log|1-z| d\arg (z)
$$ 
we rewrite the last integral as 
\begin{equation} \label{1331} 
3 c_3 \cdot \int_{\C \PP^1} d {\cal L}_2(z) \wedge d\arg (z-a)
\end{equation} 
Computing 
the differential in the sense of distributions we get
$$
d ({\cal L}_2(z) d\arg (z-a)) = 2\pi \cdot  {\cal L}_2(z) \delta(z-a)dx d y + 
d {\cal L}_2(z)\wedge d\arg (z-a)
$$
So the integral of the right hand side over $\C \PP^1$ is zero, i.e. the integral (\ref{1331}) is equal to 
$- 6 \pi c_3 \cdot {\cal L}_2(a)$. ( $c_3 = -1/6\pi$).

{\bf 3. Volumes of geodesic simplices in $SL_n(\C)/SU(n)$}. Recall 
the invariant differential $(2n-1)$--form $\omega_{D_n}$ in ${\Bbb H}_n$.  

{\bf Question}. Is it true that  
\begin{equation} \label{1221?q} 
{\rm vol}_{\omega_{D_n}}(I(x_0,...,x_{2n-1})) = \mbox{\rm constant} 
\times \psi_n(x_0,...,x_{2n-1})?
\end{equation}

One can show following the lines of s. 7.1 that the positive answer 
to this question is equivalent to the following statement: if $x_0, x_1, x_2$ 
are on the same geodesic then $\psi_n (x_0,...,x_{2n-1}) =0$.

{\bf 4. Another approach to Grassmannian polylogarithms}. 
The following construction was suggested to the author 
during the Fall of 1989, independently, 
 by M. Kontsevich and by J. Nekovar. 
A hyperplane $h$ in an $n$--dimensional 
 complex vector space $V$ determines an arrow in the space of degenerate 
non-negative definite 
hermitian forms in $V$ consisting of the forms with the kernel $h$. 
Let $h_1, ..., h_{2n}$ 
be hyperplanes in $V$. Let $C(h_1, ..., h_{2n})$ be projectivization of 
the convex hull of the arrays corresponding to these
hyperplanes. It is a simplex in $\overline {\Bbb H}_n$. The idea is to 
integrate  
the form $\omega_{D_n}$ over this simplex. 
If $n=2$ this construction provides an ideal geodesic simplex 
in the Cayley realization of the hyperbolic space, given by  
the interior part of a ball in $\R \PP^3$.  However 
the convergence of this integral for $n>2$ 
has not been established yet, 
although it does not seem to be a very difficult problem.  
If the integral is convergent, we get a function on configurations 
of $2n$ hyperplanes in $\CP^{n-1}$. 
It would be very interesting to investigate this construction further and 
compare it with our construction of the Grassmannian polylogarithms. 

{\bf 5. A  $(2n-1)$-cocycle  of $GL(\C)$}.
Consider an infinite dimensional $\C$-vector space with a given basis 
$e_1,...,e_m,...$. The group $GL(\C)$ is the 
group of automorphisms of this space moving only finite number of basis
vectors.  

Let us describe the restriction of the cocycle to the subgroup
$GL_{n+m}(\C)$ acting on the subspace generated by first $n+m$
vectors. 
Take $2n$ elements $g_1,..., g_{2n}$ of this group. Consider the corresponding
$2n$ $(m+1)$-tuples of vectors:
$$
g_1(e_{n},...,e_{n+m}), \quad  ... , \quad g_{2n}(e_{n},...,e_{n+m})
$$
The set of all  $(m+1)$--tuples of vectors form a vector space. 
Let $\Delta_{2n-1}$ be the standard simplex $\sum^{2n}_{j=1} \lambda_j =
1$.
Consider the set 
$
C_n^{m} \subset \Delta_{2n-1}
$
consisting of  $(\lambda_1,...,\lambda_{2n})$ such that 
$$
\sum_{i=1}^{2n} \lambda_i \cdot g_i(e_{n},...,e_{n+m})
$$
is an $(m+1)$-tuple of vectors {\it not} in generic position. It is a cycle of
codimension $n$. The Chow polylogarithm function evaluated 
on it provides the desired (measurable) cocycle. 

The cocycle property follows from the functional equation for the Chow
polylogarithm and the following general fact. 
The set of those $(\lambda_1,...,\lambda_{2n+1})$ such that 
$
\sum_i \lambda_i \cdot g_i(e_{n},...,e_{n+m})
$
is an $(m+1)$-tuple of vectors {\it not} in generic position is also of
codimension $n$.

The construction is consistent with the restriction to $GL_n$ just by definition.

{\bf Problem}. Show that the cohomology class of this cocycle is nontrivial and proportional to the Borel class.

\bigskip

Department of Mathematics, Brown University, Providence, RI 02912, USA. 

e-mail sasha@math.brown.edu

\end{document}